\documentclass[12pt]{report}
\usepackage{latexsym}
\usepackage{amsfonts}
\usepackage{amssymb}
\usepackage{color}

\begin{document}

\title{Connectivity of the Coset Poset}

\author{Daniel A. Ramras\\
        \vspace{.15in}\\
	Department of Mathematics\\
        Cornell University\\
	May, 2002\\
	(Revised September, 2002)\\
       }

\date{}

\maketitle

\tableofcontents

\newcommand{\normal}{\vartriangleleft}
\newcommand{\nnormal}{\ntriangleleft}
\newcommand{\R}{\mathbb{R}}
\newcommand{\Z}{\mathbb{Z}}
\newcommand{\F}{\mathbb{F}}
\newcommand{\A}{\mathcal{A}}
\newcommand{\bbC}{\mathbb{C}}
\newcommand{\abs}[1]{\left| #1\right|}
\newcommand{\ord}[1]{\Delta \left( #1 \right)}
\newcommand{\leqs}{\leqslant}
\newcommand{\geqs}{\geqslant}
\newcommand{\heq}{\simeq}
\newcommand{\iso}{\simeq}
\newcommand{\implies}{\Longrightarrow}
\newcommand{\maps}{\longrightarrow}
\newcommand{\meet}{\wedge}
\newcommand{\join}{\vee}
\newcommand{\homeo}{\cong}
\newcommand{\isom}{\cong}
\newcommand{\cross}{\times}
\newcommand{\mC}{\mathcal{C}}
\newcommand{\C}[1]{\mathcal{C}(#1)}
\newcommand{\M}[1]{\mathcal{M}(#1)}
\newcommand{\supp}[1]{\rm supp\left( #1 \right)}
\newcommand{\gen}[1]{\langle #1 \rangle}
\newcommand{\prj}{\F_p \cup \{\infty\}}
\newcommand{\gl}{GL_2 (\F_p)}
\newcommand{\spl}{SL_2 (\F_p)}
\newcommand{\psl}{PSL_2 (\F_p)}
\newcommand{\psls}{PSL_2 (\F_7)}
\newcommand{\injects}{\hookrightarrow}

\newtheorem{fact}{Fact}[section]
\newtheorem{lemma}[fact]{Lemma}
\newtheorem{theorem}[fact]{Theorem}
\newtheorem{definition}[fact]{Definition}
\newtheorem{corollary}[fact]{Corollary}
\newtheorem{proposition}[fact]{Proposition}
\newtheorem{claim}[fact]{Claim}
\newtheorem{remark}[fact]{Remark}
\newtheorem{question}[fact]{Question}

\vspace{4in}
\noindent {\bf Acknowledgements:}  

I would like to thank Professor Kenneth S. Brown for his guidance throughout the course of
this project.  Professor Brown served as my mentor for the two years I spent
on this research, and his assistance was invaluable.  In particular, Professor
Brown suggested studying the coset poset and introduced me to most of the techniques employed
in this thesis.  Two of my professors also deserve special thanks: Professor Ravi Ramakrishna, 
for teaching me algebra and providing guidance throughout my 
undergraduate career, and Professor Marshall Cohen, who taught me algebraic topology
from the simplicial viewpoint.  Finally, I would like to thank the Cornell Presidential
Research Scholars program, which provided me with funding for the first year of this project.

\chapter{Introduction}$\label{intro}$

At the intersection of topology, algebra and combinatorics lies the study
of simplicial complexes arising from finite groups.
Over the last twenty-five years, there has been a great deal of research
in this area, largely
stimulated by Quillen's influential paper~\cite{Quillen-p-subgroups}
on the $p$-group complex.  A number of other references in the area are listed in the
bibliography~\cite{Aschbacher-simple-connectivity,Aschbacher-Quillen-conjecture,
Brown-coset-poset,Kratzer-subgroups,Welker-solvable-p-subgroup,
Shareshian-shell-solv,Thevenaz-perm-reps,Webb-Quillen-complex}
This thesis represents a contribution to the theory, focusing in particular 
on a new complex (the coset poset) recently introduced by 
Brown~\cite{Brown-coset-poset}.  

The coset poset of a group finite $G$, denoted $\C{G}$, is the partially ordered set 
containing all cosets of all proper subgroups of $G$ (ordered by inclusion).  One
forms a simplicial complex on the vertex set $\C{G}$ by declaring each chain 
$$x_1 H_1<\cdots < x_n H_n$$ 
in $\C{G}$
to be a simplex.
This complex arose from Brown's analysis
of $P(G,s)$, the probabilistic zeta function
of a finite group.  This function gives the probability that $s$ elements
of $G$, chosen at random, generate $G$.  The theory of M\"{o}bius inversion allows
one to write $P(G,s)$ in the form 
$$P(G,s) = \sum_{H \leqs G} \frac{\mu (H,G)}{(G:H)^s},$$
where $\mu$ is the M\"{o}bius function of the subgroup lattice of $G$ (defined in 
Section~\ref{mu(psls)-section}).
Written as such, we see that $P(G,s)$ is defined for any complex number $s\in \bbC$.
The starting point of~\cite{Brown-coset-poset} is the observation that $P(G,-1)$ is in fact the 
Euler characteristic of $\C{G}$.  

This result motivates the study of the homotopy type of $\C{G}$.  In this thesis,
we will mainly be concerned with connectivity.
In particular we will investigate the question of simple connectivity, 
raised in Brown~\cite[Question 4]{Brown-coset-poset}.
Although we are not able to give a characterization of groups with simply connected
coset posets, we do present several results in either direction.

The structure of this thesis is as follows.  In Chapter 2, we develop the
basic topological theory of simplicial complexes and partially ordered sets.  This theory
is developed from scratch, assuming only standard topological prerequisites.  After
describing the geometric realization functors and proving some of their basic properties,
we proceed to prove a variety of results on the homotopy type of posets and 
simplicial complexes.  Most of this background material is taken from 
Bj\"{o}rner~\cite{Bjorner-top-methods}.

In Chapter 3, we introduce the coset poset and prove a variety of new results 
on its connectivity.  Some of these results extend to a more general class of posets
we call $\emph{atomized posets}$.  These are essentially posets in which any set of minimal 
elements ``generates" a well-defined element of the poset (or generates the whole poset).
In addition to $\C{G}$ this class includes, for example, 
the poset of proper, non-trivial subgroups of a finite group.  In these cases, the notion of 
generation in the poset coincides with the group-theoretic definition.

The three central methods in Chapter 3 are the introduction of the minimal cover of
an atomized poset, the use of Mayer-Vietoris sequences to analyze the homology
of $\C{G}$, and the technique (introduced in Brown~\cite{Brown-coset-poset}) of finding a
large subposet that can be described as a join.  
The minimal cover of $\C{G}$ is the simplicial complex with vertex set 
$G$ and with a simplex for each subset of a (proper) coset of $G$.  Denoting the minimal cover
by $\M{G}$, we show that $\M{G}\heq \C{G}$ (Lemma~\ref{min-cover}).  
For many purposes, $\M{G}$ is easier to analyze,
and many of our results are based on the structure of this complex.  For example,
we prove that any group not generated by $k$ elements has a $(k-1)$-connected coset poset
(Theorem~\ref{k-connected-atom}).  
If $G$ is a finite 2-generator group this theorem says only that
$\C{G}$ is connected, but with additional hypotheses we show that $\C{G}$ is simply connected 
as well (Theorem~\ref{2-gen-presentation}).

Using Mayer-Vietoris sequences, we show that certain classes of finite groups have 
non-simply connected coset posets (Theorems~\ref{M-V1} and~\ref{M-V2}).  
In addition, we give lower bounds
on the rank of the first homology of $\C{G}$ for these groups.  These results lead to a purely
group-theoretic result on the number of split factors in a chief series for certain groups
(Corollaries~\ref{chief-series1} and~\ref{chief-series2}).

In the last section of Chapter 3, we consider the coset poset of a group $G$ written
as an extension $N\injects G \maps G/N$.  
As a corollary of Brown's analysis of products~\cite[Lemma 5]{Brown-coset-poset},
we deduce a simple characterization of those direct products with simply connected coset posets.
In addition, we generalize~\cite[Lemma 5]{Brown-coset-poset} to the case of a semi-direct
product, and note a simple result about general extensions.

In the final chapter, we compute a number of interesting examples.  First,
we describe a third complex homotopy equivalent to $\C{G}$.  This complex immediately 
gives the homotopy-type of $\C{Q_8}$, where $Q_8$ denotes the Quaternion group
of eight elements. The computation provides a counterexample to the converse of an 
earlier result.  (Note that even for a small group such as $Q_8$, the coset poset 
is too large to admit direct analysis: $\C{Q_8}$ is a two-dimensional complex with 
eighteen vertices, forty-four edges and twenty-four triangles.)

Next, we turn to finite simple groups.  Brown has given a complete description
of the homotopy-type of $\C{G}$ for any finite solvable group $G$, and our results about
extensions suggest that the existence of a non-trivial normal subgroup greatly simplifies 
the coset poset.
This motivates the study of simple groups, at the other end of the spectrum.  
First we calculate the homotopy-type of $\C{A_5}$, giving two simpler proofs of an 
unpublished result of Shareshian (cited in~\cite[p. 1009]{Brown-coset-poset}).
The first proof makes use of a theorem from Chapter 3, and the second utilizes
2-transitive actions.

We also consider the coset poset of $PSL_2(\F_7)$.  We begin by describing in detail
the subgroup lattice of the finite simple groups $PSL_2(\F_p)$, following the analysis
given in Burnside~\cite{Burnside-theory-of-groups}.  
We then show that $\C{PSL_2(\F_7)}$ is simply connected and make a few other
observations about its homotopy-type.

In the final chapter we suggest several new directions for research and formulate some
conjectures on the homotopy type of $\C{G}$.


\chapter{Background Material}$\label{background}$

In this chapter we will introduce the topological definitions and 
machinery that will be used through this thesis.
We will deal mainly with simplicial complexes, and will focus in particular on
simplicial complexes arising from finite partially ordered sets (posets).  

Throughout this thesis, the word map will mean a morphism in the appropriate category.
When the category is not clear from context, we will be more explicit.

\section{Simplicial Complexes}
\begin{definition}  An abstract simplicial complex is a pair 
$S = (V,\Delta)$ where $V$ is a finite set (called the vertices of $S$) 
and $\Delta$ is a collection of subsets of $V$ (called the simplices of $S$)
satisfying the following properties:

\begin{enumerate}
	\item If $\sigma_1\subset \sigma_2$ and $\sigma_2\in \Delta$, then 
                $\sigma_1\in \Delta$,
        \item For all $v\in V$, $\{v\}\in \Delta$.
\end{enumerate}

In other words, $\Delta$ is closed under inclusion and contains all singleton-sets.

For any simplex $\sigma \in \Delta$, let $\dim{\sigma} = \abs{\sigma} - 1$.  If 
$\dim{\sigma} = n$ we call $\sigma$ an $n$-simplex.
\end{definition}

From here on, \emph{simplicial complex} will mean abstract simplicial complex. 
If $S = (V,\Delta)$ is a simplicial complex, then by condition 2, 
$V = \bigcup_{\sigma \in \Delta} \sigma$ and we may suppress $V$ from the notation and 
denote the complex simply by $\Delta$.  When necessary, we will denote the vertex
set of $\Delta$ by $V(\Delta)$.

If $\Delta$ is a simplicial complex and $\Gamma \subset \Delta$ is also a simplicial
complex, then we call $\Gamma$ a \emph{subcomplex} of $\Delta$.
For $k = 0,1,\ldots$ we define $\Delta^k$, the $k$-skeleton of $\Delta$, to be 
the subcomplex containing all simplices of dimension at most $k$.

To each simplicial complex we will now associate a topological space, called the
\emph{geometric realization} of the complex.  

For any finite set $X$, let $\R^X$ be the $\R$-vector space with basis $X$  
(so $\R^X\isom \R^{|X|}$).

If $t = \sum_{i=1}^n t_i x_i \in \R^X$ ($x_i\in X$), then the \emph{support} of $t$ is the set
$\supp{t} = \{x_i: t_i\neq 0\}$.  We will sometimes denote the point $t$ simply by
$(t_1,\ldots, t_n)$.

\begin{definition}  Let $\Delta$ be a simplicial complex with vertex set 
$V$.
Then $\abs{\Delta}$, the \emph{geometric realization} of $\Delta$, 
is the subspace of $\R^{V}$ given by
$$\abs{\Delta} = \{t = (t_1,\ldots,t_n)\in \R^V : t_i\geq 0\,\, \forall i,\,\,
				      		  \sum_{i=1}^{n} t_i = 1 
						  {\rm \,\, and \,\,} \supp{t}\in \Delta
       	         \}.$$

The geometric realization of a non-empty simplex 
$\sigma \in \Delta$ is the subspace
$\abs{\sigma} = \{\sum_{v\in V} t_v v \in |\Delta| : t_i = 0$ if $v\notin \sigma\}$.
\end{definition}

The geometric realization of a simplex thus consists of all weighted
averages of the vertices of that simplex, and the geometric realization of a complex is 
the union of the geometric realizations of its simplices.
Thus the geometric realizations of 0-, 1-, 2-, and 3-simplicexes are points,
lines, triangles and tetrahedrons respectively.

Although the same notation is used for both cardinality and geometric realization,
in context it will always be clear what is meant, or else we will be more explicit.
Also, we will sometimes refer simply to $\Delta$, and when references are made the 
topology of $\Delta$ this will always mean the topology of $|\Delta|$.

An important class of maps between simplicial complexes are the 
\emph{simplicial maps}, maps that respect the combinatorial sturcture 
of the complexes in question.  

\begin{definition} Let $\Delta$ and $\Gamma$ be simplicial complexes.  
We call a map $f:V(\Delta) \to V(\Gamma)$ \emph{simplicial} if 
for each simplex $\sigma \in \Delta$, $f(\sigma)\in \Gamma$.  (Note that $f$ need not be 
bijective; $\dim{f(\sigma)} < \dim{\sigma}$ is allowed.)

Any simplicial map $f:V(\Delta) \to V(\Gamma)$ induces a map $|f|$ on geometric realizations, 
called the \emph{geometric realization} of $f$.  For 
$t = \sum_{v\in V(\Delta)} t_v v \in |\Delta|$, we define
$|f|(t) = \sum_{v\in V(\Delta)} t_v f(v)$.
\end{definition}

A simplicial map $f:V(\Delta) \to V(\Gamma)$ induces a map $f':\Delta \to \Gamma$, and
we will usually introduce a simplicial map in the latter form.  So a map $f':\Delta \to \Gamma$
is simplicial iff it is induced by some simplicial map $f:V(\Delta) \to V(\Gamma)$.

One important fact about geometric realizations is that the realization of a simplicial
map is always continuous.  Together with the (easy) fact that $|f\circ g| = |f|\circ |g|$
(i.e. geometric realizations commute with compositions) this will mean that $|\cdot|$ is
a \emph{functor} from the category of simplicial complexes and simplicial maps to the 
category of topological spaces and continous maps.  In order to prove that the geometric
realization of a simplicial map is continuous, we will need the following basic result about
the topology of a simplicial complex.

\begin{claim} Let $\Delta$ be a simplicial complex.  Then a set $U\subset |\Delta|$
is open iff $U\cap |\sigma|$ is open in $|\sigma|$ for each $\sigma\in \Delta$.  Also, if X is any
topological space, a map $f:|\Delta| \to X$ is continuous iff its restriction to each
simplex $|\sigma|$ is continuous.
$\label{weak}$
\end{claim}

Now we show that the geometric realization of a simplicial map is continuous.  

\begin{claim}
Let $f:\Delta \to \Gamma$ be a simplicial map.  Then 
$|f|:|\Delta|\to |\Gamma|$ is continuous.
\end{claim}
{\bf Proof.}  By Claim~\ref{weak} we simply need to show that the restriction
of $|f|$
to a simplex $|\sigma|\subset |\Delta|$ is continuous.  Setting $f(\sigma) = \tau$, it
suffices to check that the map $\hat{f}:|\sigma|\to |\tau|$, 
$\hat{f} (t) = |f|(t)$, is continuous. 

Let $V\subset |\tau|$ be open, let $U = \hat{f}^{-1}(V)$ and let $\dim{\sigma} = n$.  
For any $x\in U$, we may choose $\epsilon > 0$ s.t. all points of the form 
$t = f(x) + \sum_{v\in \tau} \epsilon_v v$, with 
$\abs{\epsilon_v} < \epsilon$ for each $v$ and $t \in \tau$, are in $V$ 
(this is possible because $V$ is open in $\tau$).  
Then any point of the form $x' = x + \sum_{u\in \sigma} \frac{\epsilon_u}{n + 1} u$, with 
$\abs{\epsilon_u} < \epsilon$ and with $x'\in \abs{\sigma}$, is in $U$ since 
$\hat{f}(x') = f(x) + \sum_{u\in \sigma} \frac{\epsilon_u}{n + 1} f(u)$, and the coefficient
on any particular $v\in \tau$ has absolute value at most 
$(n + 1) \frac{\epsilon}{n + 1} = \epsilon$.
$\hfill \Box$

\vspace{.15in}
Our next goal will be to give a presentation for the fundamental group of an
arbitrary simplicial complex in terms of the combinatorial structure of its
$2$-skeleton.  In order to do this, we will first need to show that for any
simplicial complex $\Delta$ we have $\pi_1(\Delta) \homeo \pi_1(\Delta^2)$.  We will
in fact prove this latter fact for CW-complexes,
and our proof will be a simple application of the Cellular Approximation Theorem.
For standard definitions and facts about CW-complexes, we refer to 
Hatcher~\cite{Hatcher-alg-top}.

We denote the $n$-skeleton of a CW-complex $X$ by $X^n$, and we call a map $f:X\to Y$,
with $X$ and $Y$ CW-complexes, \emph{cellular} if $f(X^n) \subset Y^n$ for each $n$.

If $\Delta$ is a simplicial complex, then $\abs{\Delta}$ has an obvious CW structure,
and the definitions of the skeleta coincide.  In addition, the geometric realization
of a simplicial map is cellular.

\begin{theorem}[Cellular Approximation Theorem] $\label{cell-approx}$
Let $f:X\to Y$ be a map between two CW-complexes.  If $f$ is cellular on the $n$-skeleton
of $X$, ($n = 0, 1,\ldots$) then there exists a cellular map $f':X\to Y$ s.t.
$f\heq f'$ (rel $X^n$).
\end{theorem}

\begin{corollary}$\label{2-skeleton}$
If $\Delta$ is a (connected) simplicial complex, 
then $\pi_1(\Delta)\homeo \pi_1(\Delta^2)$.
\end{corollary}
{\bf Proof.}  Fix a basepoint $v\in V(\Delta)$.  
Consider the map $f: \pi_1(\Delta^2) \to \pi_1(\Delta)$ induced by 
the inclusion $\Delta^2 \hookrightarrow \Delta$.  We must show that this is a bijection.
First, note that each homotopy class of loops has a cellular representative, 
by Theorem~\ref{cell-approx}.
This immediately means that $f$ is surjective, since any cellular loop is 
contained in $\abs{\Delta^1}\subset \abs{\Delta^2}$.

Now, in order to prove injectivity we must show that if 
$\alpha$, $\beta : I\to \abs{\Delta^2}$ are two cellular loops based at $v$, then 
$\alpha\heq \beta$ in $\abs{\Delta}$ implies $\alpha \heq \beta$ in $\abs{\Delta^2}$.  But
this is immediate from Theorem~\ref{cell-approx}, since the homotopy in $\abs{\Delta}$ 
may be pushed down into $\abs{\Delta^2}$ keeping its restrictions $\alpha$ and $\beta$ fixed.
$\hfill \Box$

\vspace{.15in}
Next, we recall that given any presentation $\mathcal{P} = \langle S|R \rangle$ where $S$ 
is a set of generators and $R$ is a set of relators (words in the symbols of $S$ and their 
formal inverses) there is a CW-complex $K_{\mathcal{P}}$ whose fundamental group
has $\mathcal{P}$ as a presentation.  We will describe the construction in the case where
$S$ and $R$ are each finite.

This complex is formed by starting with a wedge
of circles, one for each element $s\in S$ and then a attaching 2-cells for each 
word $r\in R$.  Formally, the 1-skeleton of $K_{\mathcal{P}}$ is the disjoint union
of $\abs{S}$ copies of $I$, (we denote these copies by $I_s$, $s\in S$)
mod the relation $0_s \sim 1_s \sim 0_t \sim 1_t$ for all $s,t\in S$.  We let $q$ denote
the quotient map.

In order to attach the 2-cells, we must first choose orientations on the 1-cells (the images
of the copies of $I$ under $q$).  Given $s\in S$, we choose a map $f_s : I\to  I_s$,
where either $f_s (t) = t$ or $f_s (t) = 1-t$.  Composing $f_s$ with $q$ yields the orientation
on the 1-cell $q(I_s)$.

If $r = s_0^{\epsilon_0}\cdots s_n^{\epsilon_n}$ 
(where $s_i\in S$ and $\epsilon_i = \pm 1$) then the 2-cell corresponding to $r$
is attached via the map $\phi_r: I\to \vee_{s\in S}\,\, q(I_s)$, where 
$\phi_r|_{[\frac{i}{n+1}, \frac{i+1}{n+1}]}$ is just a linear shift of the map $f_l$.
(To be more specific, say $0\leq a < b\leq 1$ and let $L_{[a,b]} : [a,b]\to I$ be 
the homeomorphism given by $L_{[a,b]} (x) = \frac{x-a}{b-a}$.  Then 
$$\phi_r|_{[\frac{i}{n+1}, \frac{i+1}{n+1}]} 
  = f_l \circ L_{[\frac{i}{n+1}, \frac{i+1}{n+1}]}\rm{.)}$$

Note that due to the choice of orientations, there is some non-uniqueness in the
definition of $K_{\mathcal{P}}$.  We will ignore this and for any choice of orientations
we denote the complex obtained by $K_{\mathcal{P}}$.

From this definition, Van Kampen's Theorem can be used to show that $\mathcal{P}$ is 
indeed a presentation for $\pi_1 \left(K_{\mathcal{P}}\right)$.

Before we describe the presentation for $\pi_1$ of an arbitrary simplicial complex,
we need the following lemma.

\begin{lemma}$\label{contractible}$
Let $\Delta$ be a simplicial complex, and let $K\subset \Delta$ be a contractible 
subcomplex.  Then $\Delta \heq \Delta / K$.
\end{lemma}
{\bf Proof.}  Since $K$ is contractible, there is a map $C:K\cross I \to K$ satifying
$C_0 = id_K$ and $C_1 (K) = k_0$ for some $k_0\in K$.  By the Homotopy Extension Property,
this map extends to a map $F: \Delta \cross I \to \Delta$, with $F_0 = id_{\Delta}$.  The 
map $F_1$ factors through $\Delta /K$, giving a map $\phi: \Delta /K \to \Delta$.  

We claim that $\phi$ is a homotopy inverse to the quotient map $\pi: \Delta \to \Delta / K$.
First, $\phi \pi = F_1 \heq F_0 = id_{\Delta}$.  To show that $\pi \phi \heq id_{\Delta / K}$, consider
the map $\pi F:\Delta\cross I \to \Delta / K$.  This map factors through $\Delta /K \cross I$, 
yielding a map $\Phi : \Delta /K \cross I \to \Delta / K$.  It is now easy to check that 
$\Phi_0 = id_{\Delta / K}$ and $\Phi_1 = \pi \phi$.

$\hfill \Box$

\vspace{.15in}
The map $F$ in the above proof may be described explicitly; see~\cite{Bjorner-hom-comp}.

\begin{theorem}$\label{presentation}$  Let $\Delta$ be a simplicial complex, and 
let $T\subset \Delta^1$ be a maximal tree (i.e. a spanning tree).  

Then $\pi_1(\Delta)$ has a presentation with a generator for each (ordered) edge
$(u,v)$ with $\{u,v\}\in \Delta^1$, and with the following relations:
\begin{enumerate}
	\item $(u,v) = 1$ if $\{u,v\}\in T$,
	\item $(u,v)(v,u) = 1$ if $\{u,v\}\in \Delta^1$,
	\item $(u,v)(v,w)(w,u) = 1$ if $\{u,v,w\}\in \Delta^2$.
\end{enumerate}
\end{theorem}
{\bf Proof.}
We begin by simplifying the presentation given above.

For each 1-simplex $\{u,v\}\in \Delta^1 - T$, choose a (combinatorial) orientation 
(either $(u,v)$ or $(v,u)$)
and similarly for each 2-simplex.  (These orientations are just orderings of the vertices
and we do not require them to agree with 
one another in any way.)  Let $\mathcal{O}_1$ be the set of chosen oriented edges and let
$\mathcal{O}_2$ be the set of chosen oriented 2-simplices.

Let $F$ denote the free group generated by $\mathcal{O}_1$.  First, if 
$(u,v)$ is the chosen orientation for an edge $\{u,v\}\in \Delta$, then we let $(v,u)$
denote $(u,v)^{-1}$ (the inverse of $(u,v)$ in F).  
Also, for $\{u,v\}\in T$, we let $(u,v) \equiv (v,u) \equiv 1$ (i.e. $(u,v)$ and $(v,u)$
are formal symbols representing the identity in $F$).

Then the following presentation $\mathcal{P}$ clearly realizes the same group as the 
presentation given in the theorem:
$$\mathcal{P} = \gen{\mathcal{O}_1 | (u,v)(v,w)(w,u) = 1 {\rm \,\, if \,\,} 
                                   (u,v,w)\in \mathcal{O}_2}.$$

By Lemmas~\ref{2-skeleton} and~\ref{contractible} we have
$$\pi_1 (\Delta) \homeo \pi_1 (\Delta^2) \homeo \pi_1 (\Delta^2 / T).$$

Now, $\Delta^2 / T$ is a CW-complex whose cells are the images of the simplices in 
$\Delta^2$.  We claim that in fact, $\Delta^2 / T$ is (homeomorphic to) the standard 
complex for the presentation $\mathcal{P}$.  

Let $q: \Delta^2 \to \Delta^2 / T$ be the quotient map.
Then $\Delta^2 / T$ has a unique vertex $v_0 = q\left( V(\Delta)\right)$ and a 1-cell
for each generator in $\mathcal{O}_1$.  In defining the standard complex for $\mathcal{P}$,
it is necessary to choose an orientation on the 1-cells.  We will choose to orient the
1-cells of $\Delta^2 / T$ in accordance with the chosen orientations on the edges of $\Delta^2$,
as follows.  Say $(u,v)\in \mathcal{O}_1$.  Then let $f_{(u,v)}:I\to \Delta$ be the map 
$f_{(u,v)}(t) = (1-t)u + tv$.  It is clear that $f_{(u,v)}(I) = \abs{\{u,v\}}$.  
Now the orientation on 
the 1-cell corresponding to $\{u,v\}$ is simply $q\circ f_{(u,v)}$, and this agrees with 
the definition of $K_{\mathcal{P}}$.

Now consider a 2-cell $\tau = q\left(\abs{\{u,v,w\}
                                         }
                            \right)$, where $\{u,v,w\}\in \Delta$.  Assume that
$(u,v,w)$ is the chosen orientation on this simplex.  Then the attaching map for $c$ is just
(linear shifts of) the orienting map $f_{(u,v)}$ followed by $f_{(v,w)}$ followed by $f_{(w,u)}$,  
just as in the standard complex $K_{\mathcal{P}}$.  
$\hfill \Box$

\section{Homotopy Theory for Posets}$\label{poset-homotopy}$

In this section, we will describe a functor from the category of posets 
(partially ordered sets) 
to the category of simplicial complexes.  This will allow us to apply topological 
concepts to posets, and in particular we will be able to talk about the homotopy type
of a poset.  After introducing this functor, we will prove several results involving
homotopy equivalence and homotopy type of posets.

\begin{definition}$\label{poset}$
A \emph{poset} is a set $P$ together with a subset 
$R\subset P\times P$ (called a partial order on P) satisfying the following 
properties:
\begin{enumerate}
	\item $(x,x)\in R$ for all $x\in P$ (reflexive).
	\item If $(x,y)\in R$ and $(y,z)\in R$ then $(x,z)\in R$ (transitive).
	\item If $(x,y)\in R$ and $(y,x)\in R$, then $x=y$ (anti-symmetric).
\end{enumerate}

When $(x,y)\in R$, we write $x\leqslant y$.  If in addition $x\neq y$, we write $x<y$.
We call such elements $x$ and $y$ \emph{comparable}.

If $S\subset P$ has a least upper bound $j\in P$, we call $j$ the \emph{join} of $S$
and write $m = \join S$.  Similarly, if $S$ has a greatest lower bound $m\in P$, we call
$m$ the \emph{meet} of $P$ and write $m = \meet S$.
\end{definition}

From now on, we assume that all posets are finite.
Given a poset $P$, the $\emph{order complex}$ of $P$ (denoted $\ord{P}$) is the simplicial
complex whose vertices are the elements of $P$ and whose simplices are the 
$\emph{chains}$ of $P$, i.e. the sequences $x_1 < x_2 < \cdots < x_n$.  It is clear 
that $\ord{P}$ is a simplicial complex.

If $P$ and $Q$ are posets, a map $f: P\to Q$ is called \emph{order-preserving} if 
$x\leqs y \implies f(x)\leqs f(y)$.  Note that the composition of two order-preserving maps
is order-preserving.  An order-preserving map of posets induces a simplicial map between
order complexes in the obvious manner.  For simplicity, we use the same symbol to denote
the map $f:P\to Q$ and the induced map $\ord{P} \to \ord{Q}$, and we let 
$\abs{P}$ denote $\abs{\ord{P}}$ (this will be called the geometric realization of $P$).
Note that order-reversing maps also induce simplicial maps on order complexes.  [To keep
the relationship between $P$ and $\ord{P}$ functorial, one may think of an order
reversing map $f:P\to Q$ as an order-preserving map to the dual poset $Q^*$, in which
$x\leqs_{Q^*} y \iff y\leqs_Q x$.  Then simply note that $\ord{Q} = \ord{Q^*}$.]

\begin{definition} Let $P$ and $Q$ be posets.  We say $P$ and $Q$ are homotopy
equivalent ($P\heq Q$) if the geometric realizations of their order complexes are 
homotopy equivalent.
\end{definition}

We can now state an important theorem of Quillen~\cite[Proposition 1.6]{Quillen-p-subgroups}
giving conditions under which an
order-preserving map of posets induces a homotopy equivalence of order complexes.  

\begin{theorem}[Quillen, 1978] $\label{quillen}$
Let $P$ and $Q$ be posets and let $f:P\to Q$ be an order-preserving map.  For $x\in Q$, let
$Q_{\geqs x}$ be the set of all elements greater than or equal to $x$.  
If $\abs{f^{-1}\left(Q_{\geqs x}\right)}$ is contractible for each $x\in Q$, then $\abs{f}$ is a 
homotopy equivalence from $\abs{P} \to \abs{Q}$.
\end{theorem}

Note that any subset of a poset is again a poset, and the order complex of a subposet
is a full subcomplex of the order complex.

The proof of Quillen's Theorem presented here comes from~\cite{Bjorner-top-methods} and is an 
application of the Carrier Lemma.  
In order to prove this lemma, we will first need a basic fact about
maps into contractible spaces.  

Given a topological space $X$, the \emph{cone} on $X$
is the space 
$$C(X) = \left(X\cross I\right)/X\cross \{1\}.$$

\begin{lemma}$\label{extension}$  Let $X$ and $Y$ be topological spaces, and
assume that $Y$ is contractible.  Then any map $f:X\to Y$ can be extended
to a map $\hat{f}:C(X) \to Y$ satisfying $\hat{f}\circ q (x,0) = f(x)$ for
each $x\in X$ (here $q$ denotes the quotient map $X\cross I\to C(X)$).
\end{lemma}
{\bf Proof.}  Since $Y$ is contractible, there is a map $H:Y\cross I\to Y$ and a point 
$y_0\in Y$ s.t.
$H(y,0) = y$ and
$H(y,1) = y_0$ for each $y\in Y$.  Letting $f\cross id_I: X\cross I \to Y\cross I$ denote
the map $(x,t)\mapsto (f(x),t)$, the composition 
$H\circ \left(f\cross id_I\right): X\cross I \to Y$ 
factors through the quotient map $q$ and
gives the desired extension of $f$. $\hfill \Box$

\vspace{.15in}
Let $\Delta$ be a simplicial complex and let $X$ be any topological space.
We call a function $K:\Delta \to 2^X$ (the power set of $X$) a \emph{contractible system} if 
\begin{enumerate}
	\item $\sigma \subset \tau \implies K(\sigma)\subset K(\tau)$
	\item $K(\sigma)$ is contractible for each simplex $\sigma$.
\end{enumerate}
Also, if $f:\abs{\Delta} \to X$ is a continuous map, we say $f$ is \emph{carried} by
$K$ if $f\left(\abs{\sigma}\right)\subset K(\sigma)$ for each $\sigma \in \Delta$.
The following lemma
shows that maps carried by a particular system $K$ exist and are unique up to homotopy.

\begin{lemma}[Carrier Lemma] $\label{carrier}$
Let $\Delta$ be a simplicial complex, $X$ a topological space and $K$ a contractible
system.  Then there exists a continuous map $f:\abs{\Delta}\to X$ carried by $K$, and if
$g$ is also carried by $K$ then $f\heq g$.
\end{lemma}
{\bf Proof.}  First we construct a map $f$ carried by $K$.  We will construct $f$ inductively,
beginning with the vertices of $\Delta$.  If $v\in V(\Delta)$, let $f(v)$ be any point in 
$K(v)$.  Now assume that we have a map $f: \Delta^k \to X$ carried by $K$, and let $\sigma$
be a $k+1$-simplex.  Clearly, $\abs{\sigma} \homeo C\left(\sigma^k\right)$ and
$f(\sigma^k)\subset K(\sigma)$, which is contractible.  So Lemma~\ref{extension} applies,
and $f|_{\sigma^k}$ extends to a map $\abs{\sigma}\to K(\sigma)$.  Since $\sigma$ is finite,
we may extend the map $f$ to the entire $k+1$-skeleton of $\sigma$, and repeating this
process will eventually give a map from $\Delta$ to $X$ carried by $K$.

Next, say $g:\abs{\Delta}\to X$ is also carried by $K$.  We need to construct a homotopy 
$F:\abs{\Delta}\times I \to X$ with $F_0 = f$ and $F_1 = g$.  This can be done in a manner
similar to the construction of $f$ above.  Again we define $F$ inductively on successive skeleta
of $\Delta$.  Note that it suffices to check the continuity of $F$
on the set of $\abs{\sigma} \times I$, $\sigma \in \Delta$ (these sets are closed in the product
topology, and cover $\Delta \times I$).

If $v$ is a vertex of $\Delta$, then $f(v), g(v)\in C(v)$ and since $K(v)$ is contractible, it 
is path-connected, i.e. there is a map $H:I\to K(v)$ s.t. $H(0) = f(v)$ and $H(1) = g(v)$.  We 
define $F|_{\{v\}\cross I} = H\circ p$, where $p$ is the projection $\{v\}\cross I \to I$.  Now, 
if $F$ has been appropriately defined on $\Delta^k\cross I$ and $\sigma$ is a 
$k+1$-simplex of $\Delta$, then $F|_{\abs{\sigma^k}\cross I}$ extends as before to a map of 
$C\left(\abs{\sigma^k}\cross I\right)\homeo \abs{\sigma}\cross I$.
$\hfill \Box$

\vspace{.15in}
\noindent {\bf Proof of Theorem~\ref{quillen}.}
Let $P$ and $Q$ be posets, and let $f:P\to Q$ be an order-preserving map s.t. for each 
$x\in Q$, $f^{-1}\left(Q_{\geqs x}\right)$ is contractible.  For any simplex 
$\sigma\in \ord{Q}$, define 
$K(\sigma) = \abs{f^{-1}\left(Q_{\geqs \min(\sigma)}\right)} \subset \abs{P}$ 
(note that $\sigma$
is a chain in $Q$ and thus has a minimum element).  $K$ is a contractible
system by hypothesis, and Lemma~\ref{carrier} gives us a map $g:\abs{Q} \to \abs{P}$ 
carried by $K$.  We claim that $g$ is a homotopy inverse to $\abs{f}$ (which from here
on we denote simply by $f$).

First we show that $g\circ f\heq id_P$.  Consider the contractible system 
$K':P\to 2^{\abs{P}}$ given by $K'(\tau) = \abs{f^{-1} \left(Q_{\geqs \min f(\tau)}\right)}$.
We claim that $K'$
carries both $g\circ f$ and $id_P$.  Clearly $K'$ carries $id_P$, since if $\sigma$
is a chain in $P$, then $f(\sigma)$ lies above $\min f(\sigma)$ and hence 
$\sigma \subset K(\sigma)$.  Also, if $\tau \in \ord{P}$ then 
$$g\circ f(\tau) 
     \subset K\left(f(\tau)
              \right) 
             = \abs{f^{-1}\left(
                                Q_{
                                   \geqs \min\left(
                                                  f(\tau)
                                            \right)
                                  }
                          \right)
                    }
             =K'(\tau).$$ 
Lemma~\ref{carrier} now shows that $g\circ f\heq id_P$.

The proof that $f\circ g\heq id_Q$ is similar.  We define the contractible system
$K'':Q\to 2^{\abs{Q}}$ by $K''(\sigma) = \abs{Q_{\geqs \min(\sigma)}}$ and as before,
it is clear that $K''$ carries $id_Q$.  $K''$ also carries $f\circ g$, which can 
be checked as follows:
$$f\circ g(\sigma) \subset f\left(K(\sigma)\right)
                   =f \left(
                            \abs{
                                 f^{-1} \left(
                                              Q_{\geqs \min(\sigma)}
                                       \right)
                                }
                      \right)
                   \subset \abs{Q_{\geqs \min(\sigma)
                                  }
                               }
                   = K''(\sigma).$$
$\hfill \Box$

\vspace{.15in}
Our next goal is to prove the Nerve Theorem, which will be useful in our study of the
coset poset.  Before proving this theorem, need a technical result about simplicial complexes.

If $\Delta$ is a simplicial complex, the \emph{face poset} of $\Delta$ is the poset 
$P(\Delta)$ whose elements are the simplices of $\Delta$, ordered by inclusion.  
The simplicial complex $sd(\Delta) = \ord{P(\Delta)}$ is the 
\emph{first barycentric subdivision} of $\Delta$ and we have:

\begin{proposition}$\label{bary}$
If $\Delta$ is a simplicial complex, then $\abs{\Delta}$ is homeomorphic to 
$\abs{sd(\Delta)}$.
\end{proposition}
{\bf Proof.}  We will define an explicit homeomorphism 
$h:\abs{sd(\Delta)}\to \abs{\Delta}$.
First we define $h$ on the vertices of $\abs{sd(\Delta)}$, again using the convention
that the standard basis vectors of the vector spaces containing $\abs{\Delta}$ and 
$\abs{sd(\Delta)}$ are identified with the vertices of these complexes.  

For any $\sigma = \{v_1,\ldots, v_n\} \in V(sd(\Delta)) = \Delta$, we define 
$$h(\sigma):= \frac{1}{\dim{\sigma} + 1}\sum_{i=1}^n v_i.$$
(The point $h(\sigma)$ is called the \emph{barycenter} of $\sigma$.)
Now we extend $h$ by linearity (so we may view $h$ as a linear transformation  
of vector spaces),
so that for an arbitrary point $t = \sum_{i=1}^{k} t_i \sigma_i \in \R^{\Delta}$,
$$h(t)	= \sum_{i=1}^{k} t_i h(\sigma_i).$$

Next we show that $h(\abs{sd(\Delta)})\subset \abs{\Delta}$.  Consider an arbitrary
point $t\in \abs{sd(\Delta)}$.  We may write $t = \sum_{i=1}^k t_i \sigma_i$, 
and the simplices $\sigma_i$ must form a simplex in $sd(\Delta)$, i.e. they must form
a chain in the face poset $P(\Delta)$.  Without loss of generality we may assume that
$\sigma_1\subset \sigma_2\subset \cdots \subset \sigma_k$.  
Let $\sigma_k = \{v_1,\ldots, v_n\}$, and let $\sigma_i = \{v_1,\ldots,v_{n_i}\}$ (so
the $n_i$'s form an increasing sequence).

Then we have
$$h(t) = \sum_{i=1}^{k} t_i h(\sigma_i)
       = \sum_{i=1}^{k} \frac{t_i}{n_i} \sum_{j=1}^{n_i} v_j,$$
and thus the total sum of the coefficients of $h(t)$ is 
$$\sum_{i=1}^{k} \frac{t_i}{n_i} n_i = \sum_{i=1}^k t_i = 1.$$

Also, all the terms above are positive, so the coefficient of any particular $v_i$ is 
positive.  Finally, the support of $h(t)$ is contained in $\sigma_k$ and is thus a simplex
of $\Delta$.  So $h(t)\in \abs{\Delta}$ as desired.

Since $h$ is continuous on $\R^{\Delta}$, 
its restriction to $\abs{sd(\Delta)}$ is a continuous map from $\abs{sd(\Delta)}$ to
$h(\abs{sd(\Delta)})$.  Thus it remains to show that $h$ is in fact a bijection between
$\abs{sd(\Delta)}$ and $\abs{\Delta}$.  This will complete the proof, because a continuous
bijection of compact hausdorff spaces is a homeomorphism.

First we show that $h$ is injective.  We will actually prove a bit more than this.

For any poset $P$, let 
$$\widehat{\abs{P}} = \left\{\sum_{i=1}^n \lambda_i p_i \in \R^P : 
                       \lambda_i > 0, \,\,\, 0 < \sum_{i=1}^n \lambda_i \leqs 1 \,\,\, 
                       {\rm and} \,\,\, p_1< \cdots <p_n \right\}.$$
If $x\in \widehat{\abs{P}}$, we may write $x$ uniquely as $\sum_{i=1}^n \lambda_i p_i$ 
(with $p_1 <\cdots <p_n$), and we define $\nu (x) = n$ and $\lambda (x) = \sum_{i=1}^n \lambda_i$.

{\bf Claim:}  The restriction of $h$ to $\widehat{\abs{sd(\Delta)}}$ is injective.

\vspace{.1in}
\noindent {\bf Proof of Claim.} Say $s,t\in \widehat{\abs{sd(\Delta)}}$ and 
$h (s) = h (t)$.  We will show, by induction on $m(s,t) := \min(\nu(s),\nu(t))$,
that $s = t$.

First, say $m(s,t) = 1$.  Assume that $\nu(s)=1$, so that $s = \lambda (s) \sigma$ for some
$\sigma \in \Delta$.  
Letting $\sigma = \{v_1,\ldots,v_p\}$, we have 
$$h(s) = \frac{\lambda (s)}{p} \sum_{i=1}^p v_p.$$

Next, let $t = \sum_{i=1}^l t_i \tau_i$, with $\tau_1 <\cdots < \tau_l$ 
(so $l = \nu (\tau)$).  Also, let $\tau_i = \{w_1,\ldots, w_{q_i}\}$.
Then we have
$$h(t) = \sum_{i=1}^{l} \frac{t_i}{q_i} \sum_{j=1}^{q_i} w_j,$$
and we see that the support of $h(t)$ is exactly $\{w_1,\ldots,w_{q_l}\}$ 
(since $t_i>0$ for $i = 1,\ldots,l$).

Thus, since $h(s) = h(t)$, $\{w_1,\ldots,w_{q_l}\} = \sigma$, and 
(by renumbering the $v_i$) we may assume $w_i = v_i$.  Equating
the coefficients of $v_p = w_{q_l}$ in the expressions for $h(s)$ and $h(t)$ gives
$$\frac{\lambda (s)}{p} = \frac{t_l}{q_l} = \frac{t_l}{p}$$
and thus $t_l = \lambda (s)$.  So the $l$-th term of the expression for $h(t)$ is 
in fact the entire expression 
(because $t_l = \lambda (s) = \lambda (h (s)) = \lambda (h(t))$),
so $l = 1$ and $t = \lambda (s) \sigma = s$.

Now assume the Claim for any pair of points $s',t'$ (satisfying the hypotheses)
with $m(s,t) = n$ ($n\geq 1$), and consider points $s,t$ with $m(s,t) = n+1$.

By the same reasoning as above, if $s = \sum_{i=1}^k s_i \sigma_i$ and 
$t = \sum_{i=1}^l t_i \tau_i$ are in $\widehat{\abs{sd(\Delta)}}$
(where $\sigma_1\subset\cdots\subset\sigma_k$ and
$\tau_1\subset\cdots\subset\tau_k$) then we have $\sigma_k = \tau_l$.  Let this
simplex be $\{v_1,\ldots,v_p\}$, and let $\sigma_i = \{v_1,\ldots,v_{p_i}\},
\tau_i = \{v_1,\ldots,v_{q_i}\}$.  Then, since $h(s) = h(t)$, we have
$$\sum_{i=1}^k \frac{s_i}{p_i} \sum_{j=1}^{p_i} v_j 
 =\sum_{i=1}^l \frac{t_i}{q_i} \sum_{i=1}^{q_i} v_j,$$
and by equating the coefficients of $v_p = v_{p_k} = v_{q_l}$ we see that 
$\frac{s_k}{p_k} = \frac{t_l}{q_k}$, so $s_k = t_l$ (since $p_k = q_k = p$) and thus
$$h\left(\sum_{i=1}^{k-1} s_i \sigma_i\right) = h\left(\sum_{i=1}^{l-1} t_i \sigma_i\right).$$

The induction hypothesis applies and gives us that 
$\sum_{i=1}^{k-1} s_i \sigma_i = \sum_{i=1}^{l-1} t_i \sigma_i$.  Adding back in
the last terms gives $s = t$, completing the proof of the claim. 

$\hfill \Box$

\vspace{.15in}

\noindent {\bf Proof of Proposition~\ref{bary} (Continued)}
To complete the proof, we must show that $h$ is surjective.

Consider an arbitrary point $\sum_{i=1}^n a_i v_i \in \abs{\Delta}$, and assume that
$a_1\geq a_2\geq \cdots \geq a_n$.

Let $t = \sum_{i=1}^n t_i \sigma_i$, where $\sigma_i = \{v_1,\ldots,v_i\}$.  Then
$$h(t) 	= \sum_{i=1}^n t_i h(\sigma_i) 
	= \sum_{i=1}^n \frac{t_i}{i} \sum_{j=1}^i v_j
	= \sum_{i=1}^n v_i \sum_{j=i}^n \frac{t_j}{j},$$
and we want to choose $t_1,\ldots,t_n$ s.t. 
\begin{equation} \label{matrix}
\sum_{j=i}^n \frac{t_j}{j} = a_i
\end{equation}
(and s.t. $t\in \abs{sd(\Delta)}$).  Note that we may choose $t_i = 0$
for some $i$.

Equation~\ref{matrix} gives us $\frac{t_n}{n} = a_n$, i.e. $t_n = n a_n$.

Next, we have $t_{n-1} = (n-1) (a_{n-1} - \frac{t_n}{n}) = (n-1) (a_{n-1} - a_n)$,
and assuming inductively that $t_j = j(a_j - a_{j+1})$ for $j>i$ (and with the
convention that $a_{n+1} = 0$) we see that
\begin{eqnarray*}
t_i 	& = & i\left(a_i - \frac{t_n}{n} - \frac{t_{n-1}}{n-1} -\cdots 
                                         - \frac{t_{i+1}}{i+1}\right) \\
	& = & i ( a_i - a_n - (a_{n-1} - a_n) - (a_{n-2} - a_{n-1}) -\cdots \\
        &   & - (a_{i+2} - a_{i+3}) - (a_{i+1} - a_{i+2}) ) \\
	& = & i\left(a_i - a_{i+1}\right).
\end{eqnarray*}

Since we have assumed that $a_1\geq\cdots\geq a_n$, all of the coefficients determined
above are non-negative.  Thus to show that $t\in \abs{sd(\Delta)}$,
it remains only to check that their sum is exactly 1.
We have:
\begin{eqnarray*}
\sum_{i=1}^n t_i & = & \sum_{i=1}^n i(a_i - a_{i+1}) \\
	 	 & = & (a_1 - a_2) + 2(a_2 - a_3) + \cdots \\
                 &   & + (n-1)(a_{n-1} - a_n) + n a_n \\
		 & = & \sum_{i=1}^n a_i = 1,
\end{eqnarray*}
the last step following from the fact that $\sum_{i=1}^n a_i v_i \in \abs{\Delta}$.
$\hfill \Box$
\vspace{.15in}

Given a finite collection $\{V_i\}_{i\in I}$ of subsets of a set $S$, we define the
\emph{nerve} of ${V_i}$ to be the simplicial complex $\mathcal{N}(\{V_i\})$ with vertices
$I$ and with a simplex for each set $J\subset I$ s.t. $\bigcap_{j\in J} V_j \neq \emptyset$.

\begin{theorem}[Nerve Theorem]
$\label{nerve}$
Let $\Delta$ be a simplicial complex and let $\{\Gamma_i\}_{i\in I}$ be a collection of 
subcomplexes.  If $\{\Gamma_i\}$ satisfies the following two properties, then 
$\Delta\heq \mathcal{N}\left(\{\Gamma_i\}\right)$:
\begin{enumerate}
	\item $\bigcup_{i\in I} \Gamma_i = \Delta$.
	\item For any $J\subset I$, the intersection $\bigcap_{j\in J} \Gamma_j$ 
is contractible.
\end{enumerate}
\end{theorem}
{\bf Proof.}  By Proposition~\ref{bary}, it suffices to show that the subdivisions 
(or face posets) of these complexes are homotopy equivalent.  
Letting $P(\Delta)^*$ denote the dual poset to 
$P(\Delta)$ (see the remark after Definition~\ref{poset}), 
we have $\ord{P(\Delta)^*} = sd(P(\Delta))$.  We will apply Quillen's Theorem
(Theorem~\ref{quillen}) to show that $P(\Delta)^* \heq P(\mathcal{N})$, where 
$\mathcal{N} = \mathcal{N}\left(\{\Gamma_i\}_{i\in I}\right)$.

Consider the map $f: P(\Delta)^* \to P(\mathcal{N})$ induced by
$$f(\sigma) = \{i\in I: \sigma \in \Gamma_i\},$$
where $\sigma \in V\left(P(\Delta)\right) = \Delta$.  
[Note that $f(\sigma)$ is a simplex of $\mathcal{N}$ (and hence an element of 
$P(\mathcal{N})$) because $\sigma \in \bigcap_{i\in f(\sigma)} \Gamma_i$.]
Now, $f$ is clearly order preserving (since we are using the dual poset $P(\Delta)^*$)
and for any $J\subset I$ we have
\begin{eqnarray*}
	f^{-1}\left(
		    P\left(
			   \mathcal{N}
		     \right)_{\geqs J}
	      \right)			& = & \{\sigma \in \Delta : J\subset f(\sigma)\}\\
					& = & \{\sigma \in \Delta : 
    					          \sigma \in \bigcap_{j\in J} \Gamma_j\}\\
					& = & \bigcap_{j\in J} \Gamma_j,
\end{eqnarray*}
which is contractible by assumption.  Quillen's Theorem shows that $\abs{f}$ 
is a homotopy equivalence, completing the proof. $\hfill \Box$

\vspace{.15in}
In light of Quillen's Theorem and the Nerve Theorem, contractible posets and complexes are of 
great importance.  Our next goal is to give a simple but useful condition under which
a poset is contractible.  Although the result could be proven by describing an 
explicit map, the following approach is more combinatorial in nature and will lead us through
a number of other interesting results.  In particular, Corollary~\ref{order-homotopic2} will
be of use to us in Section~\ref{extensions}.

\begin{definition}  Let $P$ and $Q$ be posets.  We define the product of $P$ and $Q$
to be the set $P\cross Q$ together with the ordering 
$(x,y)\leqs (v,w) \iff x\leqs v$ and $y\leqs w$. 
\end{definition}

It is easily checked that this is a partial order on $P\cross Q$.

\begin{lemma}$\label{poset-products}$
For any finite posets $P$ and $Q$, we have $\abs{P\cross Q}\homeo \abs{P}\cross \abs{Q}$.
\end{lemma}
{\bf Proof.}  The projections $P\cross Q \to P$ and $P\cross Q \to Q$ induce a map
$f: \abs{P\cross Q}\to \abs{P}\cross \abs{Q}$, with 
$$f\left( \sum_{i=1}^n \lambda_i (p_i, q_i)\right) 
= \left(\sum_{i=1}^n \lambda_i p_i, \sum_{i=1}^n \lambda_i q_i \right).$$
Note that we may also consider $f$ to be the restriction to $\abs{P\cross Q}$
of the map $\R^{P\cross Q}\to \R^P\oplus \R^Q$ sending $(p,q)\in \R^{P\cross Q}$ to 
$(p,q)\in \R^P\oplus \R^Q$.

Since the spaces in question are compact Hausdorff, it suffices to prove that the map
$f$ is a bijection.  First we prove that $f$ is surjective.

\vspace{.1in}
\noindent {\bf Claim:} For any point $(x,y)\in \abs{P}\cross \abs{Q}$, 
$\exists$ $\lambda_1, \ldots, \lambda_n \in (0,1]$, $p_1\leqs \cdots \leqs p_n \in P$ and
$q_1\leqs \cdots \leqs q_n \in Q$ s.t. 
$$(x,y) = \left(\sum_{i=1}^n \lambda_i p_i, \sum_{i=1}^n \lambda_i q_i\right).$$

\vspace{.1in}
Assuming the Claim, we see that $f$ is surjective because if $(x,y)\in \abs{P}\cross \abs{Q}$
and $(x,y) = \left(\sum_{i=1}^n \lambda_i p_i, \sum_{i=1}^n \lambda_i q_i\right)$ 
(where this expression satisfies the conditions in the Claim) then we have 
$\sum_{i=1}^n \lambda_i (p_i, q_i) \in \abs{P\cross Q}$ and 
$$f\left(\sum_{i=1}^n \lambda_i (p_i, q_i)\right) 
= \left(\sum_{i=1}^n \lambda_i p_i, \sum_{i=1}^n \lambda_i q_i\right)
	= (x,y).$$

To prove the Claim, we will actually verify a stronger statement. 

\vspace{.1in}
\noindent {\bf Claim 2:} For any point $(x,y) \in \widehat{\abs{P}} \cross \widehat{\abs{Q}}$ 
with $\lambda (x) = \lambda (y)$ (see the proof of Proposition~\ref{bary} for an explanation
of this notation), $\exists$ $\lambda_1, \ldots, \lambda_n \in (0,1]$, 
$p_1\leqs \cdots \leqs p_n \in P$ and
$q_1\leqs \cdots \leqs q_n \in Q$ s.t. 
$$(x,y) = \left(\sum_{i=1}^n \lambda_i p_i, \sum_{i=1}^n \lambda_i q_i\right).$$
{\bf Proof.}  We proceed by induction on $\nu (x) + \nu (y) =: \nu (x,y)$.  
If $\nu = 2$ there is nothing to prove, so we may assume the claim for all points $(x,y)$ 
with $\nu (x,y) \leqs r$ (and $\lambda (x) = \lambda (y)$ ) and 
consider a point $(x,y)$ with $\nu (x,y) = r + 1$.  Writing
$$(x,y) = \left(\sum_{i=1}^k t_i p_i, \sum_{i=1}^m s_i q_i\right)$$
with $p_1<\cdots <p_n$ and $q_1<\cdots <q_n$, we may assume WLOG that $t_1 \leqs s_1$.  
We now have
$$(x,y) - (t_1 p_1, t_1 q_1) 
= \left(\sum_{i=2}^k t_i p_i, (s_1 - t_1) q_1 + \sum_{i=2}^m s_i q_i\right)$$
and it is easy to check that the right-hand side satisfies the induction hypothesis.  The claim 
now follows easily.
$\hfill \Box$

\vspace{.15in}
We must now prove that $f$ is injective.
Say $x,y\in \widehat{\abs{P\cross Q}}$.  Then we may write $x$ and $y$ uniquely in 
the form
$$x = \sum_{i=1}^n \lambda_i (p_i, q_i), \,\,\, y = \sum_{i=1}^m \gamma_i (p'_i, q'_i),$$
with $\lambda_i, \gamma_i > 0$, $(p_1, q_1) < \cdots < (p_n, q_n)$ 
and $(p'_1, q'_1) < \cdots < (p'_m, q'_m)$.

Assume that $f(x) = f(y)$, i.e. that 
$$\left(\sum_{i=1}^n \lambda_i p_i, \sum_{i=1}^n \lambda_i q_i\right) 
= \left(\sum_{i=1}^m \gamma_i p'_i, \sum_{i=1}^m \gamma_i q'_i\right).$$
Note that $\sum_{i=1}^n \lambda_i p_i = \sum_{i=1}^m \gamma_i p'_i 
  \implies \{ p_1\ldots p_n\} = \{ p'_1,\ldots p'_m\} \implies p_1 = p'_1$, 
and similarly $q_1 = q'_1$.  

We will show that $\lambda_1 = \gamma_1$.  
Since $(p_1, q_1) < (p_2, q_2)$, we may assume WLOG that $p_1 < p_2$.  Writing
$f(x) = f(y)$ in terms of the standard basis $\{(p,0)\}_{p\in P}\cup \{(0,q)\}_{q\in Q}$ 
for $\R^P \oplus \R^Q$, we see that 
$$\lambda_1 = \sum_{i: p'_i = p'_1} \gamma_i \geqs \gamma_1.$$  

If $p'_1\neq p'_2$, then $\lambda_1 = \sum_{i: p'_i = p'_1}^n \gamma_i = \gamma_1$ 
(since $p'_1 < p'_2 \leqs \cdots \leqs p'_m$), and we 
are done.  So assume
$p'_1 = p'_2$.  Then $q'_1 < q'_2$, and the coefficient on $(0, q'_1) = (0, q_1)$ in $f(x) = f(y)$
is 
$$\gamma_1 = \sum_{i: q_i = q_1}^n \lambda_i \geqs \lambda_1.$$
But earlier we saw that $\gamma_1 \leqs \lambda_1$, so we must have
$\lambda_1 = \gamma_1$.  Repeating this process will eventually show that $x = y$. 
$\hfill \Box$

\vspace{.15in} For an alternate proof of the above result, see 
Eilenberg and Steenrod~\cite{E-S-alg-top}.

\begin{lemma}$\label{order-homotopic1}$
Let $P$ and $Q$ be posets, and let $f,g:P\to Q$ be order-preserving maps.  If
$f(x)\leqs g(x)$ for each $x\in P$, then $\abs{f}\heq \abs{g}$.
(We write $f \leqs g$ when the above conditions hold.)
\end{lemma}
{\bf Proof.}  Let $\mathcal{I} = \{0,1\}$, with $0 < 1$  (so $\abs{\mathcal{I}}$ 
is homeomorphic to the unit interval $I$).  Then since $f \leqs g$, these maps induce
an order-preserving map $h$ on $P\cross \mathcal{I}$, with $h(x,0) = f(x)$, $h(x,1) = g(x)$.
Now $\abs{h}$ is a map from 
$\abs{P\cross \mathcal{I}} \isom \abs{P}\cross I\to \abs{Q}$, and is thus a homotopy
from $f$ to $g$.
$\hfill \Box$

\begin{claim}$\label{cones-contract1}$
If a poset $P$ has a minimum element, then $P$ is contractible.
\end{claim}
{\bf Proof.} Let $p_0$ be the minimum element of $P$.  If
$f_{p_0} : P\to P$ is the constant map to $p_0$,
we have $f_{p_0} \leqs id_P$ and the lemma shows that $P$ is contractible.

$\hfill \Box$

\vspace{.15in} Claim~\ref{cones-contract1} (as well as the generalization below) 
will be used implicitly throughout this thesis.
Using this result, we can extend Lemma~\ref{order-homotopic1}.  
Our proof is taken from Bj\"{o}rner~\cite{Bjorner-hom-comp}.

\begin{definition}  Let $f$ and $g$ be order preserving maps from a poset $P$ to a poset
$Q$.  If $f(x)$ and $g(x)$ are comparable for all $x\in P$, then we call $f$ and $g$ 
\emph{order-homotopic}.
\end{definition}

\begin{corollary}$\label{order-homotopic2}$
If $f$ and $g$ are order-homotopic maps from a poset $P$ to a poset $Q$,
then $\abs{f}\heq \abs{g}$.
\end{corollary}
{\bf Proof.}  We will apply the Carrier Lemma.  Let $C:\ord{P}\to 2^{\abs{Q}}$ be the
map $C(\sigma) = f(\sigma)\cup g(\sigma)$, $\sigma \in \ord{P}$.  Note that
$C(\sigma)$ always has a minimum element, since $\min{f(\sigma)} = f(\min{\sigma})$ and 
$\min{g(\sigma)} = g(\min{\sigma})$ are comparable.  So $C$ is a contractible system, and
since it carries both $\abs{f}$ and $\abs{g}$ we have $\abs{f}\heq \abs{g}$.
$\hfill \Box$

\vspace{.15in}  Using Corollary~\ref{order-homotopic2} we can generalize 
Lemma~\ref{cones-contract1}.

\begin{definition}  Let $P$ be a poset with an element $p_0$ which is comparable to 
each $x\in P$.  Then we call $P$ a \emph{cone on $p_0$}, or simply a cone.
\end{definition}

\begin{corollary}$\label{cones-contract2}$
Any cone $P$ is contractible.
\end{corollary}
{\bf Proof.}  Let $P$ be a cone on some $p_0\in P$.  Let $f_{p_0}:P\to P$ be the 
constant map to $p_0$.  Then for any
$x\in P$, $f_{p_0}(x) = p_0$ is comparable to $id_P (x) = x$.  So $f_{p_0}$ and $id$
are order-homotopic and hence $P$ is contractible.
$\hfill \Box$

\begin{remark}  Cones in the above sense are easily seen to be topological cones
at the point $p_0$ (in the sense of Lemma~\ref{extension}).  
This gives another proof of contractibility.
\end{remark}


\chapter{Connectivity of the Coset Poset}$\label{connectivity}$

In this chapter we introduce the coset poset of a finite group and study its
homotopy type, using the methods of Chaper~\ref{background}.  We will focus
mainly on the connectivity of the coset poset, and in particular on its fundamental
group.  

Let $G$ be a finite group.  The coset poset of $G$, denoted $\C{G}$, is the poset consisting of all 
left cosets of all proper subgroups of $G$, ordered by inclusion (so we allow cosets of the 
trivial subgroup, but do not allow $G$ itself).  Note that the choice of left cosets rather than
right cosets is irrelevant, since every left coset is a right coset, and vice versa:
$xH = (xHx^{-1}) x$.
The coset poset was introduced
in Brown~\cite{Brown-coset-poset}.  
This thesis is the result of an attempt to answer a question from~\cite{Brown-coset-poset},
asking, ``For what groups $G$ is $\C{G}$ simply connected?''

In this chapter, we will present several results describing conditions under which the 
coset poset is or is not simply connected.

We will also need to consider the subgroup poset of a finite group $G$.  This poset consists
of all proper, non-trivial subgroups of $G$, ordered by inclusion, and is denoted by $S(G)$.
We will often use $\C{G}$ and $S(G)$ to denote not only the posets but also the
corresponding simplicial complexes and topological spaces.

\section{Conditions for Simple Connectivity} 
In this section we will describe two results giving conditions under which the coset poset of
a finite group $G$ is simply connected (or $k$-connected for higher $k$).  Our results actually
extend to a wider class of posets, which we will introduce after discussing some topological 
preliminaries.

\begin{definition} Let $X$ be a topological space.  We say $X$ is $n$-connected, 
$n = 0, 1, 2, \ldots$, if for $0\leqs k\leqs n$ each map $f:S^{k}\to X$ extends to a map 
$\hat{f} :B^{k+1}\to X$.
We say that $X$ is simply connected if $X$ is 1-connected.
\end{definition}

\vspace{.1in}
Note that if $X$ is path-connected, then $X$ is $k$-connected iff the homotopy group 
$\pi_k (X)$ is trivial (see Hatcher~\cite{Hatcher-alg-top}).  
Also, $0$-connectivity is the same as path-connectivity.

\begin{lemma}$\label{full-skeleton}$ 
Let $\Delta$ be a simplicial complex and say there exists an $n$-connected
simplicial complex $\Gamma$ with $\Delta^{n+1} = \Gamma^{n+1}$.  
Then $\Delta$ is $n$-connected.
\end{lemma}
{\bf Proof.} Say $0\leqs k\leqs n$, and let 
$f:S^{k}\to \abs{\Delta}$ be any map.  By the Cellular Approximation Theorem,
there is a map $f':S^{k}\to \abs{\Delta}$ with $f'(S^{k})\subset \abs{\Delta^{k}
                                                                         }$ 
and with
$f\heq f'$.  Let $F:S^{k}\cross I\to \abs{\Delta}$ be a homotopy from $f$ to $f'$.

Now since $\Delta^{n+1} = \Gamma^{n+1}$ and $k\leqs n$, we can consider $f'$ as a map 
$S^{k}\to \abs{\Gamma}$.
Since $\Gamma$ is $n$-connected, $f'$ extends to a map 
$f'':B^{k+1}\to \abs{\Gamma}$.  Applying the Cellular Approximation Theorem again 
yields a map 
$f''' : B^{k+1} \to \abs{\Gamma^{k+1}} = \abs{\Delta^{k+1}}$.  
This map is homotopic (rel $S^{k}$) to
$f'$, and thus is an extension of $f'$ to $B^k$.  Combining the homotopy $F$ with
the map $f'''$ now yields an extension of $f$ to $B^{k+1}$.
$\hfill \Box$

\vspace{.15in}
Note that Lemma~\ref{full-skeleton} could have been stated, ``A simplicial complex
is $n$-connected $\iff$ its $(n+1)$-skeleton is $n$-connected."

We now introduce the concept of an \emph{atomized poset}.  Essentially, an atomized poset
is a poset in which each set of minimal elements 
``generates" a well-defined element of the poset (or generates the whole poset).

\begin{definition}  We call a poset $P$ \emph{atomized} if for each set $S\subset P$ consisting
of minimal elements, the subposet $P_{\geqs S} = \{p\in P: p\geqs s \, \, \forall s\in S\}$ either
is empty or has a minimum element.  

We call the minimal elements of $P$ \emph{atoms}, and denote the set of atoms of $P$
by $\A(P)$.  If $S\subset \A (P)$,
we say that $S$ \emph{generates} the minimum element of $P_{\geqs S}$ 
(or that $S$ generates $P$ if $P_{\geqs S}$ is empty), and we write $\gen{S}$
for the element generated by $S$ ($\gen{S} = P$ is allowed).
\end{definition}

\vspace{.15in}
It is easy to check that the coset poset and the subgroup poset (the poset of non-trivial, proper
subgroups) of a finite group are atomized.  The definition of generation coincides with the
standard group theoretic definitions (where the
coset generated by elements $x_1,\ldots x_n\in G$ 
is $x_1\gen{x_1^{-1} x_2,\ldots x_1^{-1} x_n}$).

For many purposes, the coset poset is too large to analyze.  The following lemma shows 
that, up to homotopy, we can replace any atomized poset $P$ with a smaller simplicial complex, 
$\M{P}$.  
This complex has many fewer vertices, although it has much higher dimension, and will play
a crucial role in our analysis of the coset poset.

\begin{lemma} $\label{min-cover}$
Let $P$ be an atomized poset and let $\M{P}$ denote the simplicial complex with vertex set
$\A (P)$ and with a simplex for each set $S\subset \A (P)$ s.t. $\gen{S} \neq P$.  
Then $\ord{P} \heq \M{P}$.
\end{lemma}
{\bf Proof.}  Consider the $\emph{minimal cover}$ of $P$, i.e. the collection of 
all cones $P_{\geqs x}$ with $x \in \A (P)$.  
We claim that this is a contractible cover
of $\ord{P}$.  Indeed, if $S\subset \A (P)$ and $\gen{S}\neq P$, then
$\cap_{s\in S} P_{\geqs s} = P_{\geqs \gen{S}} \heq \{\ast\}.$  So each intersection is either
empty or contractible, and the Nerve Theorem tells us that $\ord{P}$ is homotopy equivalent 
to the nerve of this cover, which is exactly $\M{P}$.
$\hfill \Box$
\vspace{.15in}

When $P = \C{G}$ for some finite group $G$, we denote $\M{\C{G}}$ by $\M{G}$.  This complex
has vertex set $G$ and a simplex for each subset of $G$ contained in a proper coset.

Our main theorem in this section is the following:

\begin{theorem}$\label{k-connected-atom}$
Let $P$ be an atomized poset such that no $k$ atoms generate $P$.  Then
$P$ is $(k-2)$-connected.
\end{theorem}
{\bf Proof.}  By the previous lemma, it suffices to check that $\M{P}$ is $(k-2)$-connected.
Since no $k$ atoms generate $P$, any $k$ atoms form a simplex in $\M{P}$.  So
$\M{P}^{k-1}$ is just the $(k-1)$-skeleton of the full simplex on the vertex set $\A (P)$.  
Since any simplex is 
contractible, Lemma~\ref{full-skeleton} shows that $\M{P}$ is $(k-2)$-connected.
$\hfill \Box$

\vspace{.15in}
We have the following immediate corollary of Theorem~\ref{k-connected}

\begin{corollary} $\label{k-connected}$
Let $G$ be a finite group in which any $k$ elements generate a proper subgroup.
Then $S(G)$ is $(k-2)$-connected and $\C{G}$ is $(k-1)$-connected.  In fact,
$S(G)$ is $(k-2)$-connected so long as any $k$ elements of prime order generate
a proper subgroup.
\end{corollary}

\vspace{.15in}
Corollary~\ref{k-connected} specializes to the following result giving conditions
under which $\C{G}$ (or any atomized poset) is simply connected.

\begin{corollary}$\label{2-generator}$
If $G$ is a finite group in which any two elements generate a proper subgroup, then
$\C{G}$ is simply connected.  More generally, if $P$ is an atomized poset and
no three atoms generate $P$, then $\abs{P}$ is simply connected.
\end{corollary}

As we will see in the next chapter, this theorem does not characterize finite groups with
simply connected coset posets.  In fact, the first non-abelian simple group, $A_5$, affords
a counter-example.

We now give an alternate proof of Corollary~\ref{2-generator}.  This proof will rely solely 
on the combinatorial structure of $\ord{P}$ itself (as opposed to that of $\M{P}$).  First 
we introduce some terminology.

By an edge path we mean a map $f:I\to \ord{P}^1$ s.t. for some $n$ the restrictions
$f|_{[k/n, (k+1)/n]}$ are just a linear shifts of the canonical 
homeomorphisms $I\to \abs{\{p_i, p_j\}}$ 
(where the sets $\{p_i, p_j\}$ are 2-simplices of $\ord{P}$ and $k = 0, \ldots , n-1$).
In a simple, closed edge path we require further that $f(0) = f(1)$ and that the $p_i$ are
otherwise distinct.  From now on we will refer to simple, closed edge paths simply as cycles.

\vspace{.15in}
\noindent {\bf Alternate Proof of Corollary~\ref{2-generator}.}

Let $P$ be an atomized poset in which no three atoms generate.  First,
note that $\ord{P}$ is connected since every element is connected by an edge
to an atom, and between any two atoms $a_1, a_2\in \A (P)$ we have the path
$a_1 < \gen{a_1, a_2} > a_2$.  We now need to show that $\pi_1 (\abs{P}) = 1$, and
it will suffice to show that any cycle is null-homotopic
(this is a consequence of the Simplicial Approximation Theorem.)  

Our first step will be to consider \emph{reduced} cycles.  
Any cycle can be represented by listing 
the vertices through which it passes.  Let $C$ be a cycle and say $C$ passes through
$p_1, p_2, \ldots , p_n = p_1$ (in that order).  For each $i$ we have either $p_i < p_{i+1}$
or $p_i > p_{i+1}$, and we call a cycle \emph{reduced} if these inequalities alternate direction 
(they must alternate at the base-point as well).  

We claim that any cycle is homotopic to a reduced cycle.  If, in a cycle $C$, some inclusion 
does not alternate, we have a chain $p_i<p_{i+1}<p_{i+2}$ (or $p_i > p_{i+1} > p_{i+2}$) and 
removing $p_{i}$ from the sequence 
yields a homotopic loop (since this chain forms a 2-simplex).  (To be precise, the subscripts
above must be taken modulo the length of the cycle; we will continue this abuse of notation 
throughout the proof.)  Repeating this process eventually yields a reduced cycle.

We can simplify any cycle even further. 
Say $C = (p_1, \ldots , p_n)$ is a reduced cycle.  We will call a vertex $p_i$ \emph{lower}
if $p_i < p_{i+1}$.  We call $C$ \emph{atomic} if every lower vertex of $C$ is an atom.  
Any reduced cycle (and hence any cycle) is homotopic to an atomic cycle, 
since we may simply replace each lower vertex $p_i$ with an atom $a_i < p_i$ (see the diagram).

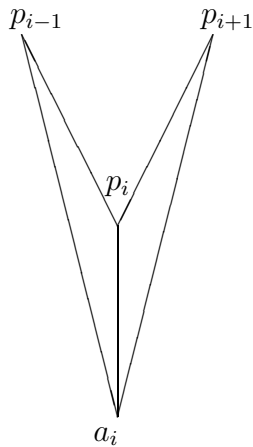
\begin{center}
\setlength{\unitlength}{3947sp}%
\begingroup\makeatletter\ifx\SetFigFont\undefined%
\gdef\SetFigFont#1#2#3#4#5{%
  \reset@font\fontsize{#1}{#2pt}%
  \fontfamily{#3}\fontseries{#4}\fontshape{#5}%
  \selectfont}%
\fi\endgroup%
\begin{picture}(2625,2760)(1051,-2911)
\thinlines
{\color[rgb]{0,0,0}\put(1801,-361){\line( 1,-4){600}}
\put(2401,-2761){\line( 0, 1){1200}}
\put(2401,-1561){\line(-1, 2){600}}
}%
{\color[rgb]{0,0,0}\put(2401,-2761){\line( 1, 4){600}}
\put(3001,-361){\line(-1,-2){600}}
}%
\put(2251,-2911){\makebox(0,0)[lb]{\smash{\SetFigFont{12}{14.4}{\rmdefault}{\mddefault}{\updefault}{\color[rgb]{0,0,0}$a_i$}%
}}}
\put(2326,-1336){\makebox(0,0)[lb]{\smash{\SetFigFont{12}{14.4}{\rmdefault}{\mddefault}{\updefault}{\color[rgb]{0,0,0}$p_i$}%
}}}
\put(2926,-286){\makebox(0,0)[lb]{\smash{\SetFigFont{12}{14.4}{\rmdefault}{\mddefault}{\updefault}{\color[rgb]{0,0,0}$p_{i+1}$}%
}}}
\put(1726,-286){\makebox(0,0)[lb]{\smash{\SetFigFont{12}{14.4}{\rmdefault}{\mddefault}{\updefault}{\color[rgb]{0,0,0}$p_{i-1}$}%
}}}
\end{picture}
\end{center}

Note that both triangles in the above diagram are 2-simplices in $\ord{P}$.

To complete the proof we must show that any atomic cycle $C$ is null-homotopic.  We will proceed
by induction on the length of $C$, which must be even (since $C$ is reduced).  The shortest 
(non-constant) atomic cycle has length 4, and can be written $C = (a_1<p>a_2<q>a_1)$.  
To show that $C$ is null-homotopic, simply consider the point $\gen{a_1, a_2}$.  
If $\gen{a_1,a_2}\in \{p, q\}$, then it is easy to check that $C$ is null-homotopic,
and otherwise the following diagram shows that $C$ is null-homotopic (again, each triangle
is a 2-simplex).

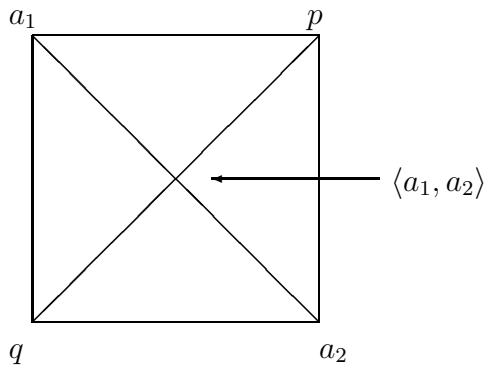
\begin{center}
\setlength{\unitlength}{3947sp}%
\begingroup\makeatletter\ifx\SetFigFont\undefined%
\gdef\SetFigFont#1#2#3#4#5{%
  \reset@font\fontsize{#1}{#2pt}%
  \fontfamily{#3}\fontseries{#4}\fontshape{#5}%
  \selectfont}%
\fi\endgroup%
\begin{picture}(2400,2235)(451,-1786)
\thinlines
{\color[rgb]{0,0,0}\put(601,239){\line( 0,-1){1800}}
\put(601,-1561){\line( 1, 1){900}}
\put(1501,-661){\line(-1, 1){900}}
}%
{\color[rgb]{0,0,0}\put(2401,239){\line(-1, 0){1800}}
\put(601,239){\line( 1,-1){900}}
\put(1501,-661){\line( 1, 1){900}}
}%
{\color[rgb]{0,0,0}\put(601,-1561){\line( 1, 0){1800}}
\put(2401,-1561){\line(-1, 1){900}}
\put(1501,-661){\line(-1,-1){900}}
}%
{\color[rgb]{0,0,0}\put(2401,239){\line( 0,-1){1800}}
\put(2401,-1561){\line(-1, 1){900}}
\put(1501,-661){\line( 1, 1){900}}
}%
{\color[rgb]{0,0,0}\put(2776,-661){\vector(-1, 0){1050}}
}%
\put(451,314){\makebox(0,0)[lb]{\smash{\SetFigFont{12}{14.4}{\rmdefault}{\mddefault}{\updefault}{\color[rgb]{0,0,0}$a_1$}%
}}}
\put(2326,314){\makebox(0,0)[lb]{\smash{\SetFigFont{12}{14.4}{\rmdefault}{\mddefault}{\updefault}{\color[rgb]{0,0,0}$p$}%
}}}
\put(451,-1786){\makebox(0,0)[lb]{\smash{\SetFigFont{12}{14.4}{\rmdefault}{\mddefault}{\updefault}{\color[rgb]{0,0,0}$q$}%
}}}
\put(2401,-1786){\makebox(0,0)[lb]{\smash{\SetFigFont{12}{14.4}{\rmdefault}{\mddefault}{\updefault}{\color[rgb]{0,0,0}$a_2$}%
}}}
\put(2851,-736){\makebox(0,0)[lb]{\smash{\SetFigFont{12}{14.4}{\rmdefault}{\mddefault}{\updefault}{\color[rgb]{0,0,0}$\gen{a_1, a_2}$}%
}}}
\end{picture}
\end{center}

Now we may assume that the length of $C$ is greater than 4, so that $C$ contains a path of the 
form $\rho = (a_1<p>a_2<q>a_3)$ (with the $a_i$ atoms).  
The diagram below shows how to construct a homotopy of $\rho$
with the path $(a_1<\gen{a_1, a_2, a_3}>a_3)$, assuming that the vertices $a_1, a_2, a_3,
p, q, \gen{a_1, a_2}, \gen{a_2, a_3}$ and $\gen{a_1, a_2, a_3}$ are all distinct.  
(Note that, by hypothesis, $\gen{a_1, a_2, a_3}\neq P$.)

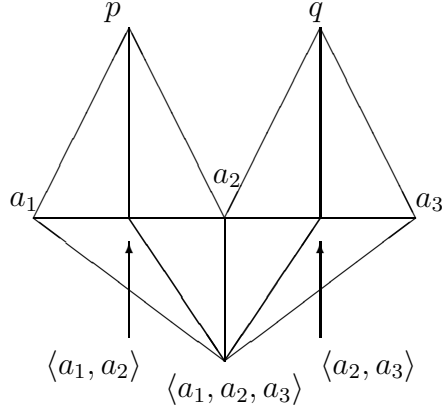
\begin{center}
\setlength{\unitlength}{3947sp}%
\begingroup\makeatletter\ifx\SetFigFont\undefined%
\gdef\SetFigFont#1#2#3#4#5{%
  \reset@font\fontsize{#1}{#2pt}%
  \fontfamily{#3}\fontseries{#4}\fontshape{#5}%
  \selectfont}%
\fi\endgroup%
\begin{picture}(2562,2598)(1051,-2740)
\thinlines
{\color[rgb]{0,0,0}\put(1201,-1561){\line( 1, 2){600}}
\put(1801,-361){\line( 0,-1){1200}}
\put(1801,-1561){\line(-1, 0){600}}
}%
{\color[rgb]{0,0,0}\put(2401,-2461){\line(-4, 3){1200}}
\put(1201,-1561){\line( 1, 0){600}}
\put(1801,-1561){\line( 2,-3){600}}
}%
{\color[rgb]{0,0,0}\put(2401,-1561){\line(-1, 2){600}}
\put(1801,-361){\line( 0,-1){1200}}
\put(1801,-1561){\line( 1, 0){600}}
}%
{\color[rgb]{0,0,0}\put(2401,-2461){\line( 0, 1){900}}
\put(2401,-1561){\line(-1, 0){600}}
\put(1801,-1561){\line( 2,-3){600}}
}%
{\color[rgb]{0,0,0}\put(3001,-1561){\line(-1, 0){600}}
\put(2401,-1561){\line( 0,-1){900}}
\put(2401,-2461){\line( 2, 3){600}}
}%
{\color[rgb]{0,0,0}\put(3601,-1561){\line(-1, 0){600}}
\put(3001,-1561){\line(-2,-3){600}}
\put(2401,-2461){\line( 4, 3){1200}}
}%
{\color[rgb]{0,0,0}\put(3001,-361){\line(-1,-2){600}}
\put(2401,-1561){\line( 1, 0){600}}
\put(3001,-1561){\line( 0, 1){1200}}
}%
{\color[rgb]{0,0,0}\put(3001,-361){\line( 1,-2){600}}
\put(3601,-1561){\line(-1, 0){600}}
\put(3001,-1561){\line( 0, 1){1200}}
}%
{\color[rgb]{0,0,0}\put(1801,-2311){\vector( 0, 1){600}}
}%
{\color[rgb]{0,0,0}\put(3001,-2311){\vector( 0, 1){600}}
}%
\put(2326,-1336){\makebox(0,0)[lb]{\smash{\SetFigFont{12}{14.4}{\rmdefault}{\mddefault}{\updefault}{\color[rgb]{0,0,0}$a_2$}%
}}}
\put(3601,-1486){\makebox(0,0)[lb]{\smash{\SetFigFont{12}{14.4}{\rmdefault}{\mddefault}{\updefault}{\color[rgb]{0,0,0}$a_3$}%
}}}
\put(1651,-286){\makebox(0,0)[lb]{\smash{\SetFigFont{12}{14.4}{\rmdefault}{\mddefault}{\updefault}{\color[rgb]{0,0,0}$p$}%
}}}
\put(2926,-286){\makebox(0,0)[lb]{\smash{\SetFigFont{12}{14.4}{\rmdefault}{\mddefault}{\updefault}{\color[rgb]{0,0,0}$q$}%
}}}
\put(1051,-1486){\makebox(0,0)[lb]{\smash{\SetFigFont{12}{14.4}{\rmdefault}{\mddefault}{\updefault}{\color[rgb]{0,0,0}$a_1$}%
}}}
\put(3001,-2536){\makebox(0,0)[lb]{\smash{\SetFigFont{12}{14.4}{\rmdefault}{\mddefault}{\updefault}{\color[rgb]{0,0,0}$\gen{a_2,a_3}$}%
}}}
\put(1276,-2536){\makebox(0,0)[lb]{\smash{\SetFigFont{12}{14.4}{\rmdefault}{\mddefault}{\updefault}{\color[rgb]{0,0,0}$\gen{a_1,a_2}$}%
}}}
\put(2026,-2686){\makebox(0,0)[lb]{\smash{\SetFigFont{12}{14.4}{\rmdefault}{\mddefault}{\updefault}{\color[rgb]{0,0,0}$\gen{a_1,a_2,a_3}$}%
}}}
\end{picture}
\end{center}

If
$\gen{a_1, a_2} = p$ or $\gen{a_2, a_3} = q$, there is no real problem.  Otherwise, it is easy
to check that if two of these vertices were equal, we would have either $a_1<q$ or, symmetrically, 
$a_3<p$.  If $a_1<q$ then the cycle $(a_1<p>a_2<q>a_1)$ is null-homotopic 
(since it is an atomic cycle of length four) and hence its sub-paths
$(a_1<p>a_2<q)$ and $(a_1<q)$ are homotopic.  This implies that $\rho$ is homotopic to the path
$(a_1<q>a_3)$.  

In any case, we see that $\rho$ is homotopic (relative to its endpoints) to a shorter
path, and hence the cycle $C$ is homotopic to a shorter cycle.  In addition, this shorter cycle
is still atomic, and repeating the process will eventually provide a null-homotopy of $C$.
$\hfill \Box$ 

\vspace{.15in}
We now present a result guaranteeing simple-connectivity of $\C{G}$ 
under weaker conditions than those of
Corollary~\ref{2-generator}.  First we need a simple lemma 
(Brown~\cite[Proposition 14]{Brown-coset-poset}).

\begin{lemma}$\label{coset-poset-connectivity}$
$\C{G}$ is connected unless $G$ is cyclic of prime-power order.
\end{lemma}
{\bf Proof.}  If $G$ is non-cylic, then any proper coset $xH$ connects to the identity via
the path $xH \geqs \{x\} \leqs \gen{x} \geqs \{1\}$.  Now say $G$ is cyclic but not of prime
power order.  If $x$ does not generate $G$, then any coset $xH$ connects to the indentity
as before.  If $x$ is a generator, then we may write $x=yz$ where $y$ and $z$ are not
generators (this is possible because $G$ decomposes as a direct product in a non-trivial manner).
Then any coset $xH$ connects to $\{1\}$ via the path
$$xH\geqs \{x = yz\} \leqs y\gen{z}\geqs \{y\}\leqs \gen{y} \geqs \{1\}.$$

If $G$ is cyclic of prime-power order, then $G$ has a unique maximal subgroup $M$, 
and $\C{G}$ is the disjoint union of the cones $\C{G}_{\leqs xM}$.  If $G\isom \Z/p^n$,
then $(G:M) = p$ and we see that $\C{G}$ is homotopy equivalent to a discrete space
with $p$ points.
$\hfill \Box$

\begin{theorem}$\label{2-gen-presentation}$
Let $G$ be a finite, non-cyclic group in which the following condition holds: 
For any $x,y\in G$ s.t. 
$\langle x,y \rangle = G$ there exists $z\in G$ s.t. $\langle z,x \rangle \neq G$,
$\langle z,y \rangle \neq G$ and $\langle zx,zy \rangle \neq G$.  Then $\C{G}$ is simply
connected.
\end{theorem}
{\bf Proof.}  Let $G$ be a group satisfying the hypotheses of the theorem.  We will show
that $\pi_1 (\M{G}) = 0$.  First, consider the 1-skeleton of $\M{G}$.  Since $G$ is non-cyclic,
any two elements $x,y$ are contained in a proper coset, namely $x\langle x^{-1} y \rangle$.
Thus in $\M{G}$ all possible edges between vertices exist.

We wish to apply Theorem~\ref{presentation} to $\M{G}$, so we must choose a maximal tree
in $\M{G}^1$.
Since all possible edges exist, we can choose our maximal tree $T$ to be the star at
the vertex $1$, i.e. $T$ consists of all edges of the form $\{1,x\}$ with $x\in G$.  In order
to show that $\pi_1 (\M{G})$ is trivial, we will use the relations given in 
Theorem~\ref{presentation} to show that each generator is trivial.

There are two cases.  First, if $x$ and $y$ do not generate $G$, then the subgroup $\gen{x,y}$
is a coset containing $1,x$ and $y$, and thus these vertices form a simplex in $\M{G}$.
So we obtain the relation $(x,y)(y,1)(1,x)=1$, which forces $(x,y)$ to be trivial.

Next, consider the case in which $\gen{x,y} = G$.  Let $z$ be the element guaranteed by the
hypotheses.  Then, as above, the generators corresponding to the edges $\{x,z^{-1}\}$ and
$\{y,z^{-1}\}$ are trivial (note that $\gen{x,z}\neq G \iff \gen{x,z^{-1}} \neq G$).  
Furthermore, since $\gen{zx, zy} \neq G$ there is a simplex $\{x,y,z^{-1}\}$ in $\M{G}$, and
we have the relation $(x,y)(y,z^{-1})(z^{-1},x) = 1$, which again forces $(x,y)$ to be trivial.  


Thus all generators for $\pi_1 (\M{G})$ are trivial, and hence $\M{G}$ and $\C{G}$ are simply
connected.
$\hfill \Box$

\vspace{.15in} In the next chapter we will apply this result to show that $\C{A_5}$ 
is simply connected.


\section{Mayer-Vietoris Sequences and $H_* (\C{G})$}
$\label{m-v}$

We will now present some results based on Mayer-Vietoris sequences for the homology of 
$\C{G}$.

From here on, we let $H_n(X)$ denote the homology of $X$ with coefficients in $\Z$,
and we let $\tilde{H}_n (X)$ denote reduced homology (again with coefficients in $\Z$).

\begin{theorem}$\label{M-V1}$
Let $G$ be a finite non-cyclic group with a cyclic maximal subgroup $M$
of prime-power order.
Let $o(M) = p^n$.  Then $H_1(\C{G})$ has rank at least $(p-1) (G:M)$, and in particular 
$\C{G}$ is not simply connected.
\end{theorem}
{\bf Proof.}  Let $\{x_i\}_{i\in I}$ be a set of (left) coset representatives for $M$.
Let $X = \ord{\C{G}}$, $Y = \ord{\C{G} - \{x_i M\}_{i\in I}}$, 
and $Z = \ord{\cup_{i\in I} \C{G}_{\leqs x_i M}}$.
Then we have $X = Y\cup Z$, since any chain in $\C{G}$ either contains $x_i M$ for some $i\in I$,
and lies in $Z$, or does not contain any coset of $M$, and lies in $Y$.  
Thus we have a Mayer-Vietoris sequence
$$\cdots \maps \tilde{H}_1(X) 
         \stackrel{f}{\maps} \tilde{H}_0 (Y\cap Z) 
         \stackrel{h}{\maps} \tilde{H}_0 (Y) \oplus \tilde{H}_0 (Z) 
         \maps \tilde{H}_0 (X)$$
on (reduced) homology groups.  

Note that $X$ is connected because $G$ is non-cyclic (Lemma~\ref{coset-poset-connectivity}),
and hence $\tilde{H}_0 (X) = 0$.
Also, $Z$ is a union of $(G:M)$ cones, one for each $i\in I$, and hence 
$\tilde{H}_0 (Z) \homeo \Z^{(G:M) - 1}$.  
In addition, we claim that $Y$ is connected so that $\tilde{H}_0(Y) = 0$.  Substituting these
values into the above sequence yields an exact sequence
$$\cdots \maps \tilde{H}_1(X) \stackrel{f}{\maps} \tilde{H}_0 (Y\cap Z) 
         \stackrel{h}{\maps} \Z^{(G:M) - 1} \maps 0,$$
and thus $h:\tilde{H}_0(Y\cap Z) \maps \Z^{(G:M) - 1}$ is a surjection.  Now, 
$$Y\cap Z = \ord{\cup_{i\in I} \C{G}_{<x_i M}} \homeo \coprod_{i\in I} \C{M}.$$  
As discussed in the proof of Lemma~\ref{coset-poset-connectivity}, the coset
poset of $M\homeo \Z/p^n$ is disconnected, and in fact it has $p$ contractible components 
(the cones under the cosets of the unique maximal subgroup).  
Thus $\tilde{H}_0 (Y\cap Z) \homeo \Z^{p (G:M) - 1}$, and
thus Im$\, f = $ Ker$\, h \homeo \Z^{p(G:M) - 1 - ((G:M) - 1)} = \Z^{(p-1)(G:M)}$ 
and we see that the rank of $H_1(X)$ is at least $(p-1)(G:M)$, as desired.

To complete the proof we must show that $Y$ is connected.
Consider a coset $xH$, where
$H\neq M$.  If $\gen{x}\neq M$, then we have a path 
$$xH\geqs \{x\} \leqs \gen{x} \geqs \{1\}$$
connecting $xH$ to the identity (note that we have assumed $G$ is non-cyclic, so
that $\gen{x}\in Y$).  Next, if $\gen{x} = M$ then choose some $g\in G$, $g\notin M$.
We now have a path
$$xH\geqs \{x\} \leqs x \gen{g} \geqs \{xg\} \leqs \gen{xg} \geqs \{1\}$$
connecting $xH$ to the identity (note that $\gen{xg} \neq M$, since if $xg\in M$ then
$g\in x^{-1}M = M$, contrary to our choice of $g$).
$\hfill \Box$

\vspace{.15in}
The above result does not extend to cyclic groups, as we will see in a moment.  
The basic problem is that the subcomplex $Y$ is no longer connected.  Nevertheless,
the coset poset of a cyclic group with maximal subgroup isomorphic to $\Z/p^n$ is still
not simply connected.  Any group of this type may be written as $\Z/q \cross \Z/p^n$, with
$q$ prime, and now~\cite[Lemma 5]{Brown-coset-poset} 
implies that $\C{G} \homeo \abs{\C{\Z/q}} \ast \abs{\C{\Z/p^n}}$
(where $*$ denotes the join).
The latter is (homotopy equivalent to) a bipartite graph with $q + p$ vertices and $pq$ edges.  
Modding out a maximal tree in this graph and applying Van Kampen's Theorem shows that 
$H_1(\C{G})\homeo \Z^{(p-1)(q-1)}$.  (For $G$ to satisfy the conclusion of the theorem,
we would need the rank of $H_1$ to be at least $(p-1)q$.)

\begin{question} Does Theorem~\ref{M-V1} characterize (non-cyclic) 
finite groups with simply connected
coset posets?
\end{question}

We expect the answer to be no because, as we will show later in this
section, Theorem~\ref{M-V1} applies only to a small class of solvable groups.  In addition,
Corollary~\ref{M-V2} may allow one to find other finite groups with non-simply connected
coset posets. 

\vspace{.15in} 

The following theorem shows a (somewhat weak) connection between the coset poset and
the subgroup poset of a finite group.  In Chapter~\ref{cp-sp}, we will discuss
other possible connections.

\begin{theorem}$\label{homology-surjection}$
Let $G$ be a finite group, and for $g\in G$, let $\C{G}_g$ denote
the poset $\C{G} - \{g\}$.

If there exists $g\in G$ with $\tilde{H}_n (\C{G}_g) = 0$,
then there is a surjection $\tilde{H}_{n+1} (\C{G}) \to \tilde{H}_n (S(G))$.
\end{theorem}

Note that for any $g, h\in G$, the posets $\C{G}_g$ and $\C{G}_h$ are isomorphic
(consider the map given by left-multiplication by $hg^{-1}$).  In particular,
$\C{G}_g \isom \C{G}_1 \isom S(G)$.

\vspace{.15in}
{\bf Proof.}  Let $X = \ord{\C{G}}$, $Y = \ord{\C{G}_{\geqs \{g\}}}$, and 
$Z =  \ord{\C{G}_g}$.  
Note that $Y$ is contractible (so $\tilde{H}_n (Y) = 0$ for all $n$) and 
that $Y\cap Z \isom S(G)$.

Since $X = Y\cup Z$, we get a Mayer-Vietoris sequence
$$\cdots \maps \tilde{H}_{n+1} (X) \stackrel{f}{\maps} \tilde{H}_n (Y\cap Z) 
     \maps \tilde{H}_n (Y)\oplus \tilde{H}_n (Z) \maps H_n (X) \maps\cdots
$$
and since $\tilde{H}_n (Y) = \tilde{H}_n (Z) = 0$, the map 
$$f:H_{n+1} (X)\to H_n (Y\cap Z) \isom H_n (S(G))$$
is a surjection.
$\hfill \Box$

\begin{corollary}$\label{M-V2}$
Let $G$ be a finite group and let $n$ be the number of path-components of $S(G)$.
Then $H_1(\C{G})$ has rank at least $n-1$.  
In particular, if $S(G)$ is disconnected then $\C{G}$ is not simply connected.
\end{corollary}
{\bf Proof.}  First, if $G$ is cyclic of prime-power order then $\C{G}$ is homotopy
equivalent to a discrete set of $p$ points, and hence $H_1(\C{G}) = 0$.  Also,
$G$ has a unique maximal subgroup in this case, and hence $S(G)$ has one path component.

Now we turn to the case in which $G$ is not cyclic of prime-power order.
In light of Theorem~\ref{homology-surjection} it will suffice to prove that
$\C{G}_1$ is connected.    
First, consider the case in which $G$ is not cyclic.
Choose some element $x\in G$, $x\neq 1$.  Then any vertex $yH$
in $Y$ can be connected to $\{x\}$ via the path 
$$yH \geqs \{y\} \leqs x\gen{x^{-1} y}\geqs \{x\}$$
(note that $G$ is non-cyclic so $x\gen{x^{-1} y}\in Y$, and $yH\neq \{1\}$ so we may assume 
$y\neq 1$).  

If $G$ is cyclic (but not of prime-power order), let $x$ be 
a generator for $G$.  Then any coset in $G$ can be written as $x^n\gen{x^m}$, for some
$n,m>0$.  Also, note that since $o(G)$ has at least two prime factors, $G$ decomposes
as a (non-trivial) direct product and hence any element in $G$ can be written as a product of two
non-generators.  Specifically, we may choose $k,l$ s.t. $x^{n-1} = x^kx^l$ and $\gen{x^k}$,
$\gen{x^l}$ are proper subgroups of $G$.  Using this decomposition, we have the following
path joining $x^n\gen{x^m}$ to $\{x\}$:
$$   	x^n\gen{x^m}
\geqs 	\{x^n = x^{k+1}x^l\}
\leqs 	x^{k+1}\gen{x^l}
\geqs 	\{x^{k+1}\}
\leqs 	x\gen{x^k}
\geqs	\{x\}.
$$
Thus $Y$ is connected and the proof is complete.
$\hfill \Box$

\vspace{.15in} Corollary~\ref{M-V2} does not characterize finite groups with simply connected
coset posets.  The quaternion group $Q_8$ provides a counter-example, since $S(Q_8)$ is
contractible (it is a cone on $\gen{-1}$) and yet 
$\pi_1 \left(\C{G}(Q_8)\right) \homeo \Z \ast \Z \ast \Z$, as will be shown
in the next chapter (Section~\ref{Q_8}).

\begin{question} For which finite groups $G$ is $S(G)$ disconnected?  In particular,
do there exist non-solvable groups with disconnected subgroup lattices?
\end{question}

We now turn to the question of which groups satisfy the conditions of Theorem~\ref{M-V1}.
In the case where $G$ is a $p$-group, these groups have been classified 
(see Brown~\cite[Theorem 4.1, p. 98]{Brown-cohomology}).

We now turn to the case in which $o(G)$ has at least two prime divisors, $p$ and $q$.
First, the following result can be found in Herstein~\cite{Herstein-remark}:

\begin{theorem}$\label{herstein}$  Let $G$ be a finite group with an abelian maximal subgroup.  
Then $G$ is solvable.
\end{theorem}

In our case, something stronger may be shown.  Our proof of the following
result is based on Herstein's proof of Theorem~\ref{herstein}.

\begin{theorem} $\label{max-abelian-p}$
Let $G$ be a finite group and assume $p$ and $q$ are distinct primes 
dividing $o(G)$.  Assume further that $G$ has a maximal subgroup $M$ which is an 
abelian $p$-group.  Then either $G\homeo M\rtimes \Z/q$ or $G \homeo Q\rtimes M$,
where $Q$ is the (unique) $q$-Sylow subgroup of $G$.
\end{theorem}

Note that solvability of $G$ follows immediately from this theorem, since
$Q$ and $M$ are each nilpotent.  

In proving the theorem, we will need a result due to Burnside 
(see~\cite[p. 289]{Robinson-theory-of-groups}) and a lemma from 
Herstein~\cite{Herstein-remark}.  
Herstein's proof is elementary, and will be included for completeness.

\begin{theorem}[Burnside's Theorem]$\label{burnside}$ 
Let $G$ be a finite group and let $P$ be a $p$-Sylow subgroup of $G$ which lies
in the center of its normalizer.  Then $G$ is $p$-nilpotent, that is there exists
a subgroup $T\normal G$ s.t. $T\cap P = \{1\}$ and $TP = G$ (hence $G\homeo T\rtimes P$).
\end{theorem}

\begin{lemma} Let $\A$ be an abelian group of automorphisms of a finite group
$G$.  Assume that each non-identity element $\alpha \in \A$ fixes only the identity in $G$.
Then for each prime $p$ dividing $o(G)$, there exists a unique Sylow $p$-subgroup
$S_p$ of $G$ which is invariant under the action of $\A$.
\end{lemma}
{\bf Proof.}  Let $\alpha \in A$ be a non-trivial automorphism of $G$.  Consider the set
$\{g^{-1}\alpha (g):g\in G\}$.  We claim that this set is all of $G$.  Indeed, if
$g^{-1}\alpha (g) = h^{-1}\alpha (h)$ then $hg^{-1} = \alpha (hg^{-1})$ and thus 
$hg^{-1}$ is fixed
by $\alpha$.  By hypothesis we now have $hg^{-1} = 1$, and thus $h = g$.  

Now, fix a prime $p$ dividing $o(G)$ and a Sylow $p$-subgroup $P$.  Then 
$\alpha (P) = x P x^{-1}$ for some $x\in G$, and (by the discussion above)
there exists $y\in G$ s.t.
$y^{-1}\alpha (y) = x^{-1}$.  So we have
$$\alpha (y P y^{-1}) = \alpha (y) \alpha (P) \alpha (y^{-1}) 
  = \alpha (y) x P x^{-1} \alpha (y^{-1}) = y P y^{-1}.$$ 
Thus $y P y^{-1}$ is a Sylow $p$-subgroup fixed by $\alpha$.  We claim that in fact, 
this is the only Sylow $p$-subgroup fixed by $\alpha$.  Assuming uniqueness, we now
complete the proof.  Letting $S_p = y P y^{-1}$, we see that for any $\beta \in \A$
$$\alpha \beta (S_p) = \beta \alpha (S_p) = \beta (S_p).$$
(Recall that $\A$ is abelian.)
Thus $\beta (S_p)$ is also fixed by $\alpha$ and must be $S_p$.  So $S_p$ is fixed by all of 
$\A$.

We will now prove uniqueness of $S_p$.  If $S'_p$ is another $p$-Sylow subgroup fixed 
by $\alpha$, we have $g S_p g^{-1} = S'_p$ for
some $g\in G$, and $g^{-1} \alpha (g)\in N(S_p)$ because
$$g^{-1} \alpha (g) S_p \alpha (g^{-1}) g = g^{-1} \alpha (g S_p g^{-1}) g 
  = g^{-1} \alpha (S'_p) g = g^{-1} S'_p g = S_p.$$
Now, note that $\alpha$ induces an automorphism of $N(S_p)$ which fixes only the identity,
and thus all elements of $N(S_p)$ have the form $n^{-1} \alpha (n)$ for some $n\in N(S_p)$.
Since $g^{-1} \alpha(g) \in N(S_p)$, this means $g^{-1} \alpha (g) = n^{-1} \alpha (n)$ for
some $n\in N(S_p)$, and hence $g = n$ (as was shown at the start of the proof).  
So $g\in N(S_p)$, which means $S'_p = S_p$.  This
shows that $S_p$ is the unique Sylow $p$-subgroup fixed by $\alpha$, completing the proof.
$\hfill \Box$

\vspace{.15in}
{\bf Proof of Theorem~\ref{max-abelian-p}.}
Let $G$ be a finite group satisfying the hypotheses of the theorem.  If $M\normal G$,
then $G/M$ has no non-trivial subgroups, and hence $G/M\homeo \Z/q$.  Thus $o(G) = q \, o(M)$,
and hence the $q$-Sylow subgroups of $G$ are isomorphic to $\Z/q$ and we have
$G\homeo M \rtimes \Z/q$.

Next, say $M\nnormal G$.  We will proceed by induction on $o(G)$.  Since $p,q|o(G)$,
the base case occurs when $o(G) = pq$.  In this case $G\homeo \Z/q\rtimes \Z/p$, since
$M\nnormal G$ (this follows from an easy argument using the Sylow theorems.)  We now
assume that for any group $H$ with $o(H)<o(G)$, $p,q|o(H)$ 
and with an abelian maximal subgroup $P$ 
which is a $p$-group, $P\nnormal H$ implies $H\homeo Q\rtimes P$ (where $Q$ is the Sylow
$q$-subgroup of $H$).

Now, since $M$ is maximal, $N(M) = M$ and since $M$ is abelian
Burnside's Theorem applies.  So there is a subgroup $T\normal G$ with $G\homeo T\rtimes M$,
and it remains to show that $T$ is a $q$-group.

We will consider two cases.  First we assume that for any $g\in G$, $g\notin M$,
$(g^{-1}Mg)\cap M = \{1\}$.
Consider the action of $M$ on $T$ by conjugation.  We claim that for any $m\in M$ ($m\neq 1$), 
the only element of $T$ fixed under conjugation by $m$ is 1.  Indeed, if $t\in T$ and 
$m^{-1}tm = t$, then $tmt^{-1} = m \in M\cap (tMt^{-1})$.  Since $m\neq 1$, 
$M\cap (tMt^{-1}) \neq \{1\}$ and hence $t\in M$, which implies $t = 1$ (since $t\in T$ 
and $T\cap M = \{1\}$).  Thus we may apply
Herstein's lemma to the (abelian) group of automorphisms of $T$ induced by $M$.  This
gives us a unique Sylow $q$-subgroup $Q$ of $T$ normalized by $M$ (note that $q|o(T)$
because $q|o(G)$ and $G\isom T\rtimes M$).  
Now we have $Q, M\subset N(Q)$, so $Q\normal G$ by maximality of $M$.  Hence $QM$ is a 
subgroup of $G$, and since $M$ is maximal we have $QM = G$.  This shows that 
$G\homeo Q\rtimes M$.

Finally we must consider the case in which $M\nnormal G$ and there exists $g\in G$, $g\notin M$,
with $(g^{-1}Mg) \cap M \neq \{1\}$.  Let $W = (g^{-1}Mg) \cap M$, and choose $w\in W$, 
$w\neq 1$.  Note that $W< M$, since $N(M) = M$ and $g\notin M$.
We now have $\gen{w} \normal M, g^{-1}Mg$ (since these groups are abelian), 
which implies $\gen{w}\normal G$ (by maximality of $M$).  

Consider $G/\gen{w}$.  Since $\gen{w}\leqs W< M$, $M/\gen{w}$ is an abelian
$p$-group, and is maximal in $G/\gen{w}$.  Also, we have $p,q|o(G/\gen{w})$, and since
$M\nnormal G$, $M/\gen{w}\nnormal G/\gen{w}$.  So
$G/\gen{w}$ is a group of smaller order satisfying the conditions of the induction hypothesis,
and we may assume that $G/\gen{w}\homeo \bar{Q} \rtimes M/\gen{w}$, where $\bar{Q}$
is the Sylow $q$-subgroup of $G/\gen{w}$.  Notice that
$p$ and $q$ are the only prime factors of $o(G/\gen{w})$, and since $\gen{w}$ is a $p$-group, 
these must be the only prime factors of $o(G)$ as well.
Finally, since $M$ is maximal it must be a Sylow $p$-subgroup of $G$, and hence 
$T$ is a $q$-group, as desired.
$\hfill \Box$

\begin{corollary} Let $G$ be a group with maximal subgroup $M\homeo \Z/p^n$, for some prime
$p$ and some $n>0$.  If $G$ is not a $p$-group, then there is a prime $q\neq p$ such that
either $G\isom \Z/p^n \rtimes \Z/q$ or $G\homeo Q\rtimes \Z/p^n$ (with $Q$ a $q$-group).
\end{corollary}

\section{Applications to Group Theory}

We now present group theoretic applications of the results in the previous section.  First
we will need to discuss a result of Brown~\cite[Proposition 11]{Brown-coset-poset}, 
describing the homotopy type of the coset poset of a solvable group.

\begin{definition}  Let $G$ be a finite group, and let $P$ be the poset of normal subgroups
of $G$ (inlcuding $\{1\}$ and $G$), ordered by inclusion.  
The maximal chains in $P$ are called \emph{chief series} for $G$.

If $\{1\} = G_0 \normal G_1 \normal \cdots \normal G_k = G$ is a chief series for $G$,
we call an element $G_i$, $i = 1,\ldots k$, 
\emph{complemented} if $G_i/G_{i-1}$ has a complement in 
$G/G_{i-1}$.
\end{definition}

\begin{theorem}$\label{solvable-coset-poset}$
Let $G$ be a finite solvable group, and let 
$$\{1\} = G_0 \normal G_1 \normal \cdots \normal G_k = G$$
be a chief series for $G$.  Then $\C{G}$ is homotopy equivalent to a bouquet of 
spheres, each of dimension $d-1$ where $d$ is the number of indices $\, i$, 
$1\leqs i \leqs k$, s.t. $G_i$ is complemented in the above chief series.
\end{theorem}

\begin{corollary}$\label{chief-series1}$
Let $G$ be a finite group with maximal subgroup $M \homeo \Z/p^n$,
for some prime $p$ and some $n>0$.  Then in any chief series for $G$, exactly one 
proper subgroup is complemented.
\end{corollary}  
{\bf Proof.}
The results of the previous section imply that $G$ is solvable,
so Theorem~\ref{solvable-coset-poset} implies that $\C{G}$ has the homotopy type of a 
bouquet of spheres.  Since $\C{G}$ is not simply connected (by Theorem~\ref{M-V1})
the dimension of these spheres must be 1.  So $d=2$, and there is exactly one complemented
factor (other than $G$ itself) in any chief series for $G$.
$\hfill \Box$

\vspace{.15in}
Similarly, we have:

\begin{corollary}$\label{chief-series2}$
Let $G$ be a solvable finite group with disconnected subgroup poset.
Then in any chief series for $G$, exactly one proper subgroup is complemented.
\end{corollary}

\section{Group Extensions}$\label{extensions}$

In this section we will consider the coset poset of a group $G$ which may be
written as an extension
$$N\injects G \maps G/N.$$
First, we characterize (non-trivial) direct products with simply connected coset posets
(Corollary~\ref{direct-products}). 
The proof employs a result from Brown~\cite[Lemma 5]{Brown-coset-poset} on direct products, 
which we generalize to semi-direct products
(Theorem~\ref{semi-direct-products}).  Lastly, we employ another result 
from~\cite{Brown-coset-poset} to show that if a group $G$ has a quotient $G/N$ 
whose coset poset is simply connected, then $\C{G}$ is simply connected as well.
These results show, in particular,
that there are infinitely many non-solvable groups with simply-connected coset posets.

\begin{definition} We call a coset $C\in \C{G\cross H}$ \emph{saturating} if it surjects
onto each factor under the natural projections from $G\cross H$ to $G$ and $H$.  The poset of
non-saturating cosets of $G\cross H$ will be denoted $\mC_0 (G\cross H)$.
\end{definition}

\begin{theorem}[Brown] 
$\label{brown-direct-products}$
For any finite groups $G$ and $H$, 
$$\mC_0 (G\cross H)\heq \C{G}*\C{H},$$
where $*$ denotes the join operation.
\end{theorem}

For the purposes of this section, we simply regard $P*Q$ as the topological space
$\abs{P} * \abs{Q}$ (for any posets $P$ and $Q$).  
For an alternate (and compatible) viewpoint, see Quillen~\cite{Quillen-p-subgroups}.

\begin{corollary}$\label{direct-products}$  Let $G$ and $H$ be non-trivial finite groups.
Then $\C{G\cross H}$ is simply connected if and only if 
at least one of these groups is not cyclic of prime-power order.  
\end{corollary}
{\bf Proof.}  If both groups are cyclic of prime-power order, it is easy to 
check that there are just two split factors in any chief series for $G\cross H$
(in the sense of Theorem~\ref{solvable-coset-poset}) and the desired result then 
follows from that theorem.

In the other direction, assume WLOG that $G$ is not cyclic of prime-power order
(so $\C{G}$ is connected).  Then $\mC_0 (G\cross H)$ is simply connected: by 
Theorem~\ref{brown-direct-products}, 
this space is (up to homotopy) the join of $\C{H}$ with the connected space $\C{G}$,
and the join of a non-empty space with a connected space is simply connected 
(see Milnor~\cite{Milnor-univ-bundles-2}).

Now we consider the effect of adding saturating cosets to $\mC_0 (G\cross H)$.  Consider
a minimal element $xK$ of the poset $\C{G\cross H} - \mC_0 (G\cross H)$.
Adding this to $\mC_0 (G\cross H)$ has the effect of coning off a copy of 
$\C{K}\isom \C{G}_{\leqs xK}$, which we
claim is connected.  Connectivity of $\C{K}$ results from the fact that $K$ surjects onto 
$G$ and $H$: since one of these is not
cyclic of prime-power order, neither is $K$ and hence $\C{K}$ is connected (by 
Lemma~\ref{coset-poset-connectivity}).  Now, coning
off a connected subset of a simply connected space yields another simply connected space (this
follows from Van Kampen's Theorem, or can be proven directly).

We may now continue in this manner, adding in minimal elements from the poset of remaining cosets,
and at each stage we simply cone off a copy of $\C{K}$, where $K\leqs G\cross H$ surjects
onto $G$ and onto $H$.  Thus we see that $\C{G\cross H}$ is simply connected.
$\hfill \Box$

\vspace{.15in}
We will now generalize Theorem~\ref{brown-direct-products} 
to the case of a semi-direct product.

For any semi-direct product $G\rtimes H$ (with $G$ and $H$ finite),
let $f:G\rtimes H \to G$ be the function $f(g,h) = g$, and let $\pi: G\rtimes H\to H$
be the quotient map.  We call a coset $xT\in \C{G\rtimes H}$ \emph{saturating} if
$\pi (xT) = H$ and the only $H$-invariant subgroup of $G$ that contains $f(T)$ is $G$
itself.

\begin{theorem}$\label{semi-direct-products}$
Let $G$ and $H$ be finite groups and consider a semi-direct product $G\rtimes H$.
Let $\mC_0 (G\rtimes H)$ be the poset of all non-saturating cosets and let $\mC_H (G)$ 
denote the poset of all cosets of proper, $H$-invariant subgroups of $G$.  
Then $\mC_0 (G\rtimes H)$ is homotopy equivalent to the join $\mC_H (G)*\C{H}$.
\end{theorem}
{\bf Proof.}    Let
$\mC^+ (H)$ denote the set of all cosets in $H$ (including $H$ itself) and let
$\mC_H^+ (G) = \mC_H (G) \cup \{G\}$.  Then
if $\mC_{00} (G\rtimes H)$ denotes the set of all proper cosets of the form 
$(g,h) I\rtimes K$ (with $I\leqs G$ invariant and $K\leqs H$), one checks that the
map
$$(x,y) I\rtimes K \mapsto (xI,yK)$$
is a well-defined poset isomorphism 
$$\mC_{00} (G\rtimes H) \stackrel{\isom}{\maps} \mC_H^+ (G) \cross \mC^+ (H) - \{(G,H)\}.$$
The latter has geometric realization homotopy equivalent to $\mC_H (G) * \C{H}$ 
(see Quillen~\cite[Proposition 1.9]{Quillen-p-subgroups}).
Finally, we have an increasing poset map $\Phi$ from $\mC_0 (G\rtimes H)$ onto 
$\mC_{00} (G\rtimes H)$ given by
$$\Phi((x,y) T) = (x,y) \hat{f}(T)\rtimes \pi (T)$$
where $\hat{f} (T)$ is the smallest invariant subgroup containing $f(T)$ (i.e.
$\hat{f} (T)$ is the intersection of all invariant subgroups containing $f(T)$).
Corollary~\ref{order-homotopic2} shows that this map is a homotopy equivalence
between $\mC_{0} (G\rtimes H)$ and 
$\Phi (\mC_{0} (G\rtimes H)) = \mC_{00} (G\rtimes H)$
(with homotopy inverse the inclusion map $i:\mC_{00} (G\rtimes H)\to \mC_{0} (G\rtimes H))$.
$\hfill \Box$

\vspace{.15in}  With this theorem in hand, one ought to be able to show that many
 semi-direct products have simply connected coset posets.  In light of the results
of Section~\ref{m-v}, the collection
will not be quite the same as in the case of direct products, but hopefully it will 
be relatively large.

We now turn to the case of a general group extension, and prove a simple corollary to
the following result of Brown~\cite[Proposition 10]{Brown-coset-poset}.

\begin{theorem}$\label{brown-extension}$
For any finite group $G$ and any normal subgroup $N\normal G$ there is a homotopy
equivalence
$$\C{G} \maps \C{G/N}*\C{G,N},$$
where $\C{G,N}$ denotes the poset of saturating cosets, i.e. cosets $C\in \C{G}$
s.t. $q(C) = G/N$.
\end{theorem}

\begin{corollary} Let $G$ be a finite group with quotient $H$, and assume that
$\C{H}$ is simply connected.  Then $\C{G}$ is simply connected as well.
\end{corollary}
{\bf Proof.}  Choose $N\normal G$ s.t. $G/N\isom H$.  Then Theorem~\ref{brown-extension}
tells us that $\C{G}\heq \C{H}*\C{G,N}$, and if $\C{H}$ is simply connected
then its join with any space is simply connected.
$\hfill \Box$

\vspace{.15in}  Presumably, it would be easy to extend the results of this section to
deal with higher connectivity, and it would be interesting to explore this further.

\chapter{Examples}$\label{examples}$

In this chapter, we will examine the homotopy type of $\C{G}$ for some specific finite groups
$G$.  Along the way, we will introduce new techniques for dealing with the coset poset.

The main goal of this chapter is to study the coset posets of two finite simple groups,
$A_5$ and $\psls$.  Brown's theorem (Theorem~\ref{solvable-coset-poset}),
gives a complete description of the homotopy type of the coset poset of any finite solvable
group, and but little is known about non-solvable groups.  Thus it makes sense
to study simple groups, at the opposite end of the spectrum.  In particular, we will
give two simple proofs of an unpublished result of Shareshian 
(see~\cite[p. 1009]{Brown-coset-poset}) 
describing the homotopy type of $\C{A_5}$, and we will show that the coset poset of 
$\psls$ is simply connected.

Before moving on to these more complicated results, we present another description of the
coset poset and use this description to calculate the homotopy type of $\C{Q_8}$, where $Q_8$
is the quaternion group.  This will show that Corollary~\ref{M-V2} does
not characterize groups with simply connected coset posets.

\section{The Quaternion Group}$\label{Q_8}$

Let $Q_8$ denote the standard quaternion group consisting of the 8 elements
$\pm 1$, $\pm i$, $\pm j$, $\pm k$.  The proper, non-trivial subgroups of $Q_8$ are just
$\gen{-1} = \{\pm 1\}$, $\gen{i} = \{\pm 1, \pm i\}$, $\gen{j} = \{\pm 1, \pm j\}$,
and $\gen{k} = \{\pm 1, \pm k\}$.  Thus $S(Q_8)$ is a cone with minimum element $\gen{-1}$,
and hence this poset is contractible (so Corollary~\ref{M-V2} does not apply to it).

We will now calculate the homotopy type of $\C{Q_8}$, and show that this poset is not 
simply connected.  To do this, we use the 
Cross-Cut Theorem~\cite[Theorem 10.8]{Bjorner-top-methods}.

\begin{definition} Let $P$ be a poset, and let $C\subset P$ be a subposet satisfying:
\begin{enumerate}
	\item C is an antichain (no two elements of C are comparable),
	\item If $\sigma$ is a chain in $P$ then there is some $c\in C$ which is comparable
		to all elements of $\sigma$,
	\item If $A\subset C$ has an upper (lower) bound in $P$, then $A$ has a meet 
		(join) in $P$.
\end{enumerate}

We call $C$ a cross-cut, and we define $\Gamma(P,C)$, the 
cross-cut complex associated to $C$, to be the simplicial complex with vertex set $C$
and with a simplex for each subset of $C$ which is bounded in $P$.
\end{definition}

\begin{lemma}[Cross-Cut Theorem] Let $P$ be a poset and let $C\subset P$ be a cross-cut.  
Then $\ord{P}\heq \Gamma(P,C)$.
\end{lemma}
{\bf Proof.}  The proof is an application of the Nerve Theorem.
Consider the subposets $D_x = P_{\geqs x}\cup P_{\leqs x}$ for $x\in P$.  We claim that
the subcomplexes $\ord{D_c}$, $c\in C$, give a contractible cover of $\ord{P}$.

Condition 2 gives $\ord{P} = \bigcup_{c\in C} \ord{D_c}$.  
Next, say $C'\subset C$.  Then $C'$ is
an antichain, so if $x\in D_{c'}$ for each $c'\in C'$, we either have $x>c'$ for all $c'\in C'$
or $x<c'$ for all $c'\in C'$.  Thus if this intersection is non-empty, $C'$ is bounded
in $P$ and hence has either a meet or a join (by condition 3).  
In either case, this meet or join is 
comparable to each element in the intersection.  So $\bigcap_{c'\in C'} D_c$ is either empty
or a cone, hence contractible.

So we have a contractible cover of $P$, and in fact we see that $\Gamma(P,C)$ is exactly
the nerve of this cover.  This completes the proof.
$\hfill \Box$

\begin{corollary}$\label{prime-complex}$
Let $G$ be a finite group and denote by $Pr(G)$ the simplicial complex
with vertex set $V = \{xH\in \C{G}: o(H)$ is prime$\}$ and with a simplex for each subset of
$V$ bounded in $\C{G}$.  Then $\C{G}\heq Pr(G)$.
\end{corollary}
{\bf Proof.}
It suffices to show that $V$ is a cross-cut in $\C{G}$ (since $Pr(G)$ is then the
cross-cut complex of $V$).  
Certainly, $V$ is an antichain since no two subgroups of prime order are comparable.
Next, any chain $xH_1<xH_2<\cdots <xH_n$ with $H_1\neq \{1\}$ lies above $xP$, where
$P$ is a subgroup of prime order in $H_1$.  If $H_1 = \{1\}$,
then letting $P$ be a prime-order subgroup of $H_2$, we see that $xH_i$ is comparable
to $xP$ for each $i$.
Finally, the third condition of a cross-cut is satisfied trivially, 
since any subset in $\C{G}$ with a lower (upper) bound has a meet (join).
$\hfill \Box$

\begin{claim} $\C{Q_8}$ has the homotopy type of a wedge of three circles, and hence
$\pi_1(\C{Q_8})\homeo \Z\ast\Z\ast\Z$.
\end{claim}
{\bf Proof.}  We will apply Corollary~\ref{prime-complex} to $Q_8$.  Note that the only
subgroup of prime order in $Q_8$ is $\gen{-1}$.  The cosets of $\gen{-1}$ are
$\{1,-1\}$, $\{i,-i\}$, $\{j,-j\}$ and $\{k,-k\}$, and it is easy to check that 
any two are contained
in a proper coset.  So $Pr(Q_8)$ is the complete graph on 4 vertices, and has
six edges.  A maximal tree in $Pr(Q_8)$ has three edges, and after modding out such a tree
we are left with a wedge of three circles.  The result now follows from Van Kampen's Theorem.
$\hfill \Box$

\section{The Coset Poset of $A_5$}

In this section we will consider the smallest example of a group whose coset poset is not
determined (up to homotopy) by Brown's theorem (Theorem~\ref{solvable-coset-poset}).
This group is $A_5$, the first non-solvable group (as well as the 
first non-abelian simple group).  Our goal will be to establish the following theorem.

\begin{theorem}$\label{A_5}$
The coset poset of $A_5$ has the homotopy type of a bouquet of 1560 two-dimensional spheres.
\end{theorem}

It can 
be checked that there are 1018 proper cosets in $A_5$, and hence $\C{A_5}$ has
1018 vertices.  Evidently, $\C{A_5}$ is far too large to admit direct analysis.  

The proof of Theorem~\ref{A_5} will proceed in stages.  First we will show that 
$\mC = \C{A_5}$ has the homotopy-type of a two-dimensional complex.  For
this portion of the proof we will work directly with $\mC$.  Then we will
show that $\mC$ is simply connected by examining $\M{A_5}$.  We will give
two proofs of simple connectivity: the first will be an application of 
Theorem~\ref{2-gen-presentation} and the second will use 2-transitive actions of 
$A_5$. To show that $\C{A_5}$ has the homotopy type of a bouquet of 2-spheres,
we appeal to the general result that a $k$-dimensional complex which is $(k-1)$-connected
is homotopy equivalent to a bouquet of $k$-spheres.

The number of spheres in the bouquet can be calculated from the Euler characteristic 
$\tilde{\chi} (\C{A_5})$.  Brown~\cite{Brown-coset-poset} has shown that 
$\tilde{\chi} (\C{A_5}) = 1560$, and viewing the 
Euler characteristic as an alternating sum of homology ranks, we see that the number
of 2-spheres is exactly 1560.

We will now examine the subgroup lattice of $A_5$, in order to describe the chains in $\mC$ 
of length greater than three (i.e. the simplices in $\ord{\mC}$ of dimension greater than two).
We will consider $A_5$ together with its standard action on the set $\{1,2,3,4,5\}$, and 
we say that a subgroup $H\leqs A_5$ has a fixed point if some element of $\{1,2,3,4,5\}$
is fixed by all elements of $H$.  

We use the notation $D_{2n}$ for the dihedral group of order $2n$ 
(so $D_{2n} = \Z/n\rtimes \Z/2$).

\begin{claim} If $M$ is a maximal subgroup of $A_5$, then $M$ is isomorphic to one of 
$D_{10}$, $A_4$, $S_3$.  There are exactly five copies of $A_4$ corresponding to 
the five subsets of $\{1,2,3,4,5\}$ of cardinality four; there are six copies of $D_{10}$, 
each normalizing a different copy of $\Z/5$; and there are ten copies of $S_3$, all conjugate
to one another in $S_5$.
\end{claim}
{\bf Proof.}  Let $M\leqs A_5$ be maximal.  Note that a subgroup
$H$ of index $n$ in $A_5$ induces a non-trivial homomorphism $A_5\maps S_n$ 
(via the action of $A_5$ on the cosets of $H$), and by simplicity of $A_5$ this map must be 
an injection.  Hence we must have $n\geqs 5$, i.e. $o(H)\leqs 12$.  

Now, any subgroup with a fixed point is contained in a copy of $A_4$, and each copy of 
$A_4$ is maximal in $A_5$ (since $o(A_4) = 12$).
So if $M$ has a fixed point, $M\homeo A_4$.

Next, say $M$ acts transitively.  Then we must have $5|o(M)$, so $o(M)$ is either
five or ten.  Now, $(15)(24)$ normalizes $\gen{(12345)}$, 
so $N(\gen{12345})\neq A_5$ has order at least ten, and in fact its order can be no larger.
Thus $N(\gen{(12345)}) \homeo D_{10}$, and hence every subgroup of order 
five is contained in a copy of $D_{10}$ (subgroups of order five are Sylow subgroups, 
hence conjugate).  So if $M$ acts transitively, $M\homeo D_{10}$.

Finally, if $M$ neither acts transitively nor has a fixed point, then the orbits of 
$\{1,2,3,4,5\}$ under the action of $M$ must have sizes two and three.  In this case, one may
check that the elements of $A_5$ permuting the sets $\{1,2\}$ and $\{3,4,5\}$ amongst themselves
form a subgroup isomorphic to $S_3$ (in fact, this is the standard embedding of $S_{n-2}$ in 
$A_n$).  This subgroup is in fact maximal in $A_5$, since any larger subgroup would be 
transitive and therefore have order divisible by 5.  So in this last case, $M\homeo S_3$.

Now we turn to the number of subgroups of each type.  Any copy of $A_4$ is maximal since
it has order twelve, and thus must be one of the obvious five copies.  Any copy of $D_{10}$
certainly normalizes a $\Z/5$, and so the only such subgroups are the normalizers of the
Sylow 5-subgroups of $A_5$.  There are six such subgroups, as is easily checked.  Finally,
we claim that any copy of $S_3$ is maximal (and hence arrises in the manner described above).
If $H\leqs A_5$ were isomorphic to $S_3$ but not maximal, then $H$ would be contained in a
copy of $A_4$.  Thus it sufficed to check that $A_4$ has no subgroup of order six.  
Any such subgroup
would be normal, and hence would have to contain all the Sylow 3-subgroups of $A_4$.  But
there are four of these, totalling nine elements.  
So all the $S_3$ arrise in the above 
manner, and are clearly conjugate in $S_5$. 
$\hfill \Box$

\begin{corollary}$\label{long-chains}$ All chains in $\mC = \C{A_5}$ of length greater than
three are of the form $\{x\}\leqs xC_2 \leqs xD_4 \leqs xA_4$ (where $C_2$ denotes the
cyclic group of order two).
\end{corollary}
{\bf Proof.}  Underlying every chain in $\mC$ is a chain of the same length in 
$S(A_5)\cup \{1\}$.  These latter chains each end in either a copy of $S_3$, a copy
of $D_{10}$, or a copy of $A_4$.  Since the orders of $S_3$ and $D_{10}$ are products of 
two distinct primes, the longest chains ending at these subgroups are of length three.
So we only need to consider chains ending at $A_4$.  The subgroups of $A_4$ are
of orders two, three and four (there are no subgroups of order six, as discussed earlier).
So any chain of length four must be of the form described in the corollary (noting
that the unique Sylow 2-subgroup of $A_4$ is isomorphic to $D_4$).
$\hfill \Box$

\begin{claim}$\label{A_5-dimension}$
Let $\bar{\mC}$ denote the poset $\mC$ with all cosets of all copies of $D_4$
removed.  Then $\ord{\bar{\mC}}$ is two-dimensional and $\bar{\mC}\heq \mC$.
\end{claim}
{\bf Proof.}  The fact that $\ord{\bar{\mC}}$ is two-dimensional follows immediately from 
Corollary~\ref{long-chains}.  To show that $\bar{\mC}\heq \mC$, we apply Quillen's Theorem
to the inclusion $i:\bar{\mC}\to \mC$.  If $H< A_5$ is not isomorphic to $D_4$, then
$i^{-1}(\mC_{\geqs xH})$ is a cone on $xH$ for any $x\in A_5$.  When $H\homeo D_4$,
we have $i^{-1}(\mC_{\geqs xH}) = \{xK: H< K < G\}$.  Any such $K$ is isomorphic to 
$A_4$, and each $D_4$ lies in a unique $A_4$ (its normalizer in $A_5$).  
So $i^{-1}(\mC_{\geqs xH})$ is a single point, and thus Quillen's Theorem
shows that $i$ is a homotopy equivalence.
$\hfill \Box$

\vspace{.15in} We now turn to simple connectivity.  First we will
show that $A_5$ satisfies the conditions of Theorem~\ref{2-gen-presentation}, thereby
showing $\C{A_5}$ is simply connected.  Then we will give a second, more conceptual
proof using 2-transitive actions of $A_5$.

\begin{claim}$\label{pi-A_5}$
For each pair $x, y\in A_5$ with $\gen{x,y} = A_5$, 
there is an element $z\in A_5$ such that $\gen{x,z}$, $\gen{y,z}$, and $\gen{zx,zy}$ are
all proper subgroups.  Theorem~\ref{2-gen-presentation} now shows that $\C{A_5}$ is simply
connected.
\end{claim}
{\bf Proof.}  We will prove this result by considering representatives of all 
possible automorphism classes of generating pairs (i.e. all sets of the form 
$\{(\phi (x), \phi (y)): \phi \in Aut(A_5), \,\, \gen{x,y} = A_5\}$) and finding, 
for each such representative pair, an appropriate element
$z$ (in some cases it is not actually necessary to choose a specific representative).  
Luckily the total number of automorphism classes is relatively small.  We will
break the problem down by considering the orders of $x$ and $y$.
\begin{itemize}
	\item {\bf o(x) = o(y) = 2:} In any group, the subgroup generated by two elements of order
		two is dihedral, and in particular two such elements cannot generate a simple 
		group.  (The latter statement is easily checked by verifying that 
		$x, y\in N(\gen{xy})$, which is a proper subgroup by simplicity.)

	\item {\bf o(x) = 2, o(y) = 3:}  Since $x$ and $y$ generate $A_5$ they cannot
		have a common fixed point.  Up to automorphism, we may assume that 
		$x = (1 2)(4 5)$.  The three-cycle $y$ must move an element from
 		each two-cycle, since otherwise $x$ and $y$ would sit in a copy of $S_3$.
		So up to automorphism, we may assume that $y = (2 3 4)$.
		Now, setting $z = (1 2 3)$ we see that $\gen{x,z}$ lies a copy of $S_3$ (its
		action has orbits $\{1,2,3\}$, $\{4,5\}$) and $\gen{y,z}$ fixes 5 and is 
		contained in an $A_4$.  Furthermore, $zx$ and $zy$ have order $2$, so 
		$\gen{zx,zy}\neq A_5$ by the above discussion.

	\item {\bf o(x) = 2, o(y) = 5:}  First, $x$ lies in a copy of $A_4$ and $y$ lies in 
		a copy of $D_{10}$.  These subgroups must intersect nontrivially in $A_5$,
		since otherwise we would have $o(A_5)\geqs o(A_4)o(D_{10}) = 120$, 
		a contradiction.  Since the only common divisor of $o(D_{10})$ and $o(A_4)$
		is two, there must be an element $z$ of order two in this intersection.
		(Note that $x$ does not lie in the intersection, since $\gen{x,y} = A_5$.)
		By choice of $z$, we see that $\gen{z,x}$ and $\gen{z,y}$ are proper subgroups.
		In addition, since $D_{10}$ is dihedral the product of an element of order
		two and an element of order five has order two, i.e. $o(zy) = 2$.  Similarly
		the product of two elements of order two in $A_4$ has order two, so $o(zx) = 2$.
		So once again we see that $\gen{zx,zy}\neq A_5$.

	\item {\bf o(x) = o(y) = 3} Let $x = (1 2 3)$.  In order for $x$ and $y$ to generate,
		$y$ must move both 4 and 5, and up to automorphism we may assume that 
		$y = (3 4 5)$.  Let $z = (2 3 4)$.  Then $x$ and $z$ fix 5,
		and $y$ and $z$ fix 1, so neither pair generates $A_5$.  Once again
		we find that $o(zx) = o(zy) = 2$.

	\item {\bf o(x) = 3, o(y) = 5:}  Let $x = (1 2 3)$.  It is easily checked that 
		there is an element $z = (a b)(4 5)\in N(\gen{y})$, with $a,b\in \{1, 2, 3\}$.
		Now, $x$ and $z$ each lie in the copy of $S_3$ with orbits
		$\{1, 2, 3\}$ and $\{4, 5\}$, and hence do not generate $A_5$.  
		As discussed above, $o(zy) = 2$ and since $zx$ contains the transposition
		$(45)$ it has order two as well.

	\item {\bf o(x) = o(y) = 5:}  The normalizers of $\gen{x}$ and $\gen{y}$ must intersect
		non-trivially, since otherwise we would have 
		$o(A_5)\geqs o(D_{10})^2$, a contradiction.  Since $x$ and $y$ generate,
		$\gen{x}\neq \gen{y}$ and hence $\gen{x}\cap\gen{y} = \{1\}$. 
		This implies that the normalizers of $\gen{x}$ and $\gen{y}$ meet in a subgroup
		of order two.  If $z$ is the non-trivial element in this intersection, then
		we see that $\gen{x,z}$ and $\gen{y,z}$ are proper and $o(zx) = o(zy) = 2$.
\end{itemize}
This covers all possible generating pairs, and thus completes the proof.  
$\hfill \Box$

\vspace{.15in}
\noindent {\bf Alternate Proof of Simple Connectivity.}
In the above proof, it was shown that for each edge 
$\{x,y\}\in \mathcal{M}^2(A_5)$ with $o(x) = 2$, there is an element $z\in G$ 
satisfying the conditions of Theorem~\ref{2-gen-presentation}.  Consider
the presentation for $\pi_1 (\M{A_5})$ given by Theorem~\ref{presentation}, using
as a maximal tree all edges $\{1,g\}$, $g\in G$ ($g\neq 1$).  Note that this
is the same presentation as was used in the proof of Theorem~\ref{2-gen-presentation},
and recalling that proof, the existence of the element $z$ implies that the generators
$(x,y) = (y,x)^{-1}$ are trivial in $\pi_1 (\M{A_5})$.  

Starting from this fact 
(that generators $(x,y)$ with $o(x) = 2$ are trivial) we will
complete the proof in a less ad hoc manner.  The same method will be used
in the next section to show that $\C{PSL_2 (\F_7)}$ is simply connected.

Let $\{g,h\}$ be any edge of $\mathcal{M}^2 (A_5)$.  
We must show that the corresponding generators
of $\pi_1 (\M{A_5})$ are trivial.
It will suffice to show that there exists an element
$z$ of order two and a subgroup $K\in \C{A_5}$ s.t. $h\equiv g \equiv z$ (mod $K$).  This
is true because then the set $\{g,h,z\}$ forms a 2-simplex in $\M{A_5}$ and then the relations
in Theorem~\ref{presentation} show that $(g,h)(h,z)(z,g) = 1$.  Since $o(z) = 2$,
we have $(h,z) = (z,g) = 1$, so $(g,h) = 1$ as desired.

Now, say $G$ acts 2-transitively on the set $X$, and assume that some element $z$ with
$o(z) = 2$ acts non-trivially.  Let $H = Stab (x)$ for some $x\in X$, and note that
for any $y\in X$, the set $S_{x,y} = \{g\in G: g\cdot x = y\}$ is just the coset $k H$, 
where $k\cdot x = y$.  Now, since $z$ acts non-trivially and the action is 2-transitive,
some conjugate of $z$ sends $x$ to $y$, i.e. there is an element of order two in 
every non-trivial coset of $H$.  Thus we can prove simple connectivity by showing that
every two-element set $\{g,h\}\subset G$ lies in a coset $kH$ where $G$ acts 2-transitively
on some set $X$, $H$ acts as the stabilizer of a point, and some element of order
two acts non-trivially.

First, consider the standard action of $A_5$ on the set $\{1,2,3,4,5\}$.  This action is
clearly 2-transitive, and the stabilizer of a point is isomorphic to $A_4$.  Every non-trivial
element acts non-trivially, so certainly elements of order two act non-trivially.

Next, consider the action of $A_5$ on its Sylow 5-subgroups (by conjugation).  Call
this set $S$.  There are six subgroups in $S$, and each of them acts transitively on 
the other five (certainly each acts non-trivially, and since the orbits must have order 
dividing 5, the action is transitive).  Now 
given any $H_1 \neq K_1$, $H_2 \neq K_2 \in S$, there is an element $y\in H_1$ with
$K_1^y = K_2$ and there is an element $z\in K_2$ with $H_1^z = H_2$.  We have
$H_1^{yz} = H_1^z = H_2$ and $K_1^{yz} = K_2^z = K_2$, showing that the action is indeed
2-transitive.  Note that the stabilizer of a subgroup in $S$ is just its normalizer,
a copy of $D_{10}$.

We claim that every two-element set $\{g,h\}$ lies in either a coset of a copy of $A_4$
or in a coset of a copy of $D_{10}$.  In other words, we claim that $g^{-1}h$ always
lies in a subgroup isomorphic to either $A_4$ or $D_{10}$.  But every cyclic subgroup of $A_5$
has order two, three, or five, so this follows immediately.
$\hfill \Box$

\section{The Coset Poset of $PSL_2 (\F_7)$}

In this section we will consider the finite simple group $G = PSL_2 (\F_7)$.  This group
has order 168, and is the next larger simple group after $A_5$.  Our goal is
to show that $\C{G}$ is simply connected.  Unfortunately, we do not know whether $\C{G}$ is 
(homotopy equivalent to) a bouquet of spheres or something more exotic.

The proof of simply connectivity will require a careful understanding of the
subgroup poset of $G$, 
and in particular we will need to compute the M\"{o}bius function of
$G$ (defined in Section~\ref{mu(psls)-section}).  
Rather than restrict ourselves to the case $p = 7$, we will give a detailed
description of the subgroups of $PSL_2 (\F_p)$ for all primes $p\geq 5$.  This computation
will be based on Burnside's exposition in~\cite{Burnside-theory-of-groups}.  
The computation of the M\"{o}bius function of 
$G$ will be based on the (brief) calculation given by Hall in his original 
paper~\cite{Hall-eulerian} on M\"{o}bius inversion.

The proof of simple connectivity will be analogous to the second proof for $A_5$ 
(given in the previous section), and will utilize 2-transitive actions of $G$.
The M\"{o}bius function of $G$ will be used to count automorphism classes of 
generating pairs.  In the case of $A_5$ this was unnecessary simply because these classes were
easy to enumerate by hand.

Once we have proven simple connectivity, we will briefly discuss other facts about the
homotopy-type of $\C{G}$.

\vspace{.15in}
We will begin by defining the groups we propose to study, and then describing a useful 
permutation representation for them.

Fix a prime $p$, and let 
$\F_p \isom \Z/p$ denote the finite field with $p$ elements.  We define
$GL_n (\F_p)$, the $n$-dimensional general linear group over $\F_p$,
to be the multiplicative group of invertible $n\cross n$ matrices with entries
in $\F_p$, and we define the $n$-dimensional special linear group, 
$SL_n (\F_p) \leqs GL_n (\F_p)$, to be the subgroup of matrices with
determinant one.  

Let $\Gamma_{n,p}$ denote the center of $SL_n (\F_p)$.
The $n$-dimensional projective special linear group, 
$PSL_n (\F_p)$, is defined to be the quotient $SL_n (\F_p) / \Gamma_{n,p}$.
The following theorem is well-known 
(see, for example,~\cite[Theorem 3.2.9]{Robinson-theory-of-groups}).

\begin{theorem}$\label{simplicity}$
If $n>2$ or $p>3$, then $PSL_n (\F_p)$ is a simple group.
\end{theorem}

\vspace{.15in}
We will mainly be interested in the cases where $n=2$ and $p>3$.
In these cases, direct computation shows that $\Gamma_{2,p} = \{I, -I\}$, 
so $o(PSL_2 (\F_p)) = \frac{1}{2} o(SL_2 (\F_p))$.
Our first goal is to determine the order of $PSL_2 (\F_p)$.

\begin{claim}$\label{order}$
If $p>3$, then $o(PSL_2 (\F_p)) = \frac{1}{2} p (p^2 -1)$
\end{claim}
{\bf Proof.} An element of $\gl$ may be thought of as a pair $(v, w)$ where $v$ and $w$
are linearly independent vectors in $\F_p^2$.  Thus in forming an element of $\gl$, we
have $p^2-1$ choices for $v$ ($v$ must simply be non-zero), and then 
$p^2 - p$ choices for $w$ ($w$ must not lie in the span of $v$).  Thus we have
$o(\gl) = (p^2 - 1)(p^2 - p)$.  

The exact sequence
$$1 \maps \spl 
    \hookrightarrow \gl
    \stackrel{det}{\maps} \F^*_p
    \maps 1$$
(where $det$ denotes the determinant homomorphism) shows that 
$$o(\spl) = \frac{(p^2 - 1)(p^2 - p)}{p-1} = (p+1)(p^2 - p) = p (p^2 - 1).$$  
As mentioned above, $\psl$ has index two in $\spl$, completing the proof.
$\hfill \Box$

\vspace{.15in}

From now on we will write $o(\psl) = 2pqr$, where $q = \frac{p-1}{2}$ and
$r = \frac{p+1}{2}$.  The main advantage of this notation is the fact that $p$, $q$ and $r$
are pair-wise relatively prime.

\vspace{.15in}
We will now describe a useful action of $\psl$.
Let $V = \F_p \oplus \F_p$ be a two-dimensional vector space over $\F_p$ and define 
an equivalence relation $\sim$ on $V-\{0\}$ by setting 
$v\sim \lambda v$ for all $v\in V-\{0\}$, $\lambda \in \F^*_p$.
We identify the set of equivalence classes with the projective line $\prj$ 
(where $\infty$ is just a formal symbol) via the map $(v_1, v_2)\mapsto \frac{v_1}{v_2}$, 
with the convention that $\frac{a}{0} = \infty$ for any $a\in \F_p$
(this map is clearly well-defined on equivalence classes, and it is not hard to check
that it is a bijection).

Then $\gl$ acts on $V - \{0\}/ \sim$ (and hence on $\prj$) by multiplication.  Restricting
to $SL_2 (\F_p)$, it is clear that this action factors through the quotient map 
$\pi:SL_2 (\F_p)\to PSL_2 (\F_p)$, and we obtain an action of $\psl$ on $\prj$.

\vspace{.15in}
It will be convenient to describe another way of writing elements in $\psl$.
Consider the group $M_p$ of (invertible) M\"{o}bius transformations $\prj \to \prj$.  These 
are maps of the form $x\mapsto \frac{ax+b}{cx+d}$ (where $x\in \prj$ and $a, b, c, d\in \F_p$)
with $ad - bc = 1$, and multiplication is defined by composition.  
We use the conventions $\frac{a\infty + b}{c\infty + d} = \frac{a}{c}$
and $\frac{z}{0} = \infty$ (in all other situations, we simply add and multiply in $\F_p$).
We write elements
of $M_p$ in the form $\frac{ax+b}{cx+d}$, where 
$x$ is an indeterminate, and it is easy to check
that the inverse of $\frac{ax + b}{cx + d}\in M_p$ is 
$\frac{dx-b}{-cx + a}$.  Since these maps
are invertible, they are bijections and hence $M_p$ acts on $\prj$.
Note that in this notation the identity element of $M_p$ is simply $x$, but to avoid
confusion we will sometimes denote the identity of $M_p$ by 1 or by $id_M$.

It is easy to check that the obvious map $\spl\maps M_p$ 
is a homomorphism with kernel $\{I, -I\}$,
and thus induces a (canonical) isomorphism $\psl \maps M_p$.  In addition, this isomorphism 
commutes with the actions on $\prj$.  It is easy to see that the action of $M_p$ on $\prj$
is transitive, and we will later see that it is in fact 2-transitive.

\subsection{The Subgroups of $PSL_2 (\F_p)$}

We will now describe the subgroups of $\psl \isom M_p$, $p>5$.  
We will refer mainly to $M = M_p$
during this analysis.  Our first goal will be to prove the following theorem describing the 
cyclic subgroups of $M$ and their normalizers.  The proof will
be broken down into a number of claims.

\begin{theorem}$\label{cyclic-subgroups}$
The set of maximal cyclic subgroups of $M$ consists of:
\begin{enumerate}
	\item $p+1$ cyclic subgroups of order $p$, each fixing a single point of $\prj$;
	\item $pr$ cyclic subgroups of order $q$, each fixing two points of $\prj$;
	\item $pq$ cyclic subgroups of order $r$, each acting freely on $\prj$.
\end{enumerate}

Each of the three collections above forms a conjugacy class in $M$.
In addition, the intersection of any two of the subgroups listed is trivial, so
each element of $M$ lies in a unique maximal cyclic subgroup.

Any cyclic subgroup of order $p$ has normalizer isomorphic to $\Z/p \rtimes \Z/q$
(for some action of $\Z/q$ on $\Z/p$, which we will not need to specify).

The normalizer of any other (non-trivial) cyclic subgroup $C\leqs M$ 
is isomorphic to either $D_{2q}$ or $D_{2r}$, depending on whether $o(C)$ divides $q$ or $r$
(we use the convention $D_{2n} = \Z/n\rtimes \Z/2$, where the generator of $\Z/2$ acts by 
inversion).
\end{theorem}
\vspace{.15in}

Note that by ``maximal cyclic subgroup" we mean a cyclic subgroup not contained in any
larger cyclic subgroup.  Also, recall that $p$, $q$ and $r$ are pair-wise relatively
prime, so in showing that none of the cyclic subgroups listed above intersect non-trivially,
we need only consider intersections between groups of the same order.

Before beginning the proof, we record two useful lemmas.

\begin{lemma}$\label{norm-conj}$
Let $G$ be a finite group and $H\leqs G$ a subgroup.  Then the number of 
conjugates of $H$ in $G$ is equal to $(G:N_G (H))$, the index of the normalizer of $H$.
\end{lemma}

This result follows easily by considering the actions of $G$ on the cosets of 
$N(H)$ (by left multiplication) and on the conjugates of $H$ (by conjugation).

\begin{lemma}$\label{normalizers}$
Let $G$ be a finite group and let $H\leqs G$ be a cyclic subgroup.
Assume that any two conjugates of $H$ intersect trivially.  Then for any
$K\leqs H$, we have $N_G (K) = N_G (H)$.
\end{lemma}
{\bf Proof.}  Since $K$ is the unique subgroup of its order in $H$, we have
$K\normal N_G (H)$.  In addition, we see that $K$ has the same number of conjugates
as $H$, since any two conjugates of $H$ intersect trivially.  
Letting $c_H$ and $c_K$ denote the number of conjugates of $H$ and $K$ in $G$,
Lemma~\ref{norm-conj} gives
$$o(N_G (K)) = \frac{o(G)}{c_K} = \frac{o(G)}{c_H} = o(N_G (H)),$$
completing the proof.
$\hfill \Box$
\vspace{.15in}

We will now examine the subgroups of order $p$ in $M$.  
Note that these are in fact Sylow 
$p$-subgroups, and hence form a single conjugacy class in $M$.
Any two of intersect trivially, since they are of prime order.

\begin{claim}$\label{p-groups}$
There are $p+1$ cyclic subgroups of order $p$ in $M$, each with normalizer
isomorphic to $\Z/p \rtimes \Z/q$.
\end{claim}
{\bf Proof.}  Consider the element $x+1\in M$.  Note that $(x+1)^n = x+n$, and hence
$o(x+1) = p$.  Thus $\gen{x+1}$ is cyclic subgroup of order $p$, and note that this
subgroup fixes $\infty\in \prj$ but does not fix any other point.

If $m\in M$, then $m\gen{x+1} m^{-1}$ will fix the point $m(\infty)$.  
Thus $N\left(\gen{x+1}\right) \leqs{\rm Stab}(\infty)$, since
$m\notin{\rm Stab}(\infty)$ implies $m\gen{x+1} m^{-1}$ does not fix $\infty$.  
There are $p+1$ points in $\prj$, giving $p+1$ distinct conjugates of $\gen{x+1}$
(recall that the action of $M$ on $\prj$ is transitive).  

To complete the proof it will suffice to show that Stab$(\infty) \iso \Z/p\rtimes \Z/q$.
Since $\frac{a\infty + b}{c\infty + d} = \frac{a}{c}$, we see that
$${\rm Stab}(\infty) = \left\{ \frac{ax + b}{cx + d}\in M: c = 0 \right\}
               = \left\{ \frac{ax + b}{a^{-1}} \right\}.$$
Since we must have $a\neq 0$, there are $\frac{1}{2} (p-1)p = pq$ elements in this
subgroup.  Note that $\gen{x+1}$ is a Sylow p-subgroup
of Stab$(\infty)$, and thus is either normal or has at least $p+1$ conjugates within
Stab$(\infty)$.  These $p+1$ conjugates would contain a total of 
$(p+1)(p-1) + 1 > pq = o\left({\rm Stab}(\infty)\right)$ elements, and thus we have 
$\gen{x+1}\normal$ Stab$(\infty)$.

Finally, note that there is a copy of $\Z/q$ in Stab$(\infty)$, generated by the element
$\frac{zx}{z^{-1}}$ with $z$ a generator of $\F_p^*$.  (Since
$\left(\frac{zx}{z^{-1}}\right)^n = \frac{z^n x}{z^{-n}}$, we have 
$\frac{z^n x}{z^{-n}} = x = id_M$ iff $z^n = z^{-n}$, i.e. iff $z^{2n} = 1$.   
Now, $o(z) = o(\F_p^*) = p-1$, so $o\left(\frac{zx}{z^{-1}}\right) = \frac{p-1}{2} = q$.)
Thus we have found a complement (in Stab$(\infty)$) to the normal subgroup $\gen{x+1}$,
and Stab$(\infty)\isom \Z/p\rtimes \Z/q$.
$\hfill \Box$
\vspace{.15in}

We now turn to the cyclic subgroups of order $q$.  We have already exhibited one such subgroup
($\gen{\frac{zx}{z^{1-}}}$ for $z$ a generator of $\F_p^*$) 
in the proof of Claim~\ref{p-groups}, and the other such subgroups are conjugates of
this one.  First we will need two simple lemmas about the action of $M$ on $\prj$.  

\begin{lemma}$\label{2-transitive}$
 The group $M$ acts 2-transitively on $\prj$.
\end{lemma} 
{\bf Proof.}
Note that, as shown in the proof of Claim~\ref{p-groups}, for each $x\in \prj$,
there is a subgroup $P_x$ of order $p$ in Stab$(x)$.  Since no element of $M$ other than
the identity acts trivially on $\prj$, we see that $P_x$ acts transitively on $\prj - \{x\}$
(because this is a set of cardinality $p$).
Now, given any $a\neq b$, $c\neq d \in \prj$, there is an element $y\in P_a$ with
$y(b) = d$ and there is an element $z\in P_d$ with $z(a) = c$.  We have
$zy(a) = z(a) = c$ and $zy(b) = z(d) = d$, showing that $M$ is indeed 2-transitive.
$\hfill \Box$

\begin{lemma}$\label{fixed-points}$
 No non-trivial element of $M$ can fix more than two points in $\prj$.
\end{lemma}  
{\bf Proof.}  Notice that fixed points (other than $\infty$) correspond to solutions 
of a quadratic polynomial: if $m = \frac{ax+b}{cx+d} \in M$, then 
$$m(x) = x \iff ax+b = cx^2 + dx \iff cx^2 + (d-a)x - b = 0.$$
If $m$ has more than two fixed points in $\F_p$, this polynomial has more than two roots and must
be zero.  Thus $m$ fixes every point of $\prj$, i.e. $m$ is the identity.

Finally, if $m (\infty) = \infty$, then $m = \frac{ax+b}{d}$ and any other fixed point $x$
satisfies $\frac{ax+b}{d} = x$, i.e. $(a-d)x - b = 0$, and if $m\neq 1$ this linear
equation can have at most one root.  
$\hfill \Box$

\vspace{.15in}
For $x,y\in \prj$, let Stab$(x,y) = $Stab$(x)\cap$Stab$(y)$.

\begin{claim}$\label{q-groups}$
For any two distinct points $x,y\in \prj$, Stab$(x,y) \isom \Z/q$.
This yields $pr$ conjugate copies of $\Z/q$, and any two intersect trivially.
Furthermore, if $t\in$ Stab$(x,y)$ ($x\neq y$), then 
$N(\gen{t}) = N($Stab$(x,y))\isom D_{2q}$.
\end{claim}
{\bf Proof.}  It is easily checked that 
Stab$(0,\infty) = \left\{ \frac{ax}{a^{-1}}: a\in \F_p \right\}$, and this
subgroup is generated by any element $\frac{zx}{z^{-1}}$ with $z$ a generator of $\F_p$.
As was shown in the proof of Claim~\ref{p-groups}, such an element has order $q$.
Since $M$ is 2-transitive, Stab$(x,y) \isom$ Stab$(0,\infty)$
for any $x\neq y \in \prj$.  There are ${p+1 \choose 2} = pr$ such stabilizers,
and two of them intersect trivially by Lemma~\ref{fixed-points}.

Finally, we must consider the normalizers of these subgroups.  
Clearly, Stab$(0,\infty)\leqs$Stab$(\{0, \infty\})$, where the latter
denotes the elements of $M$ which fix $\{0,\infty\}$ as a set.  We claim
that Stab$(\{0, \infty\}) = N($Stab$(0,\infty)$ and that 
Stab$(\{0, \infty\})\isom D_{2q}$.

Any element 
$m\in$ Stab$(\{0, \infty\})$ either fixes both points or swaps them.  In the former
case $m\in$ Stab$(0,\infty)$ and in the latter we have 
$m = \frac{b}{-b^{-1} x}$ for some $b\in \F_p^*$.  In either case we obtain $q$ elements,
so $o \left({\rm Stab}(\{0, \infty\})\right) = 2q$ and hence 
$${\rm Stab}(0,\infty)\normal{\rm Stab}(\{0, \infty\}).$$
In addition, for any
$b\in \F_p^*$ we have 
$$\left( \frac{b}{-b^{-1} x}\right)^2 
= \frac{b}{-b^{-1} \frac{b}
                        {-b^{-1} x}
          }
= x = id_M,$$
so all elements in Stab$(\{0,\infty\}) - $Stab$(0,\infty)$ have order two.  To 
show that Stab$(\{0,\infty\})$ is dihedral, all that remains to be checked is that
an element in Stab$(\{0,\infty\}) - $Stab$(0,\infty)$ conjugates any element
of Stab$(0,\infty)$ to its inverse.  This is easily verified by multiplying 
representative matrices in $\spl$.  Lemmas~\ref{norm-conj} and~\ref{normalizers} 
complete the proof.  
$\hfill \Box$
\vspace{.15in}

We now come to the cyclic subgroups of order $r$.  These arrise in a rather different
manner, since they act freely on $\prj$.

Let $K = \F_{p^2}$, and fix a basis $\beta$ for $K$ over $k = \F_p$.
For $z\in K^*$, let $L_z \in \gl$ denote the matrix (in the basis $\beta$) corresponding 
to the map $y\mapsto zy$ ($y\in K$), and let $H = \{L_z: z\in K^*\} \leqs \gl$.
Let $\tilde{H} = \phi \left( H\cap \spl \right)$, where $\phi:\spl \to M$ denotes the obvious
map.  We will show that $\tilde{H} \isom \Z/r$.

We now record two useful lemmas about the Galois extension $K/k$.

\begin{lemma}$\label{sqrt}$  Every element of $k$ has a square root in $K$.
\end{lemma}
{\bf Proof.}
This follows from the fact that any degree two extension of $k$ is isomorphic to 
$K$, but can also be seen directly, as follows.  
The (multiplicative) group $K^*$ is isomorphic to the (additive) group
$\Z / (p^2 -1)$, and squares in $K^*$ correspond to even residues in $\Z / (p^2 - 1)$.
The subgroup $k^* \leq K^*$ has order $p-1$, and maps to $\gen{p+1}\leq \Z / (p^2 -1)$.
Since $p+1$ is even, these elements all have ``square roots" in $\Z / (p^2 - 1)$ and
thus all elements of $k^*$ have square roots in $K^*$.
$\hfill \Box$

\begin{lemma}$\label{galois}$
Let $\sigma: K\to K$ denote the $p^{\rm th}$-power map, i.e.
$\sigma (z) = z^p$.  Then Gal$(K/k) = \{I, \sigma\}$, and
$det(\sigma) = -1$.
\end{lemma}
{\bf Proof.}  First, the fact that $\sigma$ is a field automorphism of $K$
follows from the binomial formula
$$(x+y)^p = \sum_{i=0}^{p} {p \choose i} x^i y^{p-i},$$
since for $i\neq 0, p$ we have $p| {p\choose i} = \frac{p!}{i! (p-i)!}$.
Next, Fermat's Little Theorem shows that for any $z\in k$, $z^p = z$ and hence
$\sigma$ is an automorphism of $K/k$.  Finally, since $deg(K/k) = 2$, $Gal(K/k)\isom \Z/2$
and hence $Gal(K/k) = \{I, \sigma\}$.

Finally we must compute $det(\sigma)$. (Note that an element of the Galois group fixes $k$
and thus is also in $\gl$, so $det(\sigma)$ makes sense).  By the Normal Basis Theorem,
there is an element $y\in K$ such that $\{y, \sigma (y)\}$ is a basis for $K$ over $k$.
Since $\sigma$ has order two, it interchanges these basis vectors and thus has determinant
$-1$.
$\hfill \Box$

\begin{claim}$\label{r-groups}$
The subgroup $\tilde{H}$ has $pq$ conjugates in $M$.
Any two of these subgroups intersect trivially, and we have:
\begin{enumerate}
	\item $\tilde{H} \isom \Z /r$; 
	\item $\tilde{H}$ acts freely on $\prj$;
	\item For any non-trivial subgroup $T\leqs \tilde{H}$,
                  $N(T) \isom D_{2r}$.
\end{enumerate}
\end{claim}
{\bf Proof.}  First note that $\tilde{H}$ is cyclic because $K^*$ is cyclic.  To 
determine the order of $\tilde{H}$ we must first determine 
$o(H\cap \spl)$, i.e.
we must determine which elements of $H$ have determinant one.  

It is a result
in Galois theory that for any $z\in K^*$, $det(L_z) = N_{K/k} (z)$ where 
$N_{K/k} (z)$ denotes
the norm of $z$, i.e. the product of its Galois conjugates in $K$.   
Lemma~\ref{galois} now shows that $det(L_z) = z \sigma (z) = z^{p+1}$.  Thus
$H\cap \spl \isom \{z\in K^*: z^{p+1} = 1\} \isom \{a\in \Z /(p^2-1): (p+1)a = 0\}$,
and the last set clearly has cardinality $p+1$.  Finally, since $L_{-1} = -I$, 
$$o(\tilde{H}) = \frac{1}{2} o\left(H\cap \spl \right) = \frac{p+1}{2} = r.$$

Next we consider the action of $\tilde{H}$ on $\prj$.  We want to show that this action
is free, i.e. that no element of $\tilde{H}$ fixes a point on $\prj$.  Consider
the action of $H\cap \spl$ on $\prj$.  If $L_y\in H\cap \spl$ 
fixes a point in $\prj$,
then $L_y (v) = yv = \lambda v$ for some $v\in K -\{0\}$ and some 
$\lambda\in k$.  But then $y = \lambda$, i.e. $y\in k$, and we have  
$$1 = det(L_y) = N(y) = y \sigma (y) = y^2$$
(since $\sigma$ fixes $k$).  Hence $y = \pm 1$ and $\phi (L_y) = \phi (-I) = id_M$.   
So $\tilde{H}$ acts freely on $\prj$.

We now consider $N_M (\tilde{H})$.  Let $y\in K^*$ be an
element with $det(L_y) = -1$.  (Such elements exist, since under the isomorphism
$K^*\stackrel{\isom}{\maps} \Z/(p^2-1)$ they correspond to elements $t\in \Z/(p^2-1)$
with $(p+1)t = \frac{p^2-1}{2}$.)  By Lemma~\ref{galois} we have 
$det(\sigma L_y) = (-1)^2 = 1$, so $\sigma L_y\in \spl$.  Let $\alpha = \phi (\sigma L_y)$.
We claim that $\gen{\tilde{H}, \alpha} \isom D_{2r}$.

To prove this, it suffices to show that $o(\alpha) = 2$ and that for any $h\in \tilde{H}$,
$\alpha h \alpha^{-1} = h^{-1}$.  First, if $v\in K$ then 
$$\left(\sigma L_y \right)^2 v = \left( \sigma L_y \sigma L_y \right) v 
= \sigma \left(y \sigma (yv) \right) = \sigma (y) yv = det(y) v = -v$$
and hence $(\sigma L_y)^2 = -I$ and $\alpha^2 = id_M$.  Thus $o(\alpha)\leqs 2$, and we will 
see that $\alpha$ is non-trivial because its action on $\tilde{H}$ is non-trivial.

Next we wish to prove that $\alpha h \alpha^{-1} = h^{-1}$ for any $h\in \tilde{H}$.
Choose $z\in K^*$ with $\phi (L_z) = h$.  Then for any $v\in K$ we have
$$\left(\sigma L_y \right) L_z \left(\sigma L_y \right)^{-1} v
 = \left(\sigma L_y L_z L_y^{-1} \sigma \right) v
 = \sigma \left( y z y^{-1} \sigma (v) \right)$$
$$ = \sigma \left( z \sigma (v) \right) 
 = \sigma (z) v
 = z^{-1} v
 = L_{z^{-1}} v,$$
where the next to last equality follows because 
$1 = det(L_z) = z \sigma (z)\implies \sigma (z) = z^{-1}$.
So $\left( \sigma L_y \right) L_z \left( \sigma L_y \right)^{-1} = L_{z^{-1}}$ and hence
$\alpha h \alpha^{-1} = h^{-1}$, as desired.

By Lemma~\ref{norm-conj}, 
the proof will be complete once we show that $\tilde{H}$ has $pq$ conjugates, any two of which
intersect trivially.
We will describe $pq$ different bases for $K$ over $k$, and changing basis in $\gl$
will yield $pq$ distinct conjugates of $\tilde{H}$.
Consider all bases of the form $\beta_{y,l} = \{\sqrt{y} + l, 1\}$ where $y,l\in \F_p$ 
and $y$ is not a square in $\F_p$ (there is ambiguity in the symbol $\sqrt{y}$ 
since there are
two roots in $K$ of the equation $x^2 = y$; for each $y$ we arbitrarily choose one of 
these roots and denote it by $\sqrt{y}$).  
Each of these is clearly a basis (since $\sqrt{y}\notin \F_p$), and there are $pq$ of
them since we have $\frac{p-1}{2} = q$ choices for $y$ and $p$ choices for $l$.

Let $H^{\beta_{y,l}}$ denote the conjugate of $H$ obtained by changing basis from $\beta$
to $\beta_{y,l}$.  By writing out the matrices involved, one finds that
$$H^{\beta_{y,l}} \cap H^{\beta_{z,n}} \cap \spl = \{I, -I\}$$
unless $y = z$ and $l = n$.  This implies that the resulting conjugates of $\tilde{H}$ 
have trivial intersection in $M$.  
$\hfill \Box$

\vspace{.15in}
To complete the proof of Theorem~\ref{cyclic-subgroups}, we must simply show that
each element of $M$ lies in one of the cyclic subgroups we have described.
This follows from a simple counting argument: the total number of
elements involved in these subgroups is
$$(p+1)(p-1) + pr(q-1) + pq(r-1) + 1 
= p^2 + 2pqr - pr - pq$$
$$= 2pqr +p^2 - \frac{p(p+1)}{2} - \frac{p(p-1)}{2}
  = 2pqr = o(M).$$

\vspace{.15in}
Our next goal will be to study subgroups $Q\leqs M$ with $Q\isom D_4$.  Following Burnside,
we call such subgroups \emph{quadratic}.  First we record some simple facts about
dihedral groups.

\begin{lemma}$\label{dihedral}$
Let $C\leqs D_{2n}$ be the cyclic subgroup of index two, and say $y\notin C$.  Then
$o(y) = 2$ and for any $c\in C$, $ycy^{-1} = c^{-1}$.  

If $H\leqs D_{2n}$, then $H$ is either cyclic or dihedral.
If $H$ is cyclic then either $H\leqs C$ or $H\isom \Z/2$, and if $H$ is dihedral
then $H = \gen{H\cap C, t}$ for any $t\in H - H\cap C$.
\end{lemma}

The proof is not difficult, and is left to the reader.

\begin{claim}$\label{quadratic-groups}$
The number of distinct quadratic subgoups in $M$ is 
$$\frac{o(M)}{12} = \frac{p (p^2 - 1)}{24}.$$
If $o(M)$ is divisible by eight,
these subgroups fall into two conjugacy classes of equal size and otherwise they 
form a single conjugacy class.
\end{claim}
{\bf Proof.}  Let $H\leqs M$ be a maximal dihedral subgroup, and assume that 
$4|o(H)$ (so $H\isom D_{2r}$ if $r$ is even and $H\isom D_{2q}$ if $q$ is even;
recall that $(r,q) = 1$).
Let $C\leqs H$ be the cyclic subgroup of index two, and let $t\in C$ be the
unique element (in $C$) of order two.  Then if $y\in H - C$, Lemma~\ref{dihedral}
shows that $\gen{y,t}$ is quadratic.  Every quadratic subgroup of $H$ may be
obtained in this manner, since by Lemma~\ref{dihedral},
any quadratic subgroup $Q\leqs H$ contains $t$ and also contains two elements of $H-C$.
So as $y$ ranges over $H-C$ the subgroups $\gen{y, t}$ range over all quadratic subgroups
of $H$, and each quadratic subgroup is counted twice.  Thus there are
$\frac{o(C)}{2}$ quadratic subgroups in $H$.

Next we claim that each quadratic subgroup $Q\leq M$ lies in exactly three conjugates of $H$,
i.e. in exactly three copies of $D_{2r}$ or $D_{2q}$, depending on the parity of $r$ and $q$.
This follows from Theorem~\ref{cyclic-subgroups}: each of the three non-trivial elements
in $Q$ lies in a different conjugate of $C$, and since $Q$ is abelian
$Q$ lies in the corresponding conjugates of $H = N_M (C)$.  Certainly
$Q$ cannot lie in a fourth conjugate of $H$, since then two of the corresponding conjugates
of $C$ would intersect non-trivially.

The total number of quadratic subgroups of $M$ is now seen to be
$$\frac{\frac{o(C)}
             {2}
        c_M (H)
       }
       {3}$$
where $c_M (H)$ denotes the number of conjugates of $H$ in $M$.  
We have either $o(C) = r = \frac{p+1}{2}$ and $c_M (H) = pq = \frac{p(p-1)}{2}$;
or $o(C) = q$ and $c_M (H) = pr = \frac{p(p+1)}{2}$.  In either case
$o(C)c_M (H) = \frac{p (p^2 -1)}{4}$ and the total number of quadratic subgroups
is
$$\frac{p (p^2 - 1)}{24} = \frac{o(M)}{12}.$$

If $o(M)$ is not divisible by eight, the quadratic subgroups are actually Sylow 2-subgroups
of $M$ and hence form a single conjugacy class.  

Now assume that $8|o(M)$.  Then 
$8|o(H)$ and we see that any quadratic subgroup $Q\leqs H$ is normalized by a copy of 
$D_8$ lying in $H$.  We claim that in fact, $N_H (Q) \isom D_8$.  This follows easily
from the Lemma~\ref{dihedral} and the fact that $D_{16}$ does not contain a normal
copy of $D_4$.  Note that this implies that $8|o(N_M (Q))$.

The above argument shows that $Q$ has $\frac{o(H)}{8}$ conjugates in $H$ and since
there are twice this many quadratic subgroups in $H$, the quadratic
subgroups of $H$ fall into two conjugacy classes of equal size.  Since each quadratic
subgroup is contained in a conjugate of $H$, every quadratic subgroup in $M$ is conjugate
to a quadratic subgroup in $H$.  This means that there are at most two conjugacy classes of 
quadratic subgroups in $M$.  

If all the quadratic subrgroups in $M$ were conjugate,
then by Lemma~\ref{norm-conj} the normalizer of any quadratic subgroup would have order
12, contradicting the fact (mentioned above) that $8|o(N_M (Q))$.  Thus there are indeed two
conjugacy classes of quadratic subgroups, as desired.
$\hfill \Box$

\vspace{.15in}

Our next goal will be to determine the normalizer of a quadratic subgroup of $M$.

\begin{theorem}$\label{S_4s-A_4s}$
Let $Q\leqs M$ be a quadratic subgroup.  
If $o(M)$ is divisible by eight then $N_Q (M)\isom S_4$ and otherwise $N_Q (M) \isom A_4$.
\end{theorem}
{\bf Proof.}  Let $N = N_Q (M)$ and let $H$ be a maximal dihedral subgroup of $M$ containing
$Q$.  First we compute $o(N)$.  
If $o(M)$ is not divisible
by eight, Claim~\ref{quadratic-groups} and Lemma~\ref{norm-conj} show that $o(N) = 12$,
and in the same manner we find $o(N) = 24$ if eight divides $o(M)$.

Now, in the case where $o(N) = 12$, consider an element $s\in N$ of order three.  
Let $C\leqs H$ be the cyclic subgroup of index two.
Then by Lemma~\ref{dihedral}, 
$Q\cap C\isom \Z/2$.  Let $a\in Q\cap C$ be the non-trivial element.
If $s$ commutes with $a$, then $s\in N(\gen{a}) = H$,
and if $s$ commuted with another non-trivial element $b\in Q$, $s$ would
lie in the intersection of $H$ and one of its conjugates.  This is impossible
(by Theorem~\ref{cyclic-subgroups}),
since then $s$ would lie in $C$ and one of its conjugates.  
Thus $N$ is a non-abelian group of order
twelve containing a normal quadratic subgroup, and any such group is isomorphic to 
$A_4$ (Dummit and Foote~\cite[p. 170]{Dummit-algebra}).

Next we must show that when $8|o(M)$, $N\isom S_4$.  Since $o(N) = 24$ in this case,
we again have an element $s\in N$ of order three, and the same argument as above
shows that $\gen{s, Q} \isom A_4$.  Now, $A_4$ has four subgroups of order three,
and thus there are at least four subgroups of order three in $N$.  In fact, since
these subgroups are Sylow 3-subgroups, there can be no more than four of them
(the number of Sylow 3-subgroups is equivalent to 1 (mod 3) and must divide 24).

Consider the action of $N$ (by conjugation) on its four subgroups of order three.
This action induces a homomorphism $f:N\to S_4$, and we claim that this map is
an injection.  Since $o(N) = o(S_4)$, this will complete the proof.  

Let $T = \ker(f)$.  
Then $T$ is the intersection of the normalizers of the Sylow 3-subgroups
of $N$, so $T$ is characteristic (and hence normal) in $N$.  No 
(non-trivial) element
of $\gen{s, Q} \isom A_4$ normalizes each of these subgroups, so $T\cap \gen{s, Q} = \{1\}$ 
and hence
$$\gen{s, T, Q} \isom T\cross \gen{s,Q} \isom T \cross A_4\leqs N.$$
Since $o(N) = 24$, we must have $o(T) \leqs 2$.  If $o(T) = 2$, then 
$\gen{Q, T} \isom D_4\cross \Z/2$.  
This group is neither cyclic nor dihedral, contradicting the fact that any 2-group in
$M$ lies in a copy of $D_{2r}$ or $D_{2q}$ and hence is either cyclic or dihedral.
This contradiction shows that $f$ is injective, completing the proof.

$\hfill \Box$

\vspace{.15in}
There is one more type of subgroup which sometimes appears in $M$.

\begin{theorem}$\label{A_5s}$
If $p\equiv \pm 1$ (mod 5) then $M$ contains exactly $\frac{o(M)}{30}$ subgroups isomorphic
to $A_5$, falling into two conjugacy classes of equal size.
\end{theorem}

Since we are mainly interested in the case $p=7$, we refer to 
Burnside~\cite[$\S 324$]{Burnside-theory-of-groups} for the proof.

\vspace{.15in}
The subgroups determined thus far exhaust all the subgroups of $M$.  We will
prove this fact in two steps.

\begin{theorem}  Let $H\leqs M$ be a subgroup whose order is divisible by $p$.
Then either $H = M$ or $H\leqs$Stab$(x)$ for some $x\in \prj$, i.e. 
$H$ lies in the normalizer of a Sylow $p$-subgroup.
\end{theorem}
{\bf Proof.}  Since $o(H)$ is divisible by $p$, there must be a subgroup $P\leqs H$ of order $p$.
This subgroup is in fact a Sylow $p$-subgroup of $H$, since $p^2$ does not
divide the order of $M$.  The number of Sylow $p$-subgroups of $H$ is equivalent to 1
(mod $p$), and hence must be either 1 or $p+1$ (since $M$ contains just $p+1$ subgroups of
order $p$.  Note that the set of Sylow $p$-subgroups of $M$ generates a normal subgroup, and
since $M$ is simple they must generate the whole group.  Thus either $H$ contains a unique
(normal) subgroup of order $p$ or $H = M$.
$\hfill \Box$

\begin{theorem}  If $H\leqs M$ is a subgroup whose order is not divisible by $p$, then one
of the following is true:
\begin{enumerate}
	\item $H$ is cyclic;
	\item $H\isom D_{2q'}$ for some $q'|q$ and $H\leqs$Stab$(x,y)$ for some $x,y\in \prj$;
	\item $H\isom D_{2r'}$ for some $r'|r$ and $H$ is contained in a copy of $D_{2r}$;
	\item $H\isom A_4$;
	\item $H\isom S_4$;
	\item $H\isom A_5$.
\end{enumerate}
\end{theorem}

Note that all possible cyclic subgroups of $M$ were already determined in 
Theorem~\ref{cyclic-subgroups}.  All dihedral subgroups were determined at this point 
as well, since any dihedral subgroup lies in the normalizer of its cyclic subgroup of 
index two.  Thus all that remains to be proved is that any subgroup of $M$ has one of the 
isomorphism types listed above.  (Note that $A_4$ is the only proper subgroup of $S_4$ or $A_5$
which is neither cyclic nor dihedral).

\vspace{.15in}
{\bf Proof.}  Consider a subgroup $H\leqs M$ whose order is not divisible by $p$.
Let $n = o(H)$.
We will break $H$ down into conjugacy classes of cyclic subgroups, beginning with
the largest cyclic subgroup in $H$.

Let $h_1$ be an element of maximum order in $H$, and let $q_1 = o(h_1)$.  By
Theorem~\ref{cyclic-subgroups}, $N_H (\gen{h_1})$ is either cyclic or dihedral,
and by maximality of $q_1$ we see that $o(N_H (\gen{h_1}))$ is either $q_1$ or 
$2q_1$.  Thus the number of conjugates of $\gen{h_1}$ in $H$ is either
$\frac{n}{q_1}$ or $\frac{n}{2q_1}$.  By Theorem~\ref{cyclic-subgroups}, no
two conjugates of $\gen{h_1}$ may intersect non-trivially.  Thus we see that
the total number non-trivial of elements of $H$ lying in conjugates of $\gen{h_1}$ is
$$\frac{n(q_1 - 1)}{\epsilon_1 q_1},$$
where $\epsilon_1 \in\{1,2\}$.

Consider the set 
$$S_1 = \{h\in H: h\notin \gen{h_1}^z \,\, {\rm for\,\, any} \,\, z\in M\}$$
and let $h_2\in S_1$ be an element with maximum order $q_2$.  As before, 
the total number of non-trivial elements of $H$ lying in conjugates of $\gen{h_2}$ is
$$\frac{n(q_2 - 1)}{\epsilon_2 q_2},$$
where $\epsilon_2 \in\{1,2\}$.  We claim that in fact, none of these
elements lie in conjugates of $\gen{h_1}$.  Theorem~\ref{cyclic-subgroups} shows
that if two cylic subgroups of $M$ intersect non-trivially, one must be contained
in the other.  Thus if $\gen{h_2}^y\cap \gen{h_1}^z$ were non-trivial for some
$y,z\in M$, then we would have $\gen{h_2}^y\leqs \gen{h_1}^z$ and 
$\gen{h_2}\leqs \gen{h_1}^{zy^{-1}}$, contradicting the fact that $h_2\in S_1$.

We may continue in this manner, defining $S_2$, $h_3$ and $\epsilon_3$ in the 
obvious ways, and eventually we find $S_m = \emptyset$ for some $m$.  
At this point we have
\begin{equation} \label{counting1}
n = 1 + \sum_{i=1}^m \frac{n(q_i - 1)}{\epsilon_i q_i},
\end{equation}
or
\begin{equation} \label{counting2}
\frac{1}{n} = 1 - \sum_{i=1}^m \frac{q_i - 1}{\epsilon_i q_i}.
\end{equation}
Our goal will be to find all the solutions to this equation, and we will break down
our analysis by considering the possible values of the $\epsilon_i$.

Letting $A = \{i: \epsilon_i = 1\}$ and $B = \{i: \epsilon_i = 2\}$,
we claim that $|A|\leqs 1$.  Indeed, noting that $q_i\geqs 2$ for each $i$,
Equation~\ref{counting2} becomes
$$\frac{1}{n} = 1 - \sum_{i\in A} \frac{q_i - 1}{q_i} - \sum_{j\in B} \frac{q_j - 1}{2q_j}$$
$$ = 1 - \sum_{i\in A} \left(1 - \frac{1}{q_i}\right) 
       - \sum_{j\in B} \left(\frac{1}{2} - \frac{1}{2q_j}\right)$$
$$ = 1 - |A| + \sum_{i\in A} \frac{1}{q_i} - \frac{1}{2} |B| 
             + \sum_{j\in B} \frac{1}{2q_j}$$
$$\leqs 1 - |A| + \frac{1}{2}|A| - \frac{1}{2} |B| + \frac{1}{4}|B| $$
$$= 1 - \frac{1}{2} |A| - \frac{1}{4} |B|.$$
Thus $|A|\leqs 1$ since $\frac{1}{n}>0$, and if $|A| = 1$, then $|B|\leqs 1$.  

We want to show that $m = |A| + |B| \leqs 3$.  This will imply that there
are only several possible cases to consider.  If $m>2$ then (by the above
discussion) we must have 
$|A| = 0$.
If $|A| = 0$, then $\epsilon_i = 2$ for all $i$, and assuming $m>3$ gives
$$\frac{1}{n} = 1 - \frac{q_1 - 1}{2q_1} - \frac{q_2 - 1}{2q_2} - \cdots$$
$$\leqs 1 - \frac{q_1 - 1}{2q_1} - \frac{q_2 - 1}{2q_2} 1 - \frac{q_3 - 1}{2q_3} 
          - \frac{q_4 - 1}{2q_4}$$
$$=1 - \frac{(q_1 -1)q_2q_3q_4 + (q_2 -1)q_1q_3q_4 
             + (q_3 -1)q_1q_2q_4 + (q_4 -1)q_1q_2q_3}
        {2q_1q_2q_3q_4}$$
$$=1 - \frac{q_1q_2q_3 + q_1q_2q_4 + q_1q_3q_4 + q_2q_3q_4 - 4q_1q_2q_3q_4}
	    {2q_1q_2q_3q_4}$$
$$= \frac{1}{2q_4} + \frac{1}{2q_3} + \frac{1}{2q_2} + \frac{1}{2q_1} - 1$$
$$\leqs \frac{1}{4} + \frac{1}{4} + \frac{1}{4} + \frac{1}{4} - 1 = 0,$$
a contradiction.

The various remaining cases can all be dealt with in a manner similar to this
computation.  In each case we either reach a contradiction or show that any group
satisfying the formula must have one of the allowed types.  To illustrate the
methods, we work two of the cases in full.  For the other
cases, we refer to Burnside~\cite[$\S 326$]{Burnside-theory-of-groups}.

First, if $m=1$ and $\epsilon_1 = 1$, Equation~\ref{counting1} gives
$$n = 1 + \frac{n(q_1 - 1)}{q_1}$$
and solving for $n$ we find that $n = q_1$.  Thus in this case $H = \gen{h_1}$,
i.e. $H$ is cyclic of order $q_1$.

The second case we will consider is $m=1$, $\epsilon_1 = 1$, and $\epsilon_2 = 2$.
In this case we have
$$\frac{1}{n} = 1 - \frac{q_1 - 1}{q_1} - \frac{q_2 - 1}{2q_2} $$
$$= \frac{2q_1q_2 - 2q_2(q_1 - 1) - q_1(q_2 - 1)}{2q_1q_2}$$
$$= \frac{2q_2 + q_1 - q_1q_2}{2q_1q_2}.$$
If $q_2\geqs 3$ then
$$\frac{1}{n} \leqs \frac{2q_2 + q_1 - 3q_1}{6q_1}$$
$$ = \frac{1}{6}-\frac{1}{2}+\frac{q_2}{3q_1}$$
$$ =\frac{1}{3} \left(\frac{q_2}{q_1} - 1\right).$$
But $q_2\leqs q_1$ by contruction, so $\frac{q_2}{q_1} - 1\leqs 0$, a contradiction.

So $q_2 = 2$ and we have 
$$\frac{1}{n} = 1 - \frac{q_1 - 1}{q_1} - \frac{2 - 1}{4} 
              = 1 - \frac{q_1 - 1}{q_1} - \frac{1}{4}$$
$$ = \frac{4q_1 - 4(q_1 -1) -q_1}{4q_1} = \frac{4-q_1}{4q_1}.$$
Thus $q_1 < 4$.  If $q_1 = 2$, then $\frac{1}{n} = \frac{4-2}{8} = \frac{1}{4}$ and
$n=4$.  It is easy to see that no group of order four could have these values
for $m$, $\epsilon_1$ and $\epsilon_2$.  So $q_1 = 3$, and we find $n = 12$.
Thus $H$ is a group of order twelve and $h_1\in H$ is an element of order three
with $N_H (\gen{h_1}) = \gen{h_1}$.  This means $\gen{h_1}$ has four conjugates in $H$, and 
$A_4$ is the only group of order twelve containing four subgroups of order three
(Dummit and Foote~\cite[p. 170]{Dummit-algebra}).
$\hfill \Box$ 

\vspace{.15in}
We now record a corollary which will help us determine the M\"{o}bius function
of $PSL_2(\F_7)$.

\begin{corollary}$\label{A_4}$
If the order of $M$ is not divisible by 60, then no two subgroups of $M$ 
intersect in a copy of $A_4$.
\end{corollary}
{\bf Proof.}  Since $60 = o(A_5)$ does not divide $o(M)$, $M$ does not contain subgroups
isomorphic to $A_5$.  Thus we must simply show that no two subgroups $H$ and $K$ isomorphic to 
$S_4$ can intersect in a copy of $A_4$.  But if $H\cap K\isom A_4$, then $H\cap K$ is normal
in both $H$ and $K$.  This is impossible, since $H$ and $K$ are maximal and $M$ is a 
simple group.
$\hfill \Box$

\subsection{The M\"{o}bius Function of $PSL_2(\F_7)$}$\label{mu(psls)-section}$

We will now describe the theory of M\"{o}bius Inversion for finite groups, following Hall's original 
paper~\cite{Hall-eulerian}.  Most of the proofs are omitted, and the interested reader
may refer to Hall's paper.  Once the necessary facts have been stated, we will proceed
to calculate the M\"{o}bius function of $PSL_2 (\F_7)$.

\begin{definition}  Let $G$ be a finite group and let $\hat{S} (G)$ denote the poset of
all subgroups of $G$, ordered by inclusion {\rm(}so $\{1\}, G\in \hat{S} (G)$.{\rm)}
Then the M\"{o}bius function of $G$, $\mu_G : \hat{S} (G) \to \Z$, is given by
$$\mu_G (H) = \left\{ \begin{array}{ll}
                         \tilde{\chi} (\ord{S (G)_{>H}}), & H<G \\
                         1, & H = G,
		      \end{array}
	      \right.
$$
where $\tilde{\chi}$ denotes the reduced Euler characteristic and $S(G)_{>1} = S(G)$.
\end{definition}

The importance of the M\"{o}bius function comes from the following theorem.  

\begin{theorem}[M\"{o}bius Inversion]
$\label{inversion}$
Let $f:\hat{S} (G) \to \Z$ be any function and let 
$g:\hat{S} (G) \to \Z$ be the function
$$g(H) = \sum_{K \leqs H} f(K).$$

Then 
$$f(G) = \sum_{H\in \hat{S} (G)} \mu_G (H) g(H).$$
\end{theorem}

We will use M\"{o}bius inversion to count automorphism classes of generating pairs in
$PSL_2 (\F_7)$.  We now take a moment to explain how this works.

An automorphism class of generating pairs is a set of the form
$$\{(\gamma (\alpha), \gamma (\beta) ): \alpha, \beta\in G, \, \,\, \gamma \in Aut(G)\}.$$
Note that in any automorphism class of generating pairs the orders of the first and
second elements remain constant. 

\begin{definition} Let $G$ be a finite group.  
For $a,b>1$, let $\Phi_{a,b} (G)$ be the number of automorphism classes
of generating pairs $(\alpha, \beta)$ with $o(\alpha) = a$, $o(\beta) = b$, and
let $\phi_{a,b} (G)$ be the number of \emph{ordered}
pairs $(\alpha, \beta)$ {\rm (}$\alpha, \beta\in G${\rm )}
with $o(\alpha) = a$, $o(\beta) = b$ and $\gen{\alpha, \beta} = G$.

Also, let $\sigma_{a,b} (G) = |\{\alpha\in G: o(\alpha) = a\}||\{\beta\in G: o(\beta) = b\}|$.
\end{definition}

Note that for any $a,b>1$ we have
\begin{equation}\label{aut-classes}
\Phi_{a,b} (G) = \frac{\phi_{a,b} (G)}{o(Aut(G))},
\end{equation}
since an automorphism which fixes a set of generators is the identity.
We now have the following important application of Theorem~\ref{inversion}.

\begin{corollary}$\label{counting-classes}$
For any finite group $G$ and any $a,b>1$,
$$\Phi_{a,b}(G) = \frac{1}{o(Aut(G))} \sum_{H\leqs G} \mu_G (H) \sigma_{a,b} (H).$$
\end{corollary}
{\bf Proof.}  This follows immediately from Theorem~\ref{inversion} and 
Equation~\ref{aut-classes} after noting that
$$\sum_{H\leqs G} \phi_{a,b} (H) = \sigma_{a,b} (G).$$
This latter statement is true because each pair $(\alpha, \beta)$ generates a unique
subgroup of $G$.
$\hfill \Box$

\vspace{.15in}
Before calculating the M\"{o}bius function of $PSL_2 (\F_7)$, we note some useful facts
from Hall~\cite{Hall-eulerian}.

\begin{lemma}$\label{mu}$
For any $H < G$, 
$$\mu_G (H) = - \sum_{K>H} \mu_G (K),$$
and (equivalently)
$$\sum_{K\geqs H} \mu_G (K) = 0.$$

If $H$ is a maximal subgroup of $G$, then $\mu_G (H) = -1$.

If $H<G$ is not the intersection of a set of maximal subgroups of $G$, 
then $\mu_G (H) = 0$.  
\end{lemma}

We now turn to the case $G = PSL_2 (\F_7)$.  
Before beginning the computaton of $\mu_G$,
we note that any maximal subgroup of $G$
is isomorphic to either $S_4$ or $\Z/7 \rtimes \Z/3$.  This
is because in this case, $p = 7, q = 3$ and $r = 4$, so 
$D_{2q} = D_6 \isom S_3$ and $D_{2r} = D_8$;  
hence these subgroups sit inside copies of $S_4$ (following 
Burnside~\cite{Burnside-theory-of-groups}, we will
call subgroups isomorphic to $S_4$ \emph{octahedral}).  By Lemma~\ref{mu},
we need only consider subgroups formed by intersecting copies of $\Z/7\rtimes \Z/3$
and of $S_4$.

Note that the only proper, non-trivial subgroups of $\Z/7\rtimes \Z/3$ are the normal
copy of $\Z/7$ and seven copies of $\Z/3$.  Every other proper subgroup of $G$ is a subgroup
of $S_4$, and each of these is isomorphic to one of the following groups:
$S_4$, $A_4$, $D_8$, $D_6 \isom S_3$, $D_4$, $\Z/4$, $\Z/3$, $\Z/2$ or the trivial group.

The following observation will be useful:
if $H, K\leqs G$ and $\exists \gamma \in Aut(G)$ such that $\gamma (H) = K$, 
then $\mu_G (H) = \mu_G (K)$ (here $G$ is any finite group).  This
follows from the fact that $\gamma$ induces a simplicial isomorphism between
$\ord{S(G)_{>H}}$ and $\ord{S(G)_{>K}}$.

Theorem~\ref{cyclic-subgroups} shows that  
if $H,K\leqs G$ are isomorphic cyclic subgroups, then $H$ and $K$ are conjugate. 
The same is true if $H\isom K\isom \Z/7\rtimes \Z/3$, $D_6$, or $D_8$.  
Unfortunately, Theorem~\ref{quadratic-groups} tell us that the quadratic subgroups of $G$ 
fall into
two conjugacy classes, and hence the same is true for octahedral subgroups and
the subgroups isomorphic to $A_4$.  
Luckily, there is an outer automorphism of $G$ interchanging these
conjugacy classes.  In order to explain this automorphism, we need the following facts.

\begin{fact}$\label{gl(3,2)}$
$GL_3 (\F_2) \isom PSL_2 (\F_7)$.
\end{fact}
{\bf Proof.}  The order of $GL_3 (F_2)$ can be calculated in a similar manner
to that of $GL_2 (F_p)$ (as in the proof of Claim~\ref{order}) and we find that
$o(GL_3 (\F_2)) = 168$.  In addition, note that for any $n$,
$GL_n (\F_2) = SL_n (\F_2)$ since $\F_2$ has only one non-zero element.  Also, 
it is not hard to check that the center of $GL_n (\F_2)$ is trivial, so
$GL_3 (\F_2) \isom PSL_3 (\F_2)$, and hence this group is simple 
(by Theorem~\ref{simplicity}).  Up to isomorphism, there is only one simple group 
of order 168, so we must have $GL_3 (\F_2) \isom PSL_2 (\F_p)$.  (For a proof of this
last fact, see Dummit and Foote~\cite[Exercise 27, p. 215]{Dummit-algebra}.)
$\hfill \Box$

\begin{fact} For any ring $R$ and any $n>1$, the map $A\mapsto (A^T)^{-1}$ is an automorphism
of $GL_n (R)$.  (This map is called the transpose-inverse.)
\end{fact}
{\bf Proof.}  First, an easy induction shows that $det (A) = det (A^T)$, and since 
$det (A^{-1}) = - det (A)$, we see that the transpose-inverse is indeed a \emph{function} from
$GL_n (R)$ to $GL_n (R)$.  Next, note that if $(A^T)^{-1} = I$, then $(A^T) = I$ and hence
$A = I$.  So all that remains to be shown is that this map is a homomorphism.

This follows easily from the fact that $(AB)^T = B^T A^T$, as we now have
$((AB)^T)^{-1} = (B^T A^T)^{-1} = (A^T)^{-1} (B^T)^{-1}.$
$\hfill \Box$

\vspace{.15in}
The group $GL_3 (\F_2)$ has a natural action on $\F_2^3 - \{0\}$ via multiplication.
We take a moment to study this action.  

\begin{lemma}$\label{gl3-2-trans}$
The action of $GL_3 (\F_2)$ on $\F_2^3 - \{0\}$ is 2-transitive,
and the stabilizer of any point is isomorphic to $S_4$.
\end{lemma}
{\bf Proof.}  That this action is 2-transitive depends heavily on the fact 
that we are working over $\F_2$, and hence the span of any vector $v\in \F_2^3 - \{0\}$
is just $\{0, v\}$.  Now, thinking of $GL_3 (\F_2)$ as the group of vector-space
automorphisms of $\F_2^3$, for any non-zero vectors $v\neq w$ there is a third
vector $u\in \F_2^3$ s.t. $\{u,v,w\}$ forms a basis for this vector space.  Now
we can find an automorphism of $\F_2^3$ sending $v$ to $v'$ and $w$ to $w'$ for any
$v'\neq w'\in \F_2^3 - \{0\}$.  This shows that the action is 2-transitive, as desired.

To compute the stabilizer of a point, we will think of $GL_3 (\F_2)$ as a group
of matrices acting by multiplication.  Then Stab$(1,0,0)$ is the 
set of all matrices of the form
$$\left[ \begin{array}{lll} 
		1 & a & b \\
		0 & c & d \\
		0 & e & f \\
	\end{array}
  \right],
$$
and it is easy to see that there are 24 such matrices in $GL_3 (\F_2)$.
All subgroups of order 24 in $GL_3 (\F_2) \isom \psls$ octahedral,
so we have Stab$(1,0,0)\isom S_4$ and now the fact that this action is transitive
completes the proof.
$\hfill \Box$

\begin{fact}$\label{oct-aut}$
There is an automorphism of $PSL_2 (\F_7)$ which interchanges the 
two conjugacy classes of octahedral subgroups. 
This automorphism is induced by the transpose-inverse
automorphism of $GL_3 (\F_2)$.
\end{fact}
{\bf Proof.}
Lemma~\ref{gl3-2-trans} shows that one conjugacy class of octahedral subgroups
arrises as stabilizers of points in $\F_2^3 - \{0\}$.
The second conjugacy class is obtained by applying the transpose-inverse.  If 
$A\in$ Stab$(1,0,0)$, then any easy computation shows that $(A^T)^{-1}$ has the form 
$$\left[ \begin{array}{lll} 
		1 & 0 & 0 \\
		a & b & c \\
		d & e & f \\
	\end{array}
  \right],
$$
for some $a,b,c,d,e,f\in \F_2$.  The matrices of this form then form an octahedral
subgroup $H$ (note that there are 24 such matrices and transpose-inverse is
an automorphism) and it remains to show that this subgroup is not the stabilizer of any point 
in $\F_2^3$.  This is easily shown by direct computation.
$\hfill \Box$

\vspace{.15in}
Note that since each quadratic subgroup has an 
octahedral normalizer, the automorphism described above also interchanges the conjugacy classes
of quadratic subgroups.  The same is true for the subgroups isomorphic to $A_4$.

\begin{theorem}$\label{mu(psls)}$
Let $G = PSL_2 (\F_7)$.  If $H<G$ then we have:
$$\mu_G (H) = \left\{ \begin{array}{ll}
			 \mu_G (H) = -1, & H\isom S_4\,\, {\rm or } \,\,\Z/7 \rtimes \Z/3 \\
			 \mu_G (H) = 1, & H\isom D_8 \,\, {\rm or }\,\, D_6 \\
			 \mu_G (H) = 2, & H\isom \Z/3 \\
			 \mu_G (H) = -4, & H\isom \Z/2 \\
			 \mu_G (H) = 0, & otherwise. \\
                      \end{array}
              \right.
$$
\end{theorem}
{\bf Proof.}  
If $H\isom S_4$ or $\Z/7 \rtimes \Z/3$, then $\mu_G (H) = -1$ by Lemma~\ref{mu}.

If $H\isom \Z_7$ then $\mu_G (H) = 0$ because $H$ is not the intersection of two maximal
subgroups of $G$ (copies of $\Z/p\rtimes \Z/q$ are stabilizers of points in $\prj$, and
intersect in copies of $\Z/q$).  

If $H\isom \Z_3$, then $H =$ Stab$(x)\cap$Stab$(y)$ for some $x,y\in \prj$, and hence
$H\leqs$ Stab$(x)$, Stab$(y)$.  Any other subgroups containing a subgroup of order three
is isomorphic to $D_6$, $A_4$ or $S_4$, and thus must contain $N_G (H) \isom D_6$.
By Lemma~\ref{mu} we have 
$$\mu_G (H) = - \left( \mu_G (Stab(x)) + \mu_G (Stab(y)) 
	              + \sum_{K\geqs N_G (H)} \mu_G (K) \right),$$
and applying the lemma again yields $\mu_G (H) = -(-1 - 1 + 0) = 2$.

If $H\isom A_4$ then $\mu_G (H) = 0$ because $H$ is not an intersection of 
maximal subgroups (Corollary~\ref{A_4}).

If $H\isom \Z/4$, then any subgroup containing $H$ also contains $N_G (H)\isom D_8$, 
and hence 
$$\mu_G (H) = - \sum_{K\geqs N_G (H)} \mu_G (K) = 0.$$

In the following calculations, we will exploit a useful counting formula.
Let $C$ and $C'$ be automorphism classes of subgroups of $G$.
Then each element $K\in C'$ contains the same number
of elements $H\in C$ (call this number $\alpha$) and each element $H\in C$
is contained in the same numer of elements of $C'$ (call this number $\beta$).
We have the following formula for $|C|$:
\begin{equation}\label{count}
|C| = \frac{|C'|\alpha}{\beta}.
\end{equation}

If $H\isom D_8$, then any subgroup containing $H$ is octahedral.  
There are $pq = 21$ copies of $D_8$ in $G$, forming
a single conjugacy class, and there are 
$\frac{o(G)}{12} = 14$ copies of $S_4$ forming a single automorphism class.
Thus we have $21 = \frac{14\cdot3}{\beta}$,
and hence $\beta = 2$.  Thus $\mu_G (H) = - (-1 - 1 + \mu_G (G)) = 1$.
For $H\isom D_6, D_4$ or $\Z/2$, the computation follows similarly.

Finally, Lemma~\ref{mu} can be used to show $\mu_G (\{1\}) = 0$.
$\hfill \Box$

\vspace{.15in}
Now that we know the M\"{o}bius function of $G$, we can give an explicit form 
of Corollary~\ref{counting-classes}.  Clearly $\sigma_{a,b} (H)$ depends only on the
isomorphism type of $H$, and we have seen that the same is true of $\mu_G (H)$.
By counting the number of subgroups of each isomorphism type we obtain the
following inversion formula (for any $a,b>1$):
$$\phi_{a,b} (G) = \sigma_{a,b} (G) - 14\sigma_{a,b} (S_4) - 8 \sigma_{a,b} (\Z/7 \rtimes \Z/3)
                    + 21 \sigma_{a,b} (D_8)$$
\begin{equation}\label{count-classes2}
   + 28 \sigma_{a,b} (D_6) + 56 \sigma_{a,b} (\Z/3) - 84 \sigma_{a,b} (\Z/2).
\end{equation}

In order to count automorphism classes, we need to know the order of $Aut(G)$.
This can be calculated using M\"{o}bius inversion.  

\begin{fact}$\label{Aut(G)}$
$o(Aut(\psls)) = 336$.
\end{fact}
{\bf Proof.}
First, note that the 
group $GL_2 (\F_7)$ acts by conjugation on $SL_2 (\F_7)$ and induces 
automorphisms of $G$ (note that $I$ and $-I$ are fixed by conjugation).  
Thus we have a map $f: GL_2 (\F_7)\to Aut(G)$.  We claim
that the kernel of this map is $\{\lambda I: \lambda \in \F_p^*\}$.  Certainly
these elements act trivially, and a simple matrix computation shows that
any element in $\ker (f)$ has this form.  This shows that 
$$o(Aut (G)) \geqs \frac{o(GL_2 (\F_7))}{o(\ker (f))} = \frac{(7^2-1)(7^2 - 7)}{6} = 336.$$
To prove that there are no other elements in $Aut(G)$, we can 
Equation~\ref{count-classes2} and the observation (made earlier) that 
$o(Aut (G))$ divides $\phi_{a,b}$ for any $a,b>1$.  Now, consider $\phi_{2,3}$.
The above formula gives
$$\phi_{2,3} (G) = \sigma_{2,3} (G) - 14\sigma_{2,3} (S_4) + 28 \sigma_{2,3} (D_6)$$
$$               = 21\cdot 56 - 14\cdot 9\cdot 8 + 28 \cdot 3 \cdot 2 = 336,$$
so $o(Aut (G)) = 336$ as desired.
$\hfill \Box$

\begin{corollary}$\label{count-aut-classes}$
For any $a,b>1$, 
$$\Phi_{a,b} (G) = \frac{\phi_{a,b} (G)}{336},$$
where $\Phi_{a,b}$ denotes the number of automorphism classes of generating
pairs with orders $a$ and $b$.
\end{corollary}

\subsection{Simple Connectivity of $\C{\psls}$}$\label{pi_1-psls}$

The proof of simple connectivity will be analagous to the second proof for $\C{A_5}$, 
and will use 2-transitive actions.
In order to minimize the amount of computation in the proof, 
we will now prove two helpful lemmas.

\begin{lemma}$\label{psls-orders}$
Let $\alpha = \frac{ax+b}{cx+d}$ be any non-trivial element of $\psls$.  
We define the trace-squared of $\alpha$ to be $tr^2 (\alpha) = (a+d)^2$ 
(note that this is the square of the
trace of either representative of $\alpha$ in $SL_2 (\F_7)$, and thus is well-defined).

The order of $\alpha$ is determined as follows:
$$o(\alpha) = \left\{ \begin{array}{ll}
			2, & tr^2 (\alpha) = 0 \\
			3, & tr^2 (\alpha) = 1 \\
			4, & tr^2 (\alpha) = 2 \\
			7, & tr^2 (\alpha) = 4 \\
		\end{array}
	\right.
$$
\end{lemma}
{\bf Proof.}  In order to prove this, it will be useful to define the discriminant
of a non-trivial element of $G = \psls$.  If $\alpha = \frac{ax+b}{cx+d}\in G$, we define
$disc(\alpha)$ to be the discriminant of the quadratic polynomial 
$(cx+d)x-(ax+b) = cx^2 + (d-a)x -b$ determined
by the equation $\frac{ax+b}{cx+d} = x$,
so that $disc(\alpha) = (d-a)^2 - 4 (-b)(c) = tr^2(\alpha) - 4$.

We will consider only the elements in $G - $Stab$(\infty)$.  It is easy to check
the result on remaining elements of $G$.
The elements of order seven in $G - $Stab$(\infty)$ are exactly those with one fixed point
in $\F_7$, i.e.
those with $disc (\alpha) = 0$ and $tr^2 (\alpha) = 4$. Similarly, elements of order
three in $G - $Stab$(\infty)$ are those with two fixed points in $\F_7$, i.e. 
those with $disc(\alpha) = z^2$ for
some $z\in \F_7^*$.  Thus $o(\alpha) = 3 \iff disc(\alpha)\in \{1,2,4\} 
\iff tr^2 (\alpha)\in \{5, 6, 1\} \iff tr^2 (\alpha) = 1$.

Next, let $A\in SL_2 (\F_7)$ be a matrix representing $\alpha\in G$, and
let $\lambda_1, \lambda_2\in \bar{\F_7}$ be the eigenvalues of $A$
(where $\bar{\F_7}$ denotes the algebraic closure).  Then $o(\alpha) = 2
\implies \lambda_1^2 = \lambda_2^2 = \pm 1$.  If $\lambda_1 = \lambda_2$ then
$A = \lambda_1 I$ and $\alpha = 1$, so $o(\alpha) = 2 \iff \lambda_1 = - \lambda_2
\iff tr(A) = \lambda_1 + \lambda_2 = 0 \iff tr^2 (\alpha) = 0$.

Finally, a similar but more lengthy calculation shows that 
$o(\alpha) = 4 \iff tr^2 (\alpha) = 2$.
$\hfill \Box$

\begin{lemma}$\label{2-3-generators}$
Let $g,h\in \psls$ be elements of orders two and three, respectively.  Then 
$\gen{g,h} = \psls \iff o(gh) = 7$.
\end{lemma}  
{\bf Proof.}
Certainly, if $o(g) = 2$, $o(h) = 3$ and $o(gh) = 7$, $\gen{g,h} = \psls$,
since no subgroup contains elements of orders two and seven.

In the other direction, recall that there is a unique automorphism class of generators
with orders two and three (see the proof of Fact~\ref{Aut(G)}), and thus the result will follow
if we show that there exist elements  $g,h\in \psls$ with $o(g) = 2$, $o(h) = 3$ and $o(gh) = 7$.
The elements $g = \frac{x-2}{x-1}$ and $h = \frac{4x}{2}$ satisfy this property, since
$tr^2 (g) = 0$, $tr^2 (h) = 1$ (note that $det(h) = det (g) = 1$) and
$hg = \frac{4 \frac{x-2}{x-1}}{2} 
          = \frac{4x - 1}{2x - 2}$ 
so $tr^2 (hg) = 4$ and
$o(hg) = 7$.
$\hfill \Box$

\begin{theorem}$\label{psls-simply-connected}$
The coset poset of $PSL_2 (\F_7)$ is simply connected.
\end{theorem}
{\bf Proof.}  We will show that the minimal cover of $G = \psls$ is simply connected.  
We will employ Theorem~\ref{presentation}, using for a maximal
trees the set of all edges $\{1,h\} \in \M{G}$ ($h\in G$).  First we will show that all
generators of $\pi_1 (\M{G})$ corresponding to edges $\{g, h\}$ with $o(g) = 2$ are trivial, 
and then we will finish the proof by considering 2-transitive actions (note that, since
$G$ is non-cyclic, there is an edge between any two vertices in $\M{G}$).  As in 
Theorem~\ref{presentation}, we will denote 
the generator corresponding to the edge $\{g,h\}$ by $(g,h)$.  [Recall that in the
statement of the theorem, there was a generator for each \emph{ordered} edge, but
since we have the relation $(u,v) = (v,u)^{-1}$ for each edge $\{u,v\}$, it suffices
to show that in each case either $(u,v)$ or $(v,u)$ is trivial.]

We now analyze the various edges $\{g, h\}$ with $o(g) = 2$ appearing in $\M{G}$.
First, note that if $\gen{g,h}\neq PSL_2 (\F_7)$, then the edge $(g,h)$ is trivial
because the set $\{1,g,h\}$ is contained in a subgroup.  Also, since
$G$ is simple two elements of order two cannot generate $G$ and
thus if $o(g) = o(h) = 2$, $(g,h)$ is trivial.  We now examine the 
various possibilities for $o(h)$, and in each case we count the number of automorphism
classes of generators (using Equation~\ref{count-classes2})
and examine an edge from each class.

\vspace{.1in}
{\bf o(g) = 2, o(h) = 3:} Letting $h = \frac{4x}{2}$ and $g = \frac{x-2}{x-1}$, the
proof of Lemma~\ref{2-3-generators} shows that $(g,h)$ is a representative 
for the unique automorphism class of generators with these orders.  

Now, notice that $g(-1) = h(-1) = -2$, so $g\equiv h$ (mod Stab$(-1)$).
Consider the element $z = \frac{3x-2}{-2x-3}$.
This element has determinant one and order two, and we have 
$z (-1) = -2$.  Thus $\{g,h,z\}$ is a 2-simplex
in $\M{G}$, and we have the relation $(g,h) = (h,z)(g,z)$.  
Note that $(g,z) = 1$ because
these two elements of order two cannot generate $G$.  Thus to prove that $(g,h) = 1$,
it suffices to show that $(h,z) = 1$.  We have 
$hz = \frac{4\frac{3x-2}{-2x-3}}{2} = \frac{5x - 1}{3x + 1}$, and hence
$o(hz) = 3$.  Lemma~\ref{2-3-generators} now implies that $\gen{z,h}\neq G$.

\vspace{.1in}
{\bf o(g) = 2, o(h) = 4:} Again, M\"{o}bius inversion shows that $\Phi_{2,4} (G) = 1$:
$$\phi_{2,4} (G) = \sigma_{2,4} (G) - 14 \sigma_{2,4} (S_4) + 21 \sigma_{2,4} (D_8)$$
$$	     = 21 \cdot 42 - 14 \cdot 9 \cdot 6 + 21 \cdot 5 \cdot 2 = 336,$$
and hence $\Phi_{2,4} = 1$.  We claim that the pair $g = \frac{-1}{x}$, 
$h = \frac{4x+1}{-x}$ is a representative for the unique automorphism class of generators.
First, note that these elements each have determinant one and have the correct orders.  To show
that $\gen{g,h} = G$, note that $gh = \frac{-1}{\frac{4x+1}{-x}} = \frac{x}{4x+1}$,
so $o(gh) = 7$.  Since no proper subgroup of $G$ contains elements of orders two and seven,
$\gen{g,h} = G$.

Next, note that $g(0) = h(0) = \infty$, and set $z = \frac{-2}{4x}$.  Then $det (z) = 1$,
$o(z) = 2$, and $z(0) = \infty$.  As in the previous case, if $\gen{z, h}\neq G$, then
$\{g,h\} = 0$.  The proof that $\gen{z, h}\neq G$ is also analagous to the previous case:
$zh = \frac{-2}{4\frac{4x+1}{-x}} = \frac{2x}{2x+4}$ and $o(zh) = 3\neq o(gh)$.  Thus
the pair $(z,h)$ does not fall into the unique automorphism class of generators.

\vspace{.1in}
{\bf o(g) = 2, o(h) = 7:} This time we find that there are three automorphism
classes of generating pairs, since any pair $(g,h)$ with $o(g) = 2$ and $o(h) = 7$ 
necessarily generates.  This means $\phi_{2,a} (G) = 21 \cdot 48 = 1008$ and 
$\Phi_{2,7} = \frac{1008}{336} = 3$.

Let $h = x+1$, so $o(h) = 7$, and let $g = \frac{b}{cx}$, so that $o(g) = 2$.
We have $hg = \frac{b}{cx} + 1 = \frac{cx+b}{cx}$
and thus $tr^2 (hg) = c^2$.  This implies that the pairs
$$\begin{array}{lll} (g_1 = \frac{-1}{x},\, h), &  (g_2 = \frac{2}{3x},\, h), 
                          & (g_3 = \frac{3}{2x},\, h) \\
  \end{array}
$$
represent the three generating automorphism classes. 

Next, say there exists an element $z\in$ Stab$(\infty)$ such that $o(z) = 3$ and
$\{h, g_i, z\}$ forms a simplex in $\M{G}$.  Then we have $(h,z) = 1$ because
$\gen{h,z}\leqs$ Stab$(\infty)$ and $(g,z) = 1$ because $o(g_i) = 2$ and $o(z) = 3$ 
and we have shown above that the generators corresponding to such edges are trivial.
We will now attempt to find such elements.

Consider the equations $h(x) = g_i (x)$.  These correspond to the equations
$x+1 = \frac{b}{-b^{-1} x}$, where $b = -1, 2$ or $3$.  Equivalently, 
[note that $\infty$ can never be a solution] we are
looking for solutions to the quadratic polynomial $x^2 + x + b^2 = 0$, 
and its discriminant is 
$$1-4b^2 = \left\{ \begin{array}{ll}
			1-4(-1^2) = 4, & b=-1 \\
			1-4\cdot 2^2 = -1, & b=2 \\
 			1-4\cdot 3^2 = 0, & b = 3. \\
		   \end{array}
	    \right.
$$
Thus solutions exist when $b = -1$ or $3$, but not when $b = 2$.

For the pairs $(g_1,h)$ and $(g_3,h)$ we may now quickly finish the proof by noting that
for any $x,y\in \F_7$ there is an element $z\in$ Stab$(\infty)$ with $o(z) = 3$
and $z(x) = y$, as can be shown by direct computation (it suffices to consider
elements of the form $\frac{2x + b}{4}$).
Taking into account the previous two paragraphs, this shows that
$(g_1,h) = (g_3,h) = 1$.  

Finally, we must show that $(g_2, h) = (\frac{2}{3x}, x+1) = 1$.  
Letting $z = \frac{2x}{4}$, we
have $o(z) = 3$ and $z\in$ Stab$(\infty)$, so it suffices to show that these three
elements lie in a proper coset, i.e. that $\gen{z^{-1} g_2, z^{-1} h}\neq G$.
We have $z^{-1} = \frac{4x}{2}$, so 
$z^{-1} g_2 = \frac{4\frac{2}{3x}}{2x} = \frac{1}{-x}$ and hence 
$o(z^{-1} g_2) = 2$.  Next, $z^{-1} h = \frac{4(x+1)}{2}$ so $o(z^{-1} h) = 3$.
Also, $z^{-1} g_2 z^{-1} h = \frac{1}{-\frac{4x+4}{2}} = \frac{2}{3x-3}$, so 
$o(z^{-1} g_2 z^{-1} h) = 4$ and hence $\gen{z^{-1} g_2, z^{-1} h}\neq G$.

\vspace{.1in}
We have now shown that all edges $\{g,h\}$ with $o(g) = 2$ correspond to trivial generators.
We will complete the proof of simple connectivity by considering 2-transitive actions of
$G$.  Note that since $G$ is not cyclic, any two
elements of $G$ are equivalent modulo some maximal subgroup.  Thus it will suffice to 
show that there is an element of order two in each coset $kM$ where $k\in G$ and 
$M\leqs G$ is maximal.  We will do this by showing that in each case the action of
$G$ on the cosets $G/M$ is two-transitive and that some element of order two acts 
non-trivially.

First, recall that each subgroup isomophic to $\Z/7\rtimes \Z/3$ is the stabilizer of 
a point in $\F_7 \cup \{\infty\}$, and this action is two-transitive
(Lemma~\ref{2-transitive}).  Also, since no element of $G$ other than the identity acts 
trivially on $\F_7 \cup \{\infty\}$, there are certainly elements of order two which act
non-trivially.  

Next, Lemma~\ref{gl3-2-trans} shows that the action of $G$ on one of its conjugacy classes
of octahedral subgroups is 2-transitive, and Fact~\ref{oct-aut} shows that
there is an automorphism interchanging the two conjugacy classes, so in fact $G$
acts 2-transitively on each conjugacy class of octahedral subgroups.  The stabilizer of a 
subgroup is just its normalizer, and the octahedral subgroups are self-normalizing 
(because they are maximal and $G$ is simple).  
Finally, we need to check that there are elements of order two which act
non-trivially.  This follows immediately from the fact that in $GL_3 (\F_2)$, only
the identity acts trivially on $\F_2^3 - \{0\}$ (and this action has octahedral stabilizers).
$\hfill \Box$

\vspace{.15in}

There are essentially two barriers preventing the extension of this result to 
$\psl$ for primes $p>7$.  First, it is not clear how to generalize the
ad hoc portion of the proof, in which we showed that all edges $\{g,h\}$ with
$o(g) = 2$ corresponded to trivial generators of $\pi_1 (\M{G})$.  Also, for large $p$ 
the only 2-transitive action of $\psl$ is its standard action on $\prj$.  This follows
from~\cite[Exercises 38 and 39, p. 58]{Lang-algebra}, 
and the essential reason is that the subgroups
of $\psl$ are too small, and thus have large index.  (Whenever a group $G$ acts
2-transitively on a set $X$, we get a transitive action of $G$ on the set
$X\cross X - \{(x,x):x\in X\}$, and thus this set can be no larger than 
$G$.  It is possible that there are other small primes for which one gets 
multiple 2-transitive actions.)

\vspace{.15in}
There are a few more facts about the homotopy-type of $\C{\psls}$ that
may be proven by elementary methods.
It particular, it is possible to show that $\C{\psls}$ is homotopy equivalent to a 
three-dimensional complex.  We will briefly explain the process.

First, note that any chain of length five lies under an octahedral subgroup, since the
only other type of maximal subgroup is $\Z/7\rtimes \Z/3$.  It is easy to check that
every chain of length five contains a coset $xH$ where $H\isom D_4$ or $H\isom \Z/4$.
But cosets of copies of $D_4$ and $\Z/4$ may be removed from $\C{\psls}$
without changing the homotopy type, as can be shown by applying Quillen's Theorem
to the inclusion map (see the proof of Claim~\ref{A_5-dimension}).  This shows
that $\C{\psls}$ is homotopy equivalent to a three-dimensional complex.

We can now show that $H_2 (\C{\psls})$ is non-trivial, and in fact has rank at 
least $17\cdot 168$.  This follows immediately from a computation of the Euler characteristic
of $\C{\psls}$.  It is shown in Brown~\cite[Table I]{Brown-coset-poset}
that $\tilde{\chi} (\C{\psls}) = 17\cdot 168$, and
since the only even dimension in which $\C{\psls}$ has homology is dimension two, the
result follows.

It is possible to eliminate various other simplices from $\ord{\C{G}}$ without 
changing the homotopy type.  For
example, all cosets of all subgroups isomorphic to $A_4$ and $\Z/7$ may be removed
(by applying Quillen's Theorem to the inclusion map),
and after removing the cosets $xA_4$ we may remove all simplices of the form 
$\{x\}\subset x C_3 \subset x D_6 \subset x S_4$
or $\{x\}\subset x C_3 \subset x S_4$ (where $C_3$ denotes a cyclic subgroup of order three). 
The latter type of chain is contained in a unique 
maximal chain (obtained by adding the coset $x N(C_3)$) and thus removing these simplices 
constitutes an ``elementary colapse'' (see Brown~\cite{Brown-rewriting}).  Unfortunately,
it does not seem possible to elimanate all chains of length four, and thus we are
unable to prove that $\C{\psls}$ is homotopy equivalent to a wedge of two-spheres.  In
fact, we conjecture that $H_3 (\C{\psls})$ is non-trivial (see Chapter~\ref{cp-sp}).

\chapter{Conjectures and Directions for Further Research}$\label{cp-sp}$

In this final chapter, we briefly discuss some conjectures about the coset poset
and some potential directions further research could take.  Other such remarks
have been scattered through the previous chapters, but these seemed not to fit in elsewhere.

There appear to be certain concrete connections between the subgroup poset and the
coset poset of a finite group $G$; in some sense, the structure of the coset poset
depends heavily on the structure of the subgroup poset.

There is a general theme in the proofs that the coset posets $A_5$ and $\psls$ are
homotopy equivalent to lower dimensional subposets: 
we may remove cosets from $\C{G}$ without changing
the homotopy type whenever the corresponding subgroups can be removed from $S(G)$
(without changing the homotopy type).
Initially, the coset poset has dimension one more than the dimension of $S(G)$ (due
to the fact that singleton cosets are included in $\C{G}$), and as we remove cosets
in this fashion, the relationship between dimensions remains intact.
This suggests the idea that the coset poset should have ``homotopy dimension'' 
exactly one more than the
``homotopy dimension'' of the subgroup poset, where by homotopy dimension we mean the
lowest dimension of a simplicial complex with the same homotopy type.  (In some cases,
the subgroup poset is contractible and this relationship does not hold.)  For solvable
groups with non-contractible subgroup posets, 
this relationship follows immediately from Brown's determination of the homotopy type of $\C{G}$
(Theorem~\ref{solvable-coset-poset}) and the following similar result of
Kratzer and Th\'{e}venaz~\cite[Corollaire 4.10]{Kratzer-subgroups}:

\setcounter{fact}{0}

\begin{theorem}
If $G$ is a finite solvable group, then $S(G)$ is homotopy equivalent
to a wedge of spheres of dimension $d-2$, where $d$ is the length of
a chief series for $G$ (in the sense of Theorem~\ref{solvable-coset-poset}).  The number of
spheres is zero if and only if some subgroup in a chief series for $G$ does not have a complement.
\end{theorem}

It is worth mentioning at this point that the coset poset can be viewed as a homotopy colimit,
in the sense of~\cite{Welker-hom-colim}.  The underlying category is the poset $S(G) \cup \{1\}$,
and the space above a subgroup $H$ is a discrete set of $(G:H)$ points (thought of as
the cosets of $H$ in $G$).  The map associated to the inclusion $H<K$ is just $xH\mapsto xK$,
and the fact that this diagram has homotopy colimit homeomorphic to $\C{G}$ is an immediate
consequence of~\cite[Proposition 4.1]{Welker-hom-colim}.  Originally I had hoped that, 
given a poset
$P\heq S(G)$, one could find a diagram over $C(P)$ (the cone on $P$, i.e. $P$ with an extra
element $0$ satisfying $0<p$ for each $p\in P$) with homotopy colimit $H\heq \C{G}$.  
The desired dimension result would then follow immediately 
from~\cite[Proposition 4.1]{Welker-hom-colim}.  
Unfortunately, the literature on homotopy colimits 
does not seem to contain results of this nature.

\vspace{.15in}
In Theorem~\ref{homology-surjection}, we saw that in certain cases there is a surjection
from $H_{n+1} (\C{G})$ to $H_n (S(G))$, and I wonder whether the hypothesis (that 
$\C{G}_g$ has trivial $n$-dimensional homology) is truly necessary.  
It is almost certain that the hypothesis is not always satisfied (for example I doubt
that $\psls$ satisfies the hypothesis when $n = 2$), 
but the theorem may remain true for other reasons.  In particular, one ought to be able to
settle this question for solvable groups, due to the two
theorems discussed above.  

\vspace{.15in}

When applied to the coset poset of $G = \psls$, these conjectures 
leads me to believe that $H_3 (\C{G})$
is probably non-zero, because $H_2 (S(G))$ is non-zero.  The latter fact follows
from two previous observations and a result of Shareshian: from Theorem~\ref{mu(psls)}
we know that the reduced Euler characteristic of $S(G)$ is zero, and 
Shareshian~\cite[Lemma 3.11]{Shareshian-shell-solv}
has shown that $H_1 (S(G))\neq 0$ (this follows easily by considering the graph formed
by subgroups of type $\Z/3$ and $\Z/7\rtimes \Z/3$ and noting that the edges of this graph
are maximal simplices in $S(G)$).  The proof (sketched at the end of Section~\ref{pi_1-psls})
that $\C{G}$ is homotopy equivalent to a three-dimensional complex  also shows
that $S(\psls)$ is homotopy equivalent to a two-dimensional complex, and thus $H_2 (S(G))$
must be non-zero in order for $\tilde{\chi} (S(G))$ to be zero.

Even if $H_3 (\C{G})$ is zero, our first conjecture about dimensions suggests that
$\C{G}$ is not a bouquet of two-spheres, since $S(G)$ is not homotopy equivalent to 
a one-dimensional complex.  This would yield the only known example of a group whose
coset poset is not a bouquet of equal-dimensional spheres.  (I do not know whether
there are other examples of groups $G$ with $S(G)$ not a bouquet of spheres; it is possible
that $\psls$ is the only such group known.)

\vspace{.15in}
On a different note, it would be interesting to study the coset poset of an infinite group.
Many of the results in this thesis probably extend (in some manner) to this more general
setting.
For example, the Nerve Theorem holds for infinite simplicial complexes (under appropriate
conditions) and this should allow one to define the minimal cover of an infinite group.

\bibliography{references}

\bibliographystyle{plain}

\end{document}